\documentclass[11pt,titlepage]{article}

\usepackage{latexsym}
\usepackage{amsmath}
\usepackage{amssymb}
\usepackage{amscd}
\usepackage{amsthm}

\theoremstyle{plain}
 \newtheorem{theorem}{Theorem}[subsection]
 \newtheorem{proposition}[theorem]{Proposition}
 \newtheorem{lemma}[theorem]{Lemma}
 \newtheorem{corollary}[theorem]{Corollary}
 \newtheorem*{theo3.3.3}{Theorem 3.3.3}
 \newtheorem*{theo3.2.7}{Theorem 3.2.7}
 \newtheorem*{theo3.4.1}{Theorem 3.4.1}
 \newtheorem*{theo4.2.5}{Theorem 4.2.5}

\theoremstyle{definition}
 \newtheorem{definition}[theorem]{Definition}
 \newtheorem{example}[theorem]{Example}
 \newtheorem{conjecture}[theorem]{Conjecture}

\theoremstyle{remark}
 \newtheorem{remark}[theorem]{Remark}

\makeindex

\begin{document}

\begin{titlepage}
\begin{center}
\LARGE\bfseries
Derived categories of coherent sheaves on rational homogeneous manifolds\\
\mdseries\normalsize
\vspace{6cm}
Christian B\"ohning\\
Mathematisches Institut der Universit\"at Bayreuth\\
Universit\"atsstra\ss e, D-95440 Bayreuth, Germany\\
e-mail address: \ttfamily{boehning@btm8x5.mat.uni-bayreuth.de}
\end{center}
\end{titlepage}
\bfseries\small
\emph{Abstract.}\mdseries \: One way to reformulate the celebrated
theorem of Beilinson is that $(\mathcal{O}(-n),\dots , \mathcal{O})$
and $(\Omega^n(n), \dots , \Omega^1 (1), \mathcal{O})$ are strong
complete exceptional sequences in $D^b(Coh\,\mathbb{P}^n)$, the
bounded derived category of coherent sheaves on $\mathbb{P}^n$. In a
series of papers (\cite{Ka1}, \cite{Ka2}, \cite{Ka3}) M. M. Kapranov
generalized this result to flag manifolds of type $A_n$ and
quadrics. In another direction, Y. Kawamata has recently proven
existence of complete exceptional sequences on toric varieties
(\cite{Kaw}).\\
Starting point of the present work is a conjecture of F. Catanese
which says that on every rational homogeneous manifold $X=G/P$, where
$G$ is a connected complex semisimple Lie group and $P\subset G$ a
parabolic subgroup, there should exist a complete strong exceptional
poset (cf. def. 2.1.7 (B)) and a bijection of the elements of the
poset with the Schubert varieties in $X$ such that the partial order
on the poset is the order induced by the Bruhat-Chevalley order
(cf. conjecture 2.2.1 (A)).
 An answer to this question would also be
of interest with regard to a conjecture of B. Dubrovin (\cite{Du},
conj. 4.2.2) which has its source in considerations concerning a
hypothetical mirror partner of a projective variety $Y$: There is a
complete exceptional sequence in $D^b(Coh\, Y)$ if and only if the
quantum cohomology of $Y$ is generically semisimple (the complete form
of the conjecture also makes a prediction about the Gram matrix of such a
collection). A proof of this conjecture would also support 
M. Kontsevich's homological mirror conjecture, one of the most
important open problems in applications of complex geometry to physics today (cf. \cite{Kon}).\\
The goal of this work will be to provide further evidence for 
F. Catanese's conjecture, to clarify some aspects of it and to supply
new techniques. In section 2 it is shown among other things
that the length of every complete exceptional sequence on $X$ must
be the number of Schubert varieties in $X$ and that one can find a
complete exceptional sequence on the product of two varieties once
one knows such sequences on the single factors, both of which follow
from known methods developed by Rudakov, Gorodentsev, Bondal et
al. Thus one reduces the problem to the case $X=G/P$ with $G$
simple. Furthermore it is shown that the conjecture holds true for the
sequences given by Kapranov for Grassmannians and quadrics. One
computes the matrix of the bilinear form on the Grothendieck $K$-group
$K_{\circ}(X)$ given by the Euler characteristic with respect to the
basis formed by the classes of structure sheaves of Schubert varieties
in $X$; this matrix is conjugate to the Gram matrix of a complete
exceptional sequence. Section 3 contains a proof of theorem 3.2.7
which gives complete exceptional sequences on quadric bundles over
base manifolds on which such sequences are known. This enlarges substantially the
class of varieties (in particular rational homogeneous manifolds) on 
which those sequences are known to exist. In the
remainder of section 3 we consider varieties of isotropic flags in a
symplectic resp. orthogonal vector space. By a theorem due to Orlov
(thm. 3.1.5) one reduces the problem of finding complete exceptional
sequences on them to the case of isotropic Grassmannians. For 
these, theorem 3.3.3 gives generators of the derived category which
are homogeneous vector bundles; in special cases those can be used to
construct complete exceptional collections. In subsection 3.4 it is
shown how one can extend the preceding method to the orthogonal case
with the help of theorem 3.2.7. In particular we prove theorem 3.4.1
which gives a generating set for the derived category of coherent
sheaves on the Grassmannian of isotropic 3-planes in a 7-dimensional
orthogonal vector space. Section 4 is dedicated to providing the
geometric motivation of Catanese's conjecture and it contains an alternative
approach to the construction of complete exceptional sequences on
rational homogeneous manifolds which is based on a theorem of M. Brion
(thm. 4.1.1) and cellular resolutions of monomial ideals \`{a} la
Bayer/Sturmfels. We give a new proof of the theorem of Beilinson on
$\mathbb{P}^n$ in order to show that this approach might work in
general. We also prove theorem 4.2.5 which gives a concrete
description of certain functors that have to be investigated in this
approach.        
\vspace{0.5cm}\\
\bfseries\small
\emph{Zusammenfassung.}\mdseries \: Eine Art, den ber\"uhmten Satz von
Beilinson (\cite{Bei}) zu formulieren, ist die folgende: $(\mathcal{O}(-n),\dots , \mathcal{O})$
und $(\Omega^n(n), \dots , \Omega^1 (1), \mathcal{O})$ stellen 
vollst\"andige (starke) exzeptionelle Folgen in
$D^b(Coh\,\mathbb{P}^n)$, der beschr\"ankten derivierten Kategorie
koh\"arenter Garben auf dem $\mathbb{P}^n$, dar. M. M. Kapranov
verallgemeinerte dieses Ergebnis in einer Reihe von Arbeiten
(\cite{Ka1}, \cite{Ka2}, \cite{Ka3}) auf Fahnenmannigfaltigkeiten vom
Typ $A_n$ und Quadriken. In einer anderen Richtung hat Y. Kawamata
k\"urzlich die Existenz vollst\"andiger exzeptioneller Folgen f\"ur
torische Variet\"aten bewiesen (\cite{Kaw}).\\
Ausgangspunkt der vorliegenden Arbeit ist eine Vermutung von
F. Catanese, die besagt, da\ss\ auf jeder rational-homogenen
Mannigfaltigkeit $X=G/P$, wobei $G$ eine zusammenh\"angende
halbeinfache komplexe Liegruppe und $P\subset G$ eine paraboli-\\sche
Untergruppe bezeichnet, eine vollst\"andige starke exzeptionelle
partiell geordnete Menge (vgl. Def. 2.1.7 (B)) und
eine Bijektion zwischen den Elementen dieser Menge und den
Schubertvariet\"aten in $X$ existieren sollte, so da\ss\ die partielle
Ordnung gerade die von der Bruhat-Chevalley Ordnung induzierte ist (vgl. Vermutung
2.2.1 (A)). Eine Antwort hierauf ist auch von Interesse in Hinsicht auf
eine Vermutung von B. Dubrovin (\cite{Du}, conj. 4.2.2), deren
Motivation aus Betrachtung eines hypothetischen Spiegelpartners einer
projektiven Variet\"at $Y$ entspringt: Es gibt eine vollst\"andige
exzeptionelle Folge in $D^b(Coh\, Y)$ dann und nur dann, wenn die
Quantenkohomologie von $Y$ generisch halbeinfach ist (die
vollst\"andige Form der Vermutung macht auch eine Aussage \"uber die
Form der Grammatrix einer solchen Folge). Ein Beweis dieser Vermutung
w\"urde auch weiteren R\"uckhalt f\"ur die Richtigkeit von
M. Kontsevichs homologischer Spiegelvermutung geben, die heutzutage
eines der wichtigsten offenen Probleme in den Anwendungen der
komplexen Geometrie auf die Physik darstellt (\cite{Kon}).\\
Ziel dieser Arbeit soll es sein, weitere Belege f\"ur die 
oben genannte Vermutung von F. Catanese zu liefern, einige ihrer
Aspekte zu kl\"aren und neue Techniken bereitzustellen. In Abschnitt 2
wird unter anderem gezeigt, da\ss\ die L\"ange jeder vollst\"andigen
exzeptionellen Folge auf $X$ die Anzahl der Schubertvariet\"aten in
$X$ sein mu\ss\ und da\ss\ man auf dem Produkt zweier Variet\"aten
eine vollst\"andige exzeptionelle Folge angeben kann, sobald man
solche Folgen f\"ur die Faktoren kennt, was beides unmittelbar aus
bekannten Methoden von Rudakov, Gorodentsev, Bondal u.a. folgt. Damit
reduziert man das Problem auf $X=G/P$ mit $G$ einfach. Es wird au\ss
erdem gezeigt, da\ss\ die Vermutung f\"ur die von Kapranov f\"ur
Grassmannsche Mannigfaltigkeiten und Quadriken angegebenen Folgen
richtig ist. Es wird die Matrix der durch die Eulercharakteristik
gegebenen Bilinearform auf der Grothendieckschen $K$-Gruppe
$K_{\circ}(X)$ in der Basis der Klassen der Strukturgarben von
Schubertvariet\"aten in $X$ berechnet, die zu der Grammatrix einer
vollst\"andigen exzeptionellen Folge konjugiert ist. In Abschnitt 3
wird der Satz 3.2.7 bewiesen, der vollst\"andige exzeptionelle Folgen
auf Quadrikenb\"undeln \"uber solchen Basismannigfaltigkeiten liefert,
auf denen man vollst\"andige exzeptionelle Folgen bereits kennt. Damit
wird die Klasse von Variet\"aten (insbesondere rational-homogenen
Mannigfaltigkeiten), auf denen man solche Folgen kennt, wesentlich
erweitert. Im Rest des Abschnitts 3 werden Variet\"aten isotroper
Fahnen in einem symplektischen bzw. orthogonalen Vektorraum
betrachtet. Das Problem, auf diesen vollst\"andige exzeptionelle
Folgen zu konstruieren, reduziert sich mit einem Satz von Orlov (Satz
3.1.5) auf isotrope Grassmannsche. Satz 3.3.3 gibt f\"ur letztere im
symplektischen Fall Erzeugendensysteme der derivierten Kategorie an,
die aus homogenen Vektorb\"undeln bestehen und aus denen man in
Spezialf\"allen vollst\"andige exzeptionelle Folgen konstruieren
kann. Im Unterabschnitt 3.4 wird gezeigt, wie sich die Methode unter
Benutzung von Satz 3.2.7 auf den orthogonalen Fall ausdehnen l\"a\ss
t. Insbesondere wird Satz 3.4.1 bewiesen, der ein Erzeugendensystem
f\"ur die derivierte Kategorie koh\"renter Garben auf der
Grassmannschen der dreidimensionalen isotropen Unterr\"aume in einem
7-dimensionalen orthogonalen Vektorraum liefert. Abschnitt 4 ist der
geometrischen Motivation f\"ur Cataneses Vermutung gewidmet und enth\"alt einen alternativen Zugang zur Konstruktion
vollst\"andiger exzeptioneller Folgen auf rational-homogenen
Mannigfaltigkeiten, der auf einem Satz von M. Brion (Satz 4.1.1) und
zellul\"aren Aufl\"osungen von Monomidealen \`{a} la Bayer/Sturmfels
basiert. Es wird ein neuer Beweis f\"ur den Satz von Beilinson auf dem
$\mathbb{P}^n$ gegeben, um zu zeigen, da\ss\ dieser Zugang im
allgemeinen Fall funktionieren k\"onnte. \"Uberdies wird Satz 4.2.5
bewiesen, der eine konkrete Beschreibung gewisser Funktoren liefert,
die bei diesem Ansatz studiert werden m\"ussen.

\mdseries\normalsize
\newpage

\tableofcontents
\newpage

\mdseries
\normalsize
\section{Introduction}

The concept of derived category of an Abelian category $\mathcal{A}$,
which gives a transparent and compact way to handle the totality of
cohomological data attached to $\mathcal{A}$ and puts a given object
of $\mathcal{A}$ and all of its resolutions on equal footing, was
conceived by Grothendieck at the beginning of the 1960's and their
internal structure was axiomatized by Verdier through the notion of
triangulated category in his 1967 thesis (cf. \cite{Ver1},
\cite{Ver2}). Verdier's axioms for distinguished triangles still allow
for some pathologies (cf. \cite{GeMa}, IV.1, 7) and in \cite{BK} it
was suggested how to replace them by more satisfactory ones, but since
the former are in current use, they will also be the basis of this
text. One may consult \cite{Nee} for foundational questions on
triangulated categories.\\
However, it was only in 1978 that people laid hands on ``concrete''
derived categories of geometrical significance (cf. \cite{Bei} and
\cite{BGG2}), and A. A. Beilinson constructed strong complete
exceptional sequences of vector bundles for $D^b(Coh\,
\mathbb{P}^n)$, the bounded derived category of coherent sheaves on
$\mathbb{P}^n$. The terminology is explained in section 2, def. 2.1.7, below, but
roughly the simplification brought about by Beilinson's theorem is
analogous to the construction of a semi-orthonormal basis $(e_1, \dots
, e_d)$ for a vector
space equipped with a non-degenerate (non-symmetric) bilinear form
$\chi$ (i.e., $\chi (e_i,e_i)=1\: \forall i$, $\chi (e_j, e_i)=0\:
\forall j>i)$).\\
Beilinson's theorem represented a spectacular breakthrough and, among other things, his technique was applied to the study of
moduli spaces of semi-stable sheaves of given rank and Chern classes
on $\mathbb{P}^2$ and $\mathbb{P}^3$ by Horrocks, Barth/Hulek,
Dr\'{e}zet/Le Potier (cf. \cite{OSS}, \cite{Po} and references
therein).\\
Recently, A. Canonaco has obtained a generalization of Beilinson's
theorem to weighted projective spaces and applied it to the study of
canonical projections of surfaces of general type on a 3-dimensional
weighted projective space (cf. \cite{Can}, cf. also \cite{AKO}). \\
From 1984 onwards, in a series of papers \cite{Ka1}, \cite{Ka2},
\cite{Ka3}, M. M. Kapranov found strong complete exceptional
sequences on Grassmannians and flag varieties of type $A_n$ and on
quadrics. Subsequently, exceptional sequences alongside with
some new concepts introduced in the meantime such as helices,
mutations, semi-orthogonal decompositions etc. were intensively
studied, 
especially in Russia, an account of which can be found in
the volume \cite{Ru1} summarizing a series of seminars conducted by
A. N. Rudakov in Moscow (cf. also \cite{Bo}, \cite{BoKa},
\cite{Or}). Nevertheless, despite the wealth of new techniques
introduced in the process, many basic questions concerning exceptional
sequences are still very much open. These fall into two main
classes: first questions of existence: E.g., do complete exceptional sequences always exist on
rational homogeneous manifolds? (For toric varieties existence of
complete exceptional sequences was proven
very recently by Kawamata, cf. \cite{Kaw}.) Secondly, one often does not know if
basic intuitions derived from semi-orthogonal linear algebra hold true
in the framework of exceptional sequences, and thus one does not have
enough flexibility to manipulate them, e.g.: Can every exceptional bundle on a variety
$X$ on which complete exceptional sequences are known to exist
(projective spaces, quadrics...) be included in a complete exceptional
sequence?\\
To round off this brief historical sketch, one should not forget to
mention that derived
categories have proven to be of geometrical significance in a lot of
other contexts, e.g. through Fourier-Mukai transforms and the
reconstruction theorem of Bondal-Orlov for smooth projective varieties
with ample canonical or anti-canonical class (cf. \cite{Or2}), in the
theory of perverse sheaves and the generalized Riemann-Hilbert
correspondence (cf. \cite{BBD}), or in the recent proof of
T. Bridgeland that birational Calabi-Yau threefolds have equivalent
derived categories and in particular the same Hodge numbers
(cf. \cite{Brid}). Interest in derived categories was also extremely
stimulated by M. Kontsevich's proposal for homological mirror symmetry
(\cite{Kon}) on the one side and by new applications to minimal model
theory on the other side.\\
Let me now describe the aim and contents of this work. Roughly
speaking, the problem is
to give as concrete as possible a description of the (bounded) derived
categories of coherent sheaves on rational homogeneous manifolds
$X=G/P$, $G$ a connected complex semisimple Lie group, $P\subset G$ a
parabolic subgroup. More precisely, the following set of main questions and
problems, ranging from the modest to the more ambitious, have served
as programmatic guidelines:
\begin{itemize}
\item[\textbf{P 1.}]
Find generating sets of $D^b(Coh\, X)$ with as few elements as
possible. (Here a set of elements of $D^b(Coh\, X)$ is called a
\emph{generating set} if the smallest full triangulated subcategory
containing this set is equivalent to $D^b(Coh\, X)$). 
\end{itemize}
We will see in subsection 2.3 below that the number of elements in a
generating set is always bigger or equal to the number of Schubert
varieties in $X$.\\
In the next two problems we mean by a complete exceptional sequence an ordered tuple $(E_1,
\dots , E_n)$ of objects $E_1, \dots , E_n$ of $D^b(Coh\, X)$ which
form a generating set and such that moreover $R^{\bullet}\mathrm{Hom}(E_i,
E_j)=0$ for all $i>j$, $R^{\bullet}\mathrm{Hom}(E_i,
E_i)=\mathbb{C}$ (in degree $0$) for all $i$. If in addition all
extension groups in nonzero degrees between the elements $E_i$ are
absent we speak of a strong complete exceptional sequence. See section
2, def. 2.1.7, for further discussion. 
\begin{itemize}
\item[\textbf{P 2.}] 
Do there always exist complete exceptional sequences in $D^b(Coh\,
X)$?
\item[\textbf{P 3.}]
Do there always exist strong complete exceptional sequences in\\ 
$D^b(Coh\, X)$?
\end{itemize}
Besides the examples found by Kapranov mentioned above, the only other
substantially different example I know of in answer to \textbf{P 3.}
is the one given by A. V. Samokhin in \cite{Sa} for the Lagrangian
Grassmannian of totally isotropic 3-planes in a 6-dimensional
symplectic vector space.\\
In the next problem we mean by a complete strong exceptional poset a set of objects 
$\{E_1 ,\dots , E_n\}$ of $D^b(Coh\, X)$ that generate $D^b(Coh\, X)$
and satisfy
\\ $R^{\bullet}\mathrm{Hom}(E_i,
E_i)=\mathbb{C}$ (in degree $0$) for all $i$ and such that all
extension groups in nonzero degrees between the $E_i$ are absent,
together with a partial order $\le$ on $\{E_1 ,\dots , E_n\}$ subject
to the condition: $\mathrm{Hom}( E_j, E_i)=0$ for $j\ge i,\: j\neq
i$ (cf. def. 2.1.7 (B)). 
\begin{itemize}
\item[\textbf{P 4.}]
Catanese's conjecture: On any $X=G/P$ there exists a complete strong
exceptional poset $(\{ E_1, \dots , E_n\} , \le )$ together with a
bijection of the elements of the poset with the Schubert varieties in
$X$ such that $\le $ is the partial order induced by the
Bruhat-Chevalley order (cf. conj. 2.2.1 (A)).
\end{itemize}

\begin{itemize}
\item[\textbf{P 5.}]
Dubrovin's conjecture (cf. \cite{Du}, conj. 4.2.2; slightly modified
afterwards in \cite{Bay}; cf. also \cite{B-M}): The (small) quantum
cohomology of a smooth projective variety $Y$ is generically
semi-simple if and only if there exists a complete exceptional
sequence in $D^b(Coh\, Y)$ (Dubrovin also relates the Gram matrix of
the exceptional sequence to quantum-cohomological data but we omit
this part of the conjecture).   
\end{itemize}
Roughly speaking, quantum cohomology endows the usual cohomology space with
complex coefficients $H^{\ast}(Y)$ of $Y$ with a new commutative associative
multiplication $\circ_{\omega}: H^{\ast}(Y)\times H^{\ast}(Y)\to
H^{\ast}(Y)$ depending on a complexified K\"ahler class $\omega\in
H^2(Y, \mathbb{C})$, i.e. the imaginary part of $\omega$ is in the
K\"ahler cone of $Y$ (here we assume $H^{odd}(Y)=0$ to avoid working
with supercommutative rings). The condition that the quantum
cohomology of $Y$ is generically semi-simple means that for generic
values of $\omega$ the resulting algebra is semi-simple. The validity
of this conjecture would provide further evidence for the famous
homological mirror conjecture by Kontsevich (\cite{Kon}). However, we
will not deal with quantum cohomology in this work.\\
Before stating the results, a word of explanation is in order to
clarify why we
narrow down the focus to rational homogeneous manifolds:
\begin{itemize}
\item
Exceptional vector bundles need not always exist on an arbitrary
smooth projective variety; e.g., if the canonical class of $Y$ is
trivial, they never exist (see the explanation following definition
2.1.3).
\item
$D^b(Coh \, Y)$ need not be finitely generated, e.g., if $Y$ is an
Abelian variety (see the explanation following definition
2.1.3).
\item
If we assume that $Y$ is Fano, then the Kodaira vanishing theorem
tells us that all line bundles are exceptional, so we have at least
some \emph{a priori} supply of exceptional bundles.
\item
Within the class of Fano manifolds, the rational homogeneous spaces
$X=G/P$ are distinguished by the fact that they are amenable to
geometric, representation-theoretic and combinatorial methods alike. 
\end{itemize}
Next we will state right away the main results obtained, keeping the
numbering of the text and adding a word of explanation to each.\\
Let $V$ be a $2n$-dimensional symplectic vector space and
$\mathrm{IGrass}(k , V)$ the Grassmannian of $k$-dimensional isotropic
subspaces of $V$ with tautological subbundle $\mathcal{R}$. $\Sigma^{\bullet}$
denotes the Schur functor (see subsection 2.2 below for explanation). 
\begin{theo3.3.3}
The derived category $D^b(Coh (\mathrm{IGrass}(k , V)))$ is generated by the bundles $\Sigma^{\nu } 
\mathcal{R}$, where $\nu$ runs over Young diagrams $Y$ which satisfy
\begin{gather*}
\left( \mathrm{number}\;\mathrm{of}\;\mathrm{columns}\;\mathrm{of}\;  Y \right)\le 2n-k\, ,\\
k \ge \left( \mathrm{number}\;\mathrm{of}\;\mathrm{rows}\;\mathrm{of}\;  Y \right) \ge 
\left( \mathrm{number}\;\mathrm{of}\;\mathrm{columns}\;\mathrm{of}\;  Y \right) - 2 (n-k)\, .
\end{gather*}
\end{theo3.3.3}
This result pertains to \textbf{P 1}. Moreover, we will see in
subsection 3.3 that \textbf{P 2.} for isotropic flag manifolds of type $C_n$
can be reduced to \textbf{P 2.} for isotropic Grassmannians. Through
examples 3.3.6-3.3.8 we show that theorem 3.3.3 gives a set of bundles
which is in special cases manageable enough to obtain from it a
complete exceptional sequence. In general, however, this last step
is a difficult combinatorial puzzle relying on Bott's theorem for the
cohomology of homogeneous bundles and Schur complexes derived from
tautological exact sequences on the respective Grassmannians.\\
For the notion of semi-orthogonal decomposition in the next theorem we
refer to definition 2.1.17 and for the definition of the spinor
bundles $\Sigma$, $\Sigma^{\pm}$ of the orthogonal vector bundle
$\mathcal{O}_{\mathcal{Q}}(-1)^{\perp }/\mathcal{O}_{\mathcal{Q}}(-1)$
  we refer to subsection 3.2.
\begin{theo3.2.7}
Let $X$ be a smooth projective variety with $H^1(X ; \mathbb{Z}/2\mathbb{Z})=0$, $\mathcal{E}$ an orthogonal
vector bundle of rank $r+1$ on $X$ (i.e., $\mathcal{E}$ comes equipped
with a quadratic form $q\in\Gamma (X, \mathrm{Sym}^2
\mathcal{E}^{\vee})$ which is non-degenerate on each fibre), 
$\mathcal{Q}\subset \mathbb{P}(\mathcal{E})$ the
associated quadric bundle, and let $\mathcal{E}$ carry a spin structure.\\
Then there is a semiorthogonal decomposition
\begin{gather*}
D^b(\mathcal{Q})=\left\langle  D^b(X)\otimes \Sigma (-r+1) , D^b(X)\otimes
\mathcal{O}_{\mathcal{Q}}(-r+2) ,\right. \\
\left. \ldots , D^b(X)\otimes \mathcal{O}_{\mathcal{Q}}(-1) , D^b(X) 
\right\rangle
\end{gather*}
for $r+1$ odd and
\begin{gather*}
D^b(\mathcal{Q})=\left\langle  D^b(X)\otimes \Sigma^{+}(-r+1) , D^b(X)\otimes
\Sigma^{-}(-r+1) ,\right. \\ \left. D^b(X)\otimes
\mathcal{O}_{\mathcal{Q}}(-r+2) , \ldots , D^b(X)\otimes
\mathcal{O}_{\mathcal{Q}}(-1) , D^b(X)
\right\rangle 
\end{gather*}
for $r+1$ even.
\end{theo3.2.7}
This theorem is an extension to the relative case of  a theorem of
\cite{Ka2}. It enlarges substantially the class of varieties (especially
rational-homogeneous varieties) on which complete exceptional
sequences are proven to exist (\textbf{P 2}). It will
also be the substantial ingredient in subsection 3.4: Let $V$ be a
7-dimensional orthogonal vector space, $\mathrm{IGrass}(3 , V)$ the
Grassmannian of isotropic 3-planes in $V$, $\mathcal{R}$ the
tautological subbundle on it; $L$ denotes the ample generator of
$\mathrm{Pic} (\mathrm{IGrass}(3 , V))\simeq \mathbb{Z}$ (a square
root of $\mathcal{O}(1)$ in the Pl\"ucker embedding). For more information cf. subsection 3.4.
\begin{theo3.4.1}
The derived category $D^b(Coh\, \mathrm{IGrass}(3 , V))$ is generated
as triangulated category by the following 22 vector bundles:
\begin{gather*}
\bigwedge\nolimits^2 \mathcal{R} (-1), \:\mathcal{O}(-2),\:
\mathcal{R}(-2)\otimes L, \: \mathrm{Sym}^2 \mathcal{R}(-1)\otimes L ,
\:
\mathcal{O}(-3)\otimes L ,\\ \bigwedge\nolimits^2
\mathcal{R}(-2)\otimes L ,\: \Sigma^{2 ,1}\mathcal{R}(-1)\otimes L ,
\: \mathcal{R}(-1) ,\: \mathcal{O}(-2)\otimes L ,\: \mathcal{O}(-1)
,\\
\mathcal{R}(-1)\otimes L ,\; \bigwedge\nolimits^2
\mathcal{R}(-1)\otimes L , \: \Sigma^{2 , 1}\mathcal{R}\otimes L , \:
\mathrm{Sym}^2 \mathcal{R}^{\vee} (-2)\otimes L ,\:
\bigwedge\nolimits^2 \mathcal{R} , \: \mathcal{O} , \\
\Sigma^{2 , 1}\mathcal{R} , \: \mathrm{Sym}^2 \mathcal{R}^{\vee}(-2) ,
\: \mathcal{O}(-1)\otimes L , \: \mathrm{Sym}^2 \mathcal{R}^{\vee}
(-1) , \: \bigwedge\nolimits^2 \mathcal{R}\otimes L , \:
\mathcal{R}\otimes L .
\end{gather*}
\end{theo3.4.1}
This result pertains to \textbf{P 1.} again. One should remark that
\textbf{P 2.} for isotropic flag manifold of type $B_n$ or $D_n$ can
again be 
reduced to isotropic Grassmannians. Moreover, the method of subsection
3.4 applies to all orthogonal isotropic Grassmannians alike, but since
the computations tend to become very large, we restrict our
attention to a particular case.\\ 
Beilinson proved his theorem on $\mathbb{P}^n$ using a resolution of
the structure sheaf of the diagonal and considering the functor $Rp_{2\ast}( p_1^{\ast}(-) \otimes^L 
\mathcal{O}_{\Delta})\simeq \mathrm{id}_{D^b(Coh\, \mathbb{P}^n)}$  (here $p_1 ,p_2: \mathbb{P}^n\times
  \mathbb{P}^n\to \mathbb{P}^n$ are the projections onto the two factors). The situation is complicated on general
  rational homogeneous manifolds $X$ because resolutions of the
  structure sheaf of the diagonal $\Delta\subset X\times X$ analogous
  to those used in \cite{Bei}, \cite{Ka1}, \cite{Ka2}, \cite{Ka3} to
  exhibit complete exceptional sequences, are not known. The
  preceding theorems are proved by ``fibrational techniques''. Section
  4 outlines an alternative approach: In fact, M. Brion (\cite{Bri}) constructed,
  for any rational homogeneous manifold $X$, a degeneration of the
  diagonal $\Delta_X$ into $\mathfrak{X}_0$, which is a union, over
  the Schubert varieties in $X$, of the products of a Schubert variety
  with its opposite Schubert variety (cf. thm. 4.1.1). It turns out that it is
  important to describe the functors $Rp_{2\ast}( p_1^{\ast}(-) \otimes^L 
\mathcal{O}_{\mathfrak{X}_0})$ which, contrary to what one might
expect at first glance, are no longer isomorphic to the identity
functor (one might hope to reconstruct the identity out of $Rp_{2\ast}( p_1^{\ast}(-) \otimes^L 
\mathcal{O}_{\mathfrak{X}_0})$ and some infinitesimal data attached to
the degeneration). For $\mathbb{P}^n$ this is
accomplished by the following 
\begin{theo4.2.5}
Let $\{ pt\}=L_0\subset L_1\subset \dots \subset
L_n=\mathbb{P}^n$ be a full flag of projective linear subspaces of
$\mathbb{P}^n$ (the Schubert varieties in $\mathbb{P}^n$) and denote
by $L^j$ the Schubert variety opposite to $L_j$.\\
For $d\ge 0$ one has in $D^b(Coh\, \mathbb{P}^n)$
\begin{gather*}
Rp_{2\ast}( p_1^{\ast}(\mathcal{O}(d)) \otimes^L 
\mathcal{O}_{\mathfrak{X}_0})\simeq \bigoplus_{j=0}^n
\mathcal{O}_{L_j}\otimes H^0(L^j, \mathcal{O}(d))^{\vee}/ H^0(L^{j+1},
\mathcal{O}(d))^{\vee} \, .
\end{gather*}
Moreover, one can also describe completely the effect of $Rp_{2\ast}( p_1^{\ast}(-) \otimes^L 
\mathcal{O}_{\mathfrak{X}_0})$ on morphisms (cf. subsection 4.2 below).
\end{theo4.2.5}
The proof uses the technique of cellular resolutions of monomial
ideals of Bayer and Sturmfels (\cite{B-S}). We also show in subsection
4.2 that Beilinson's theorem on $\mathbb{P}^n$ can be recovered by our
method with a proof that uses only $\mathfrak{X}_0$ (see remark 4.2.6).\\
It should be added that we will not completely ignore the second part
of \textbf{P 4.} concerning Hom-spaces: In section 2 we show that the
conjecture in \textbf{P 4.} is valid in full for the complete strong
exceptional sequences found by Kapranov on Grassmannians and
quadrics (cf. \cite{Ka3}). In remark 2.3.8 we discuss a possibility
for relating the Gram matrix of a strong complete exceptional
sequence on a rational homogeneous manifold with the
Bruhat-Chevalley order on Schubert cells.\\ 
Additional information about the content of each section can be found
at the beginning of the respective section.\\
\textbf{Acknowledgements.} I would like to thank my thesis advisor
Fabrizio Catanese for posing the problem and several discussions on
it. Special thanks also to Michel Brion for filling in my insufficient
knowledge of representation theory and algebraic groups on
a number of occasions and for fruitful suggestions and discussions.

\section{Tools and background: getting off the ground}

This section supplies the concepts and dictionary that will be used throughout the text. We state 
a conjecture due to F. Catanese which was the motivational backbone of this work and discuss its relation 
to work of M. M. Kapranov. Moreover, we prove some results that are useful in the study of the derived 
categories of coherent sheaves on rational homogeneous varieties, but
do not yet tackle the problem of constructing 
complete exceptional sequences on them: This will be the subject matter of sections 3 and 4.   

\subsection{Exceptional sequences}

Throughout the text we will work over the ground field $\mathbb{C}$
of complex numbers.\\
The classical theorem of Beilinson (cf. \cite{Bei}) can be stated as follows.

\begin{theorem}
Consider the following two ordered sequences of sheaves on $\mathbb{P}^n=\mathbb{P}(V)$, $V$ an $n+1$ 
dimensional vector space:
\begin{eqnarray*}
\mathfrak{B}=(\mathcal{O}(-n), \ldots , \mathcal{O}(-1), \, \mathcal{O} ) \\
\mathfrak{B}'= \left( \Omega^n (n), \ldots , \Omega^1(1) ,\, \mathcal{O}\right) .
\end{eqnarray*}
Then $D^b(Coh \mathbb{P}^n)$ is equivalent as a triangulated category to the homotopy category of bounded 
complexes of sheaves on $\mathbb{P}^n$ whose terms are finite direct sums of sheaves in $\mathfrak{B}$ (and 
the same for $\mathfrak{B}$ replaced with $\mathfrak{B}'$).\\
Moreover, one has the following stronger assertion: If $\Lambda= \bigoplus_{i=0}^{n+1}\wedge^i V$ and $S= 
\oplus_{i=0}^{\infty} \mathrm{Sym}^i V^{\ast}$ are the $\mathbb{Z}$-graded exterior algebra of $V$, resp. 
symmetric algebra of $V^{\ast}$, and $K_{[0, n]}^b\Lambda$ resp. $K_{[0, n]}^b S$ are the homotopy categories 
of bounded complexes whose terms are finite direct sums of free modules $\Lambda [i]$, resp. $S[i]$, for 
$0\le i\le n$, and whose morphisms are homogeneous graded of degree $0$, then 
\[
K_{[0, n]}^b\Lambda \simeq D^b(Coh \mathbb{P}^n)\quad K_{[0, n]}^b S \simeq D^b(Coh \mathbb{P}^n)
\]
as triangulated categories, the equivalences being given by sending $\Lambda [i]$ to $\Omega^i(i)$ and 
$S[i]$ to $\mathcal{O}(-i)$ ($\Lambda [i]$, $S[i]$ have their generator in degree $i$). 
\end{theorem}

One would like to have an analogous result on any rational homogeneous variety $X$, i.e. a rational projective 
variety with a transitive Lie group action or equivalently (cf. \cite{Akh}, 3.2, thm. 2) a coset manifold 
$G/P$ where $G$ is a connected semisimple complex Lie group (which can
be assumed to be simply connected) and $P
\subset G$ is a parabolic subgroup. However, to give a precise meaning to this wish, one should first try to 
capture some formal features of Beilinson's theorem in the form of suitable definitions; thus we will recall next 
a couple of notions which have become standard by now, taking theorem 2.1.1 as a model.\\
Let $\mathcal{A}$ be an Abelian category.

\begin{definition}
A class of objects $\mathcal{C}$ \emph{generates} $D^b(\mathcal{A})$ if the smallest full triangulated 
subcategory containing the objects of $\mathcal{C}$ is equivalent to
$D^b(\mathcal{A})$. If $\mathcal{C}$ is a set, we will also call
$\mathcal{C}$ a \emph{generating set} in the sequel.
\end{definition}

Unravelling this definition, one finds that this is equivalent to saying that, up to isomorphism, every object 
in $D^b(\mathcal{A} )$ can be obtained by successively enlarging $\mathcal{C}$ through the following operations: 
Taking finite direct sums, shifting in $D^b(\mathcal{A})$ (i.e., applying the translation functor), and taking a 
cone $Z$  of a morphism $u : X\to Y$ between objects already
constructed: This means we complete $u$ to a 
distinguished triangle $X\stackrel{u}{\longrightarrow} Y \longrightarrow Z \longrightarrow X[1]$.\\
The sheaves $\Omega^i(i)$ and $\mathcal{O}(-i)$ in theorem 2.1.1 have the distinctive property of being 
``exceptional''.

\begin{definition}
An object $E$ in $D^b(\mathcal{A})$ is said to be \emph{exceptional} if
\[
\mathrm{Hom} (E , E) \simeq \mathbb{C} \quad \mathrm{and}\quad \mathrm{Ext}^i(E , E) =0 \; \forall i \neq 0 .
\]
\end{definition}

If $Y$ is a smooth projective variety, exceptional objects need not always exist (e.g., if $Y$ has trivial 
canonical class this is simply precluded by Serre duality since then
$\mathrm{Hom}(E, E)\simeq \mathrm{Ext}^n(E, E)\neq 0$).\\
What is worse, $D^b(Coh\, Y)$ need not even possess a finite generating
set: In fact we will see in subsection 2.3 below that if $D^b(Coh\,
Y)$ is finitely generated, then $A(Y)=\bigoplus_{r=0}^{\dim
  Y}A^r(Y)$, the Chow ring of $Y$ of algebraic cycles modulo rational equivalence, is finitely generated as an abelian
group (here $A^r(Y)$ denotes the group of cycles of codimension $r$ on
$Y$ modulo rational equivalence). But, for instance, if $Y$ is an Abelian variety, $A^1(Y)\simeq
\mathrm{Pic}\, Y$ is not finitely generated.\\
Recall that a vector bundle $\mathcal{V}$ on a rational homogeneous variety $X=G/P$ is called 
$G$-\emph{homogeneous} if there is a $G$-action on $\mathcal{V}$ which lifts the $G$-action on $X$ and is 
linear on the fibres. It is well known that this is equivalent to saying that $\mathcal{V}\simeq G\times_{\varrho}V
$, where $\varrho : P\to \mathrm{GL} (V)$ is some representation of the algebraic group $P$ and $G\times_{\varrho}V$ 
is the quotient of $G\times V$ by the action of $P$ given by $p\cdot (g,v):= (gp^{-1}, \varrho(p)v)$, $p\in P$, 
$g\in G$, $v\in V$. The projection to $G/P$ is induced by the
projection of $G\times V$ to $G$; this construction gives a 1-1
correspondence between representations of the subgroup $P$ and homogeneous vector bundles over $G/P$ (cf. 
\cite{Akh}, section 4.2).\\
Then we have the following result (mentioned briefly in a number of places, e.g. 
\cite{Ru1}, 6., but without a precise statement or proof).

\begin{proposition}
Let $X=G/P$ be a rational homogeneous manifold with $G$ a simply connected semisimple group, and let 
$\mathcal{F}$ be an exceptional sheaf on $X$. Then $\mathcal{F}$ is a $G$-homogeneous bundle. 
\end{proposition}
\begin{proof}
Let us first agree that a deformation of a coherent sheaf $\mathcal{G}$ on a complex space $Y$ is a triple 
$(\tilde{\mathcal{G}}, S, s_0)$ where $S$ is another complex space (or germ), $s_0\in S$, $\tilde{\mathcal{G}}$ 
is a coherent sheaf on $Y\times S$, flat over $S$, with $\tilde{\mathcal{G}}\mid_{Y\times\{ s_0\}} \simeq 
\mathcal{G}$ and $\mathrm{Supp}\,\tilde{ \mathcal{G}}\to S$
proper. Then one knows that, for the deformation with base a complex
 space germ, there is a versal deformation and its tangent space at
 the marked point is $\mathrm{Ext}^1(\mathcal{G} , 
\mathcal{G} )$ (cf. \cite{S-T}).\\
Let $\sigma : G\times X\to X$ be the group action; then $(\sigma^{\ast}\mathcal{F}, G , \mathrm{id}_G)$ is a 
deformation of $\mathcal{F}$ (flatness can be seen e.g. by embedding
$X$ equivariantly in a projective space 
(cf. \cite{Akh}, 3.2) and noting that the Hilbert polynomial of $\sigma^{\ast}\mathcal{F}\mid_{\{g\}\times X} = 
\tau_g^{\ast}\mathcal{F}$ is then constant for $g\in G$; here $\tau_g : X\to X$ is the automorphism induced 
by $g$). Since $\mathrm{Ext}^1(\mathcal{F}, \mathcal{F})=0$ one has by the above that $\sigma^{\ast}\mathcal{F}
$ will be locally trivial over $G$, i.e. $\sigma^{\ast}\mathcal{F}\simeq \mathrm{pr}_2^{\ast}\mathcal{F}$
 locally 
over $G$ where $\mathrm{pr}_2 : G\times X\to X$ is the second projection ($\mathcal{F}$ is ``rigid''). In 
particular $\tau_g^{\ast}\mathcal{F}\simeq \mathcal{F}$ $\forall g\in G$.\\
Since the locus of points where $\mathcal{F}$ is not locally free is a proper algebraic subset of $X$ and 
invariant under $G$ by the preceding statement, it is empty because $G$ acts transitively. Thus $\mathcal{F}$ 
is a vector bundle satisfying $\tau_g^{\ast}\mathcal{F}\simeq \mathcal{F}$ $\forall g\in G$. Since $G$ is 
semisimple and assumed to be simply connected, this is enough to imply that $\mathcal{F}$ is a $G$-homogeneous bundle 
(a proof of this last assertion due to A. Huckleberry is presented in \cite{Ot2} thm. 9.9).    
\end{proof}
\begin{remark}
In proposition 2.1.4 one must insist that $G$ be simply connected as an example in \cite{GIT}, ch.1, \S 3 shows
: The exceptional bundle $\mathcal{O}_{\mathbb{P}^n}(1)$ on $\mathbb{P}^n$ is $\mathrm{SL}_{n+1}$-homogeneous, 
but not homogeneous for the adjoint form $PGL_{n+1}$ with its action $\sigma : PGL_{n+1}\times\mathbb{P}^n \to 
\mathbb{P}^n$ since the $\mathrm{SL}_{n+1}$-action on
$H^0(\mathcal{O}_{\mathbb{P}^n}(1))$ does not factor through $PGL_{n+1}$.

\end{remark}

\begin{remark}
It would be interesting to know which rational homogeneous manifolds
$X$ enjoy the property that exceptional objects 
in $D^b(Coh\,  X)$ are actually just shifts of exceptional sheaves. It is straightforward to check 
that this is true on $\mathbb{P}^1$. This is because, if $C$ is a curve, $D^b(Coh\, C)$ is not very 
interesting: In fancy language, the underlying abelian category is \emph{hereditary} which means 
$\mathrm{Ext}^2(\mathcal{F}, \mathcal{G})=0\; \forall \mathcal{F}, \mathcal{G} \in \mathrm{obj} \, (Coh\, C)$.
 It is easy to see (cf. \cite{Ke}, 2.5) that then every object $Z$ in $D^b(Coh \, C)$ is isomorphic to the
 direct sum of shifts of its cohomology sheaves $\bigoplus_{i\in\mathbb{Z}} H^i(Z)[-i]$ whence morphisms 
between objects $Z_1$ and $Z_2$ correspond to tuples $(\varphi_i, e_i)_{i\in\mathbb{Z}}$ with $\varphi_i : 
H^i (Z_1) \to H^i (Z_2)$ a sheaf morphism and $e_i\in \mathrm{Ext}^1( H^i(Z_1), H^{i-1}(Z_2))$ an extension class
. Exceptional objects are indecomposable since they are simple.\\ 
The same property holds on $\mathbb{P}^2$ (and more generally on any Del Pezzo 
surface) by \cite{Gor}, thm. 4.3.3, and is conjectured to be true on $\mathbb{P}^n$ in general (\cite{Gor}, 3.2.7).
\end{remark}

The sequences $\mathfrak{B}$ and $\mathfrak{B}'$ in theorem 2.1.1 are examples of complete strong exceptional 
sequences (cf. \cite{Ru1} for the development of this notion).

\begin{definition}
\begin{itemize}
\item[(A)]
An $n$-tuple $(E_1, \ldots , E_n)$ of exceptional objects in $D^b(\mathcal{A})$ is called an \emph{exceptional 
sequence} if
\[
\mathrm{Ext}^l (E_j , E_i )= 0 \; \quad\forall 1\le i < j \le n \;\;\mathrm{and}\;\; \forall l\in \mathbb{Z} \, .
\] 
If in addition 
\[
\mathrm{Ext}^l (E_j , E_i)= 0 \; \quad \forall 1\le i, j \le n \;\; \mathrm{and} \;\; \forall l\neq 0
\]
we call $(E_1, \ldots , E_n)$ a \emph{strong exceptional sequence}. The sequence is \emph{complete} if 
$E_1, \ldots , E_n$ generate $D^b( \mathcal{A})$. 
\item[(B)]
In order to phrase conjecture 2.2.1 below precisely, it will be 
convenient to introduce also the following terminology: A set of
exceptional 
objects $\{ E_1 , \dots , E_n\}$ in $D^b(\mathcal{A})$ that generates
$D^b(\mathcal{A})$ and such that $\mathrm{Ext}^l (E_j , E_i)=0$ for
all $1\le i , j \le n$ and all $l\neq 0$ will be called a \emph{complete
strong exceptional set}. A partial order $\le $ on a complete strong 
exceptional set is \emph{admissible} if $\mathrm{Hom}(E_j
, E_i)=0$ for all $j\ge i ,\: i\neq j$. A pair $( \{ E_1 , \dots ,
E_n\} , \le )$ consisting of a complete strong exceptional set and an
admissible partial order on it will be called  a \emph{complete strong
exceptional poset}.
\item[(C)]
A \emph{complete very strong exceptional poset} is a pair $( \{ E_1 , \dots ,
E_n\} , \le )$ where $\{ E_1 , \dots , E_n\}$ is a complete strong
exceptional set and $\le$ is a partial order on this set such that $\mathrm{Hom}(E_j
, E_i)=0$ unless $i \ge j$.   
\end{itemize}
\end{definition}
Obviously every complete strong exceptional sequence is a complete
strong exceptional poset (with the partial order being in fact a total
order). I think it might be possible that for complete strong
exceptional posets in $D^b( Coh \, X)$ which consist of vector bundles, $X$ a rational homogeneous
manifold, the converse holds, i.e. any admissible partial order can be
refined to a total order which makes the poset into a complete strong
exceptional sequence. But I cannot prove this.\\
Moreover, every complete very strong exceptional poset is in
particular a complete strong exceptional poset. If we choose a total
order refining the partial order on a complete very strong exceptional
poset, we obtain a complete strong exceptional sequence.\\
Let me explain the usefulness of these concepts by first saying what kind of analogues of Beilinson's theorem 2.1.1 
we can expect for $D^b(\mathcal{A})$ once we know the existence of a complete strong exceptional set. 

Look at a complete strong exceptional set $\{ E_1, \ldots , E_n \}$ in $D^b(\mathcal{A})$ consisting 
of objects $E_i$, $1\le i \le n$, of $\mathcal{A}$. If $K^b(\{ E_1 ,\ldots , E_n\} )$ denotes the homotopy 
category of bounded complexes in $\mathcal{A}$ whose terms are finite
direct sums of the $E_i$'s, it is clear
 that the natural functor
\[
\Phi_{(E_1,\ldots , E_n)} : K^b(\{ E_1 ,\ldots , E_n\} ) \to D^b( \mathcal{A})
\]
(composition of the inclusion $K^b(\{ E_1 ,\ldots , E_n\} ) \hookrightarrow K^b(\mathcal{A})$ with the 
localization $Q : K^b(\mathcal{A}) \to D^b(\mathcal{A})$) is an equivalence; indeed $\Phi_{(E_1,\dots , E_n)}$ 
is essentially surjective because $\{ E_1,\ldots , E_n \}$ is complete and $\Phi_{(E_1,\dots , E_n)}$ is fully 
faithful because $\mathrm{Ext}^p (E_i , E_j)=0$ for all $p> 0$ and all $i$ and $j$ implies 
\begin{gather*}
\mathrm{Hom}_{K^b(\{ E_1 ,\ldots , E_n\} )} (A , B) \simeq \mathrm{Hom}_{D^b(\mathcal{A})} ( \Phi_{(E_1,\dots ,
 E_n)}\, A , \Phi_{(E_1,\dots , E_n)}\, B) \\ \quad \forall A,\, B \in \mathrm{obj}\, K^b(\{ E_1 ,\ldots , E_n\} )
\end{gather*}
(cf. \cite{AO}, prop. 2.5).\\
Returning to derived categories of coherent sheaves and dropping the hypothesis that the $E_i$'s be objects of the 
underlying Abelian category, we have the following stronger theorem of A. I. Bondal:

\begin{theorem}
Let $X$ be a smooth projective variety and $(E_1, \ldots , E_n)$ a strong complete exceptional sequence in 
$D^b( Coh \, X)$. Set $E:= \bigoplus_{i=1}^n E_i$, let $A:= \mathrm{End} (E) =\bigoplus_{i, j} \mathrm{Hom} (E_i, 
E_j)$ be the algebra of endomorphisms of $E$, and denote $\mathrm{mod}-A$ the category of right modules over $A$ 
which are finite dimensional over $\mathbb{C}$.\\
Then the functor
\[
R\mathrm{Hom}^{\bullet} (E , - ) : D^b( Coh\, (X)) \to D^b( \mathrm{mod}-A)
\] 
is an equivalence of categories (note that, for any object $Y$ of
$D^b(Coh\, (X))$, $R\mathrm{Hom}^{\bullet} (E , Y)$ 
has a natural action from the right by 
$A=\mathrm{Hom} (E , E)$).\\
Moreover, the indecomposable projective modules over $A$ are (up to isomorphism) exactly the $P_i:=
\mathrm{id}_{E_i}\cdot A$, $i=1, \ldots , n$. We have $\mathrm{Hom}_{D^b( Coh \, (X))} (E_i , E_j )$ $ \simeq 
\mathrm{Hom}_A (P_i , P_j)$ and an equivalence 
\[
K^b( \{ P_1, \ldots , P_n\} ) \stackrel{\sim}{\longrightarrow} D^b (\mathrm{mod}- A)
\] 
where $K^b( \{ P_1, \ldots , P_n\} )$ is the homotopy category of complexes of right $A$-modules whose terms are 
finite direct sums of the $P_i$'s. 
\end{theorem}
For a proof see \cite{Bo}, \S\S 5 and 6. Thus whenever we have a strong complete exceptional sequence in 
$D^b (Coh \, (X))$ we get an equivalence of the latter with a homotopy category of projective modules over the 
algebra of endomorphisms of the sequence. For the sequences $\mathfrak{B}$, $\mathfrak{B}'$ in theorem 2.1.1 
we recover Beilinson's theorem (although the objects of the module categories $K^b( \{ P_1, \ldots , P_n\} )$ that 
theorem 2.1.8 produces in each of these cases will be different from the objects in the module categories $K^b_{ 
[0,n]} S$, resp. $K^b_{[0,n]}\Lambda$, in theorem 2.1.1, the morphisms correspond and the respective module categories 
are equivalent).\\
Next suppose that $D^b(Coh\, X)$ on a smooth projective variety $X$ 
is generated by an exceptional sequence $(E_1,\ldots , E_n)$ that is not necessarily strong. Since extension groups 
in nonzero degrees between members of the sequence need not vanish in this case, one cannot expect a description 
of $D^b(Coh\, X)$ on a homotopy category level as in theorem 2.1.8. But still the existence of $(E_1, \ldots, E_n)$ 
makes available some very useful computational tools, e.g. Beilinson type spectral sequences. To state the result, we 
must briefly review some basic material on an operation on exceptional sequences called \emph{mutation}. Mutations 
are also needed in subsection 2.2 below. Moreover, the very concept of exceptional sequence as a weakening of the 
concept of strong exceptional sequence was first introduced because strong exceptionality is in general not 
preserved by mutations, cf. \cite{Bo}, introduction p.24 (exceptional sequences are also more flexible in other 
situations, cf. remark 3.1.3 below). 

For $A,\, B\in \mathrm{obj}\, D^b(Coh\, X)$ set $\mathrm{Hom}^{\times} (A , B):= \bigoplus_{k\in\mathbb{Z}} 
\mathrm{Ext}^k (A, B)$, a graded $\mathbb{C}$-vector space. For a graded $\mathbb{C}$-vector space $V$, $(V^{\vee})^i
:= \mathrm{Hom}_{\mathbb{C}}(V^{-i}, \mathbb{C})$ defines the grading of the dual, and if $X\in \mathrm{obj}\, D^b( 
Coh \, X)$, then $V\otimes X$ means $\bigoplus_{i\in\mathbb{Z}} V^i \otimes X[-i]$ where $V^i\otimes X[-i]$ is the 
direct sum of $\dim V^i$ copies of $X[-i]$.  
\begin{definition}
Let $(E_1 , E_2)$ be an exceptional sequence in $D^b( Coh\, X)$. The \emph{left mutation} $L_{E_1} E_2$ (resp. the 
\emph{right mutation} $R_{E_2} E_1$) is the object defined by the distinguished triangles
\begin{gather*}
L_{E_1} E_2 \longrightarrow \mathrm{Hom}^{\times} (E_1 , E_2) \otimes E_1 \stackrel{\mathrm{can}}{\longrightarrow } 
 E_2 \longrightarrow L_{E_1} E_2 [1] \\
\mathrm{( resp. } \quad R_{E_2} E_1[-1] \longrightarrow E_1 \stackrel{\mathrm{can}'}{\longrightarrow} \mathrm{Hom}^{ 
\times} (E_1 , E_2)^{\vee} \otimes E_2 \longrightarrow R_{E_2} E_1 \quad \mathrm{).}
\end{gather*} 
Here $\mathrm{can}$ resp. $\mathrm{can}'$ are the canonical morphisms (``evaluations''). 
\end{definition}

\begin{theorem}
Let $\mathfrak{E}= (E_1, \ldots , E_n)$ be an exceptional sequence in $D^b(Coh \, X)$. Set, for $i=1, \ldots , n-1$, 
\begin{gather*}
R_i \mathfrak{E} := \left( E_1 , \ldots , E_{i-1}, E_{i+1} , R_{E_{i+1}} E_i , E_{i+2}, \ldots , E_n \right)\, ,  \\
L_i \mathfrak{E} := \left( E_1 , \ldots , E_{i-1}, L_{E_i} E_{i+1} , E_i , E_{i+2}, \ldots , E_n \right)\, .
\end{gather*}
Then $R_i\mathfrak{E}$ and $L_i\mathfrak{E}$ are again exceptional sequences. $R_i$ and $L_i$ are inverse to each 
other; the $R_i$'s (or $L_i$'s) induce an action of $Bd_n$, the Artin
braid group on $n$ strings, on the class of exceptional 
sequences with $n$ terms in $D^b( Coh \, X)$. If moreover $\mathfrak{E}$ is complete, so are all the $R_i 
\mathfrak{E}$'s and $L_i \mathfrak{E}$'s.   
\end{theorem}

For a proof see \cite{Bo}, \S 2.\\
We shall see in example 2.1.13 that the two exceptional sequences 
$\mathfrak{B}$, $\mathfrak{B}'$ of theorem 2.1.1 are closely related
through a notion that we will introduce next:

\begin{definition}
Let $(E_1, \ldots , E_n)$ be a complete exceptional sequence in $D^b(Coh \, X)$. For $i=1, \ldots , n$ define 
\begin{gather*}
E_i^{\vee}: = L_{E_1} L_{E_2} \ldots L_{E_{n-i}} E_{n-i+1} \, ,\\
^{\vee}E_i:= R_{E_n} R_{E_{n-1}} \ldots R_{E_{n-i+2}} E_{n-i+1} \, .
\end{gather*}
The complete exceptional sequences $\left( E_1^{\vee}, \ldots , E_n^{\vee} \right)$ resp. $\left( ^{\vee}E_1 , 
\ldots , ^{\vee}E_n \right)$ are called the \emph{right} resp. \emph{left dual} of $(E_1 ,\ldots , E_n)$. 
\end{definition} 

The name is justified by the following 

\begin{proposition}
Under the hypotheses of definition 2.1.11 one has 
\[
\mathrm{Ext}^k (^{\vee}E_i , E_j)= \mathrm{Ext}^k (E_i , E_j^{\vee}) = \left\{ {\mathbb{C}\quad\mathrm{if}\; i+j=n+1 , 
\: i=k+1 }\atop{0\quad \mathrm{otherwise} \quad\quad\hspace{2cm}}    \right.
\]
Moreover the right (resp. left) dual of $(E_1, \ldots , E_n)$ is uniquely (up to unique isomorphism) defined by these 
equations.
\end{proposition}
 
The proof can be found in \cite{Gor}, subsection 2.6.

\begin{example}
Consider on $\mathbb{P}^n=\mathbb{P}(V)$, the projective space of
lines in the vector space $V$, the complete exceptional sequence $\mathfrak{B}'= (\Omega^n(n), \ldots ,
 \Omega^1(1), \mathcal{O})$ and for $1\le p\le n$ the truncation of the $p$-th exterior power of the Euler sequence 
\[
0\longrightarrow \Omega^p_{\mathbb{P}^n}\longrightarrow \left( \bigwedge\nolimits^p V^{\vee}\right)  \otimes 
\mathcal{O}_{\mathbb{P}^n}(-p)\longrightarrow \Omega^{p-1}_{\mathbb{P}^n} \longrightarrow 0 \, .
\]
Let us replace $\mathfrak{B}'$ by $\left( \Omega^n(n), \ldots , \Omega^2(2), \mathcal{O}, R_{\mathcal{O}}\Omega^1(1) 
\right)$, i.e., mutate $\Omega^1(1)$ to the right across $\mathcal{O}$. But in the exact sequence 
\[
0\longrightarrow \Omega^1(1)\longrightarrow V^{\vee}  \otimes 
\mathcal{O}\longrightarrow \mathcal{O}(1) \longrightarrow 0
\] 
the arrow $\Omega^1(1) \to V^{\vee}\otimes \mathcal{O}$ is nothing but the canonical morphism $\Omega^1(1)\to$ $ 
\mathrm{Hom}(\Omega^1(1), \mathcal{O})^{\vee}\otimes \mathcal{O}$ from definition 2.1.9. Therefore $R_{\mathcal{O} 
}\Omega^1(1)\simeq \mathcal{O}(1)$.\\
Now in the mutated sequence $(\Omega^n(n), \ldots , \Omega^2(2), \mathcal{O}, \mathcal{O}(1))$ we want to mutate in 
the next step $\Omega^2(2)$ across $\mathcal{O}$ and $\mathcal{O}(1)$ to the right. In the sequence 
\[
0\longrightarrow \Omega^2(2) \longrightarrow \bigwedge\nolimits^2 V^{\vee}\otimes \mathcal{O} \longrightarrow 
\Omega^1 (2) \longrightarrow 0
\]  
the arrow $\Omega^2(2) \to \bigwedge^2 V^{\vee} \otimes \mathcal{O}$ is again the canonical morphism $\Omega^2(2) \to$
 $\mathrm{Hom} (\Omega^2(2), \mathcal{O})^{\vee} \otimes\mathcal{O}$ and $R_{\mathcal{O}}\Omega^2(2)\simeq 
\Omega^1(2)$ and then
\[
0\longrightarrow \Omega^1(2)\longrightarrow V^{\vee}  \otimes 
\mathcal{O}(1)\longrightarrow \mathcal{O}(2) \longrightarrow 0
\]
gives $R_{\mathcal{O}(1)} R_{\mathcal{O}} \Omega^2(2) \simeq \mathcal{O}(2)$.\\
Continuing this pattern, one transforms our original sequence $\mathfrak{B}'$ by successive right mutations into 
$(\mathcal{O}, \mathcal{O}(1), \mathcal{O}(2), \ldots ,\mathcal{O}(n))$ which, looking back at definition 2.1.11 and 
using the braid relations $R_i R_{i+1} R_i = R_{i+1} R_i R_{i+1}$, one identifies as the left dual of $\mathfrak
{B}'$. 
\end{example}

Here is Gorodentsev's theorem on generalized Beilinson spectral sequences.

\begin{theorem}
Let $X$ be a smooth projective variety and let $D^b(Coh \, X)$ be generated by an exceptional sequence 
$(E_1, \ldots, E_n)$. Let $F : D^b( Coh \, X)\to \mathcal{A}$ be a covariant cohomological functor to some Abelian 
category $\mathcal{A}$.\\
For any object $A$ in $D^b (Coh\, X)$ there is a spectral sequence
\begin{gather*}
E_1^{p , q} = \bigoplus_{i+j=q} \mathrm{Ext}^{n+i-1} (^{\vee}E_{n-p}, A) \otimes F^j( E_{p+1}) \\
            = \bigoplus_{i+j=q} \mathrm{Ext}^{-i} (A, E_{n-p}^{\vee})^{\vee}  \otimes F^j( E_{p+1}) \; \Longrightarrow 
\; F^{p+q} (A)
\end{gather*} 
(with possibly nonzero entries for $0\le p, q \le n-1$ only).  
\end{theorem}

For the proof see \cite{Gor}, 2.6.4 (actually one can obtain $A$ as a convolution of a complex over $D^b( Coh\, X)$ 
whose terms are computable once one knows the $\mathrm{Ext}^i(^{\vee}E_j , A)$, but we don't need this).\\
In particular, taking in theorem 2.1.14 the dual exceptional sequences in example 2.1.13 and for $F$ the functor 
that takes an object in $D^b(Coh\, \mathbb{P}^n)$ to its zeroth cohomology sheaf, we recover the classical Beilinson 
spectral sequence.

It is occasionally useful to split a derived category into more manageable building blocks before starting 
to look for complete exceptional sequences. This is the motivation for giving the following definitions.

\begin{definition}
Let $\mathcal{S}$ be a full triangulated subcategory of a triangulated category $\mathcal{T}$. The \emph{right 
orthogonal} to $\mathcal{S}$ in $\mathcal{T}$ is the full triangulated subcategory $\mathcal{S}^{\perp}$ of 
$\mathcal{T}$ consisting of objects $T$ such that $\mathrm{Hom}(S, T)=0$ for all objects $S$ of $\mathcal{S}$. 
The \emph{left orthogonal} $^\perp\mathcal{S }$ is defined similarly. 
\end{definition}

\begin{definition}
A full triangulated subcategory $\mathcal{S}$ of $\mathcal{T}$ is \emph{right-} (resp. \emph{left-})
\emph{admissible} if for every $T\in \mathrm{obj}\, \mathcal{T}$ there is a distinguished triangle
\begin{gather*}
S\longrightarrow T \longrightarrow S' \longrightarrow S[1] \quad \mathrm{with}\quad S\in \mathrm{obj}\, 
\mathcal{S}\, ,\; S'\in \mathrm{obj}\, \mathcal{S}^{\perp}\\
\mathrm{( resp.}\quad S''\longrightarrow T \longrightarrow S \longrightarrow S''[1] \quad \mathrm{with}\quad S\in
 \mathrm{obj}\, \mathcal{S}\, ,\; S''\in \mathrm{obj}\, ^\perp\mathcal{S}\, \mathrm{)}
\end{gather*}
and \emph{admissible} if it is both right- and left-admissible.
\end{definition}

Other useful characterizations of admissibility can be found in \cite{Bo}, lemma 3.1 or \cite{BoKa}, prop. 1.5.

\begin{definition}
An $n$-tuple of admissible subcategories $(\mathcal{S}_1, \ldots , \mathcal{S}_n)$ of a triangulated category 
$\mathcal{T}$ is \emph{semi-orthogonal} if $\mathcal{S}_j$ belongs to $\mathcal{S}_i^{\perp}$ whenever 
$1\le j < i\le n$. If $\mathcal{S}_1, \ldots, \mathcal{S}_n$ generate $\mathcal{T}$ one calls this a 
\emph{semi-orthogonal decomposition of} $\mathcal{T}$ and writes
\[
\mathcal{T} =  \left\langle \mathcal{S}_1 ,\ldots , \mathcal{S}_n \right\rangle \, .
\]
\end{definition}

To conclude, we give a result that describes the derived category of coherent sheaves on a product of varieties. 

\begin{proposition}
Let $X$ and $Y$ be smooth, projective varieties and \[(\mathcal{V}_1, \ldots ,\mathcal{V}_m )\] resp. \[
(\mathcal{W}_1,\ldots ,\mathcal{W}_n)\] be (strong) complete exceptional sequences in $D^b( Coh (X))$ resp. 
$D^b(Coh (Y) )$ where $\mathcal{V}_i$ and $\mathcal{W}_j$ are vector
bundles on $X$ resp. $Y$. Let $\pi_1$ resp. 
 $\pi_2$ be the projections of $X\times Y$ on the first resp. second factor and put $\mathcal{V}_i\boxtimes 
\mathcal{W}_j:= \pi_1^{\ast}\mathcal{V}_i \otimes \pi_2^{\ast}\mathcal{W}_j$. Let $\prec $ be the 
lexicographic order on $\{ 1, \ldots , m\} \times \{ 1, \ldots ,n\}$. Then
\[
(\mathcal{V}_i\boxtimes \mathcal{W}_j )_{(i,j)\in\{ 1, \ldots , m\} \times \{ 1, \ldots ,n\} }
\]
is a (strong) complete exceptional sequence in $D^b( Coh (X\times Y))$ where $\mathcal{V}_{i_1}\boxtimes 
\mathcal{W}_{j_1}$ precedes $\mathcal{V}_{i_2}\boxtimes\mathcal{W}_{j_2}$ iff $(i_1 , j_1)\prec (i_2 , j_2)$
. 
\end{proposition}
\begin{proof}
The proof is a little less straightforward than it might be expected at first glance since one does not 
know explicit resolutions of the structure sheaves of the diagonals on $X\times X$ and $Y\times Y$.\\
First, by the K\"unneth formula,
\begin{gather*}
\mathrm{Ext}^k (\mathcal{V}_{i_2}\boxtimes \mathcal{W}_{j_2}, \mathcal{V}_{i_1}\boxtimes \mathcal{W}_{j_1})
\simeq H^k(X\times Y, (\mathcal{V}_{i_1}\otimes \mathcal{V}_{i_2}^{\vee})\boxtimes ( \mathcal{W}_{j_1}
\otimes \mathcal{W}_{j_2}^{\vee} ) ) \\
\simeq \bigoplus\limits_{k_1+k_2=k} H^{k_1} (X , \mathcal{V}_{i_1}\otimes \mathcal{V}_{i_2}^{\vee} ) 
\otimes  H^{k_2} (Y , \mathcal{W}_{j_1}\otimes \mathcal{W}_{j_2}^{\vee} )\\
\simeq \bigoplus\limits_{k_1+k_2=k} \mathrm{Ext}^{k_1} (\mathcal{V}_{i_2},\mathcal{V}_{i_1}) 
\otimes  \mathrm{Ext}^{k_2} (\mathcal{W}_{j_2}, \mathcal{W}_{j_1} )
\end{gather*}
whence it is clear that $(\mathcal{V}_i\boxtimes \mathcal{W}_j )$ will be a (strong) exceptional sequence 
for the ordering $\prec$ if $(\mathcal{V}_i)$ and $(\mathcal{W}_j)$
are so.\\
Therefore we have to show that $(\mathcal{V}_i\boxtimes \mathcal{W}_j )$ generates $D^b( Coh (X\times Y))$ (see 
also \cite{BoBe}, lemma 3.4.1). By \cite{Bo}, thm. 3.2, the triangulated subcategory $\mathcal{T}$ of 
$D^b( Coh (X\times Y))$ generated by the $\mathcal{V}_i\boxtimes \mathcal{W}_j$'s is admissible, and thus by 
\cite{Bo}, lemma 3.1, it suffices to show that the right orthogonal $\mathcal{T}^{\perp}$ is zero. Let $Z\in 
\mathrm{obj}\, \mathcal{T}^{\perp}$ so that we have
\[
\mathrm{Hom}_{D^b( Coh (X\times Y))} (\mathcal{V}_i\boxtimes\mathcal{W}_j , Z) = 0\quad \forall i\in \{ 1, 
\ldots , m\}\, ,\; \forall j\in \{ 1, \ldots , n\}\, .
\]
But
\begin{gather*}
\mathrm{Hom}_{D^b( Coh (X\times Y))} (\mathcal{V}_i\boxtimes\mathcal{W}_j , Z)\hspace{6cm} \\ \simeq 
\mathrm{Hom}_{D^b( Coh (X\times Y))} \left( \pi_1^{\ast}\mathcal{V}_i , R \mathcal{H}om_{D^b( Coh (X\times Y))}^{
\bullet } (\pi_2^{\ast} \mathcal{W}_j , Z ) \right)  \\
\simeq \mathrm{Hom}_{D^b( Coh (X))} \left( \mathcal{V}_i , R\pi_{1 \ast} R \mathcal{H}om_{D^b( Coh (X\times Y))}^
{\bullet } (\pi_2^{\ast} \mathcal{W}_j , Z ) \right) 
\end{gather*}
using the adjointness of $\pi_1^{\ast}=L\pi_1^{\ast}$ and $R\pi_{1 \ast}$. But then 
\[
 R\pi_{1 \ast} R \mathcal{H}om_{D^b( Coh (X\times Y))}^{\bullet } (\pi_2^{\ast} \mathcal{W}_j , Z )=0 \quad 
\forall j\in \{ 1, \ldots , m\} 
\] 
because the $\mathcal{V}_i$ generate $D^b(Coh (X))$ and hence there is no non-zero object in the right 
orthogonal to $\langle \mathcal{V}_1 ,\ldots ,\mathcal{V}_n \rangle$. Let $U\subset X$ and $V\subset Y$ be 
affine open sets. Then
\begin{gather*}
0= R\Gamma \left( U , R\pi_{1 \ast} R \mathcal{H}om_{D^b( Coh (X\times Y))}^{\bullet } (\pi_2^{\ast}
 \mathcal{W}_j , Z ) \right)\\
\simeq R\mathrm{Hom}^{\bullet} \left( \mathcal{W}_j , R\pi_{2 \ast} ( Z \mid_{U\times Y} ) \right)    
\end{gather*}
whence $R\pi_{2 \ast} ( Z \mid_{U\times Y} )=0$ since the $\mathcal{W}_j$ generate $D^b(Coh (Y))$ (though in 
general $R\pi_{2 \ast} ( Z \mid_{U\times Y} )$ will be a complex of quasi-coherent sheaves, one can write it
as the direct limit over its subcomplexes with coherent terms and, using that the direct limit commutes with 
$R^{\bullet} \mathrm{Hom}$, conclude that $R\pi_{2 \ast} ( Z \mid_{U\times Y} )=0$).
 Therefore we get 
\[
R\Gamma (U\times V , Z)=0 \, .
\]
But $R^i\Gamma (U\times V , Z)= \Gamma (U\times V , H^i(Z) )$ and thus all cohomology sheaves of $Z$ are zero, 
i.e. $Z=0$ in $D^b(Coh (X\times Y))$. 
\end{proof}

\begin{remark}
This proposition is very useful for a treatment of the derived categories of coherent sheaves on rational 
homogeneous spaces from a systematic point of view. For if $X=G/P$ with $G$ a connected semisimple complex 
Lie group, $P\subset G$ a parabolic subgroup, it is well known that one has a decomposition
\[
X\simeq S_1/P_1 \times \ldots \times S_N/P_N 
\]
where $S_1 ,\ldots , S_N$ are connected simply connected simple complex Lie groups and $P_1 , \ldots , P_N$ 
corresponding parabolic subgroups (cf. \cite{Akh}, 3.3, p. 74). Thus for the construction of complete exceptional
 sequences on any $G/P$ one can restrict oneself to the case where $G$ is simple.
\end{remark}

\subsection{Catanese's conjecture and the work of Kapranov}

First we fix some notation concerning rational homogenous varieties and their Schubert varieties that will 
remain in force throughout the text unless otherwise stated. References for this are \cite{Se2}, \cite{Sp}.\\
\begin{quote}
$G$ is a complex semi-simple Lie group which is assumed to be connected and simply connected with Lie algebra 
$\mathfrak{g}$.\\
$H\subset G$ is a fixed maximal torus in $G$ with Lie algebra the Cartan subalgebra $\mathfrak{h}\subset 
\mathfrak{g}$.\\
$R\subset \mathfrak{h}^{\ast}$ is the root system associated to $(\mathfrak{g}, \mathfrak{h})$ so that $
\mathfrak{g}=\mathfrak{h}\oplus\bigoplus_{\alpha\in R} \mathfrak{g}^{\alpha}$ with $\mathfrak{g}^{\alpha}$ the 
eigen-subspace of $\mathfrak{g}$ corresponding to $\alpha\in \mathfrak{h}^{\ast}$. Choose a base $S= \{ \alpha_1 ,
\ldots ,\alpha_r\}$ for $R$; $R^{+}$ denotes the set of positive roots w.r.t. $S$, $R^{-}:=-R^{+}$, so that 
$R=R^{+}\cup R^{-}$, and $\varrho$ is the half-sum of the positive roots.\\
$\mathrm{Aut}(\mathfrak{h}^{\ast})\supset W:=\langle s_{\alpha}\, \mid\, s_{\alpha}$ the reflection with vector 
$\alpha$ leaving $R$ invariant$\rangle \simeq N(H)/H$ is the Weyl group of $R$.\\
Let $\mathfrak{b}:=\mathfrak{h}\oplus \bigoplus_{\alpha >0}\mathfrak{g}^{\alpha }$,  $\mathfrak{b}^{-}:=
\mathfrak{h}\oplus \bigoplus_{\alpha < 0}\mathfrak{g}^{\alpha }$ be opposite Borel subalgebras of $\mathfrak{g}$ 
corresponding to $\mathfrak{h}$ and $S$, and $\mathfrak{p}\supset \mathfrak{b}$ a parabolic subalgebra 
corresponding uniquely to a subset $I\subset S$ (then $\mathfrak{p}=\mathfrak{p}(I)= \mathfrak{h}\oplus 
\bigoplus_{\alpha\in R^{+}}\mathfrak{g}^{\alpha}\oplus \bigoplus_{\alpha \in R^{-}(I)} \mathfrak{g}^{\alpha}$ 
where $R^{-}(I):= \{ \alpha \in R^{-}\, \mid \, \alpha = \sum_{i=1}^r k_i \alpha_i$ with $k_i\le 0$ for all $i$ and 
$k_j=0$ for all $\alpha_j\in I\}$). Let $B^{+}$, $B^{-}$, $P=P(I)\supset B$ be the corresponding connected 
subgroups of
$G$  with Lie algebras $\mathfrak{b}$, $\mathfrak{b}^{-}$, $\mathfrak{p}$.\\
$X:=G/P$ is the rational homogeneous variety corresponding to $G$ and $P$.\\
$l(w)$ is the length of an element $w\in W$ relative to the set of generators $\{ s_{\alpha}\, \mid \, \alpha\in S
\}$, i.e. the least number of factors in a decomposition 
\[
w=s_{\alpha_{i_1}} s_{\alpha_{i_2}}\ldots s_{\alpha_{i_l}}\, , \; \alpha_{i_j}\in S\, ;
\]
A decomposition with $l=l(w)$ is called reduced. One has the Bruhat order $\le $ on $W$, i.e. $x\le w$ for 
$x, \, w\in W$ iff $x$ can be obtained by erasing some factors of a reduced decomposition of $w$.\\
$W_P$ is the Weyl group of $P$, the subgroup of $W$ generated by the simple reflections $s_{\alpha}$ with $\alpha 
\notin I$. In each coset $wW_P\in W/W_P$ there exists a unique element of minimal length and $W^P$ denotes the 
set of minimal representatives of $W/W^P$. One has $W^P=\{ w\in W\, \mid \, l(ww')=l(w)+l(w') \; \forall w'\in 
W_P\}$.\\
For $w\in W^P$, $C_w$ denotes the double coset $BwP/P$ in $X$, called a Bruhat cell, $C_w\simeq \mathbb{A}^{l(w)}$. 
Its closure in $X$ is the Schubert variety $X_w$. $C_w^{-}=B^{-}wP/P$ is the opposite Bruhat cell of codimension 
$l(w)$ in $X$, $X^w=\overline{{C}_w^{-}}$ is the Schubert variety opposite to $X_w$. \\ There is the extended version
 of the Bruhat decomposition
\[
G/P= \bigsqcup_{w\in W^P} C_{w}
\]  
(a paving of $X$ by affine spaces) and for $v,\; w\in W^P$: $v\le w \Leftrightarrow X_v\subseteq X_w$; we denote 
the boundaries $\partial X_w:= X_w \backslash C_{w}$, $\partial X^{w}: =X^{w}\backslash C_w^{-}$, which have pure 
codimension $1$ in $X_w$ resp. $X^w$. 
\end{quote}
Moreover, we need to recall some facts and introduce further notation
concerning representations of the subgroup $P=P(I)\subset G$, which
will be needed in subsection 3 below. References are \cite{A}, \cite{Se2}, \cite{Sp}, \cite{Ot2},
\cite{Stei}.
\begin{quote}
The spaces $\mathfrak{h}_{\alpha}:=[ \mathfrak{g}^{\alpha},
\mathfrak{g}^{-\alpha}]\subset \mathfrak{h}$, $\alpha\in R$, are
1-dimensional, one has $\mathfrak{g}=\bigoplus_{\alpha\in
  S}\mathfrak{h}_{\alpha}\oplus \bigoplus_{\alpha\in
  R^{+}}\mathfrak{g}^{\alpha}\oplus \bigoplus_{\alpha\in
  R^{-}}\mathfrak{g}^{\alpha}$ and there is a unique
$H_{\alpha}\in\mathfrak{h}_{\alpha}$ such that $\alpha (H_{\alpha})=2$.
\\
Then we have the weight lattice $\Lambda := \{
\omega\in\mathfrak{h}^{\ast}\, \mid \, \omega
(H_{\alpha})\in\mathbb{Z}\; \forall \alpha\in R\}$ (which one
identifies with the character group of $H$) and the set of
dominant weights $\Lambda^{+}:= \{\omega\in\mathfrak{h}^{\ast}\, \mid \, \omega
(H_{\alpha})\in\mathbb{N}\; \forall \alpha\in R\}$. $\{\omega_1 ,\dots
, \omega_r \}$ denotes the basis of $\mathfrak{h}^{\ast}$ dual to the
basis $\{ H_{\alpha_1},\dots , H_{\alpha_r} \}$ of $\mathfrak{h}$. The
$\omega_i$ are the fundamental weights. If $(\cdot ,\cdot )$ is the
inner product on $\mathfrak{h}^{\ast}$ induced by the Killing form,
they can also be characterized by the equations $2(\omega_i ,
\alpha_j) /(\alpha_j , \alpha_j) =\delta_{ij}$ (Kronecker delta). It
is well known that the irreducible finite dimensional representations
of $\mathfrak{g}$ are in one-to-one correspondence with the
$\omega\in\Lambda^{+}$, these $\omega$ occurring as highest weights.\\
I recall the Levi-Mal\v{c}ev decomposition of $P(I)$
(resp. $\mathfrak{p}(I)$): The algebras 
\[
\mathfrak{s}_{P}:=\bigoplus_{\alpha\in S\backslash
  I}\mathfrak{h}_{\alpha}\oplus \bigoplus_{\alpha\in
  R^{-}(I)}(\mathfrak{g}^{\alpha}\oplus \mathfrak{g}^{-\alpha})
\]  
resp.
\[
\mathfrak{l}_{P}:=\bigoplus_{\alpha\in S}\mathfrak{h}_{\alpha}\oplus \bigoplus_{\alpha\in
  R^{-}(I)}(\mathfrak{g}^{\alpha}\oplus \mathfrak{g}^{-\alpha})
\]
are the semisimple resp. reductive parts of $\mathfrak{p}(I)$
containing $\mathfrak{h}$, the
corresponding connected subgroups of $G$ will be denoted $S_P$ resp. $L_P$. The
algebra
\[
\mathfrak{u}_P:= \bigoplus_{\alpha\in R^{-}\backslash R^{-}(I)} \mathfrak{g}^{-\alpha}
\] 
is an ideal of $\mathfrak{p}(I)$,
$\mathfrak{p}(I)=\mathfrak{l}_P\oplus \mathfrak{u}_P$, and the
corresponding normal subgroup $R_u (P)$ is the unipotent radical of
$P$. One has 
\[
P=L_P\ltimes R_u(P)\, ,
\]
the Levi-Mal\v{c}ev decomposition of $P$. The center $Z$ of the Levi
subgroup $L_P$ is $Z=\{ g\in H\, |\, \alpha (g)=1 \: \forall \alpha\in
S\backslash I\}$. It corresponds to the Lie algebra
$\bigoplus_{\alpha\in I}\mathfrak{h}_{\alpha}$ and is isomorphic to
the torus $\left( \mathbb{C}^{\ast}\right)^{|I|}$. One has 
\[
P=Z \cdot S_P \ltimes R_u(P)\, .
\]
Under the hypothesis that $G$ is simply connected, also $S_P$ is
simply connected. \\
If $r: P \to \mathrm{GL} (V)$ is an irreducible finite-dimensional
representation, $R_u(P)$ acts trivially, and thus those $r$ are in
one-to-one correspondence with irreducible representations of the
reductive 
Levi-subgroup $L_P$ and as such possess a well-defined highest weight
$\omega\in\Lambda$.
 Then the irreducible finite
dimensional representations of $P(I)$ correspond bijectively to
weights $\omega\in \mathfrak{h}^{\ast}$ such that $\omega$ can be
written as $\omega =\sum_{i=1}^{r} k_i \omega_i$, $k_i\in\mathbb{Z}$,
such that $k_j\in\mathbb{N}$ for all $j$ such that $\alpha_j\notin
I$. We will say that such an $\omega$ is the highest weight of the
representation $r: P\to \mathrm{GL}(V)$. \\
The homogeneous vector bundle on $G/P$ associated to $r$ will be
$G\times_r V:= G\times V/\{(g,v)\sim (gp^{-1} , r (p) v )\, ,\; p\in
P,\: g\in G ,\: v\in V \}$ as above. However, for a character $\chi :
H\to \mathbb{C}$ (which will often be identified with $d\chi
\in\mathfrak{h}^{\ast}$), $\mathcal{L}(\chi )$ will denote the
homogeneous line bundle on $G/B$ whose fibre at the point $e\cdot B$
is the one-dimensional representation of $B$ corresponding to the
character $-\chi$. This has the advantage that $\mathcal{L}(\chi )$
will be ample iff $d\chi =\sum_{j=1}^r k_j \omega_j$ with $k_j > 0$,
$k_j\in\mathbb{Z}$ for all $j$, and it will also prove a reasonable
convention in later applications of Bott's theorem. 
\end{quote}

The initial stimulus for this work was a conjecture due to
F. Catanese. This is variant (A) of conjecture 2.2.1. Variant (B) is a
modification of (A) due to the author, but closely related.  

\begin{conjecture}
\begin{itemize}
\item[(A)]
On any rational homogeneous variety $X=G/P$ there exists a complete
strong exceptional poset (cf. def. 2.1.7 (B)) and a bijection of the
elements of the poset with the Schubert varieties in $X$ such that the
partial order of the poset is the one induced by the Bruhat-Chevalley
order.   
\item[(B)]
For any $X=G/P$ there exists a strong complete exceptional sequence $\mathfrak{E}=(E_1, \ldots , E_n)$ in 
$D^b(Coh \, X)$ with $n=| W^P|$, the number of Schubert varieties in
$X$ (which is the topological Euler characteristic of $X$).\\
Moreover, since there is a natural partial order $\le_{\mathfrak{E}}$ on the set of objects in $\mathfrak{E}$ by 
defining that $E' \le_{\mathfrak{E}} E$ for objects $E$ and $E'$ of $\mathfrak{E}$ iff there are objects $F_1, \ldots 
, F_r$ of $\mathfrak{E}$ such that $\mathrm{Hom}( E' , F_1)\neq 0$, $\mathrm{Hom}( F_1 , F_2)\neq 0$, $\ldots$, 
$\mathrm{Hom}( F_r , E)\neq 0$ (the order of the exceptional sequence $\mathfrak{E}$ itself is a total order 
refining $\le_{\mathfrak{E}}$), there should be a relation between the Bruhat order on $W^P$ and $\le_{\mathfrak{E}}$ 
(for special choice of $\mathfrak{E}$).\\
If $P=P(\alpha_i)$, some $i\in\{ 1,\ldots , r\}$, is a maximal parabolic subgroup in $G$ and $G$ is simple, then one 
may conjecture more precisely: There exists a strong complete exceptional sequence $\mathfrak{E}=(E_1, \ldots , E_n)$
 in $D^b(Coh \, X)$ and a bijection 
\[
b : \{ E_1 ,\ldots , E_n\} \to \{ X_w\, |\, w\in W^P\}
\]
such that 
\[
\mathrm{Hom}(E_i , E_j)\neq 0 \: \iff \: b(E_j)\subseteq b(E_i)\, .
\] 
\end{itemize} 
\end{conjecture}
We would like to add the following two questions:
\begin{itemize}
\item[(C)]
Does there always exist on $X$ a complete very strong exceptional
poset (cf. def. 2.1.7 (C)) and a bijection of the
elements of the poset with the Schubert varieties in $X$ such that the
partial order of the poset is the one induced by the Bruhat-Chevalley
order?
\item[(D)]
Can we achieve that the $E_i$'s in (A), (B) and/or (C) are homogeneous
vector bundles?
\end{itemize}
It is clear that, if the answer to (C) is positive, this implies
(A). Moreover, the existence of a complete very strong exceptional
poset entails the existence of a complete strong exceptional sequence.\\ 
For $P$ maximal parabolic, part (B) of conjecture 2.2.1 is stronger
than part (A). We will concentrate on that case in the following.\\
In the next subsection we will see that, at least upon adopting the right point of view, it is clear that the number of 
terms in any complete exceptional sequence in $D^b(Coh\, X)$ must equal the number of Schubert varieties in $X$.\\
To begin with, let me show how conjecture 2.2.1 can be brought in line with results of Kapranov obtained in 
\cite{Ka3} (and \cite{Ka1}, \cite{Ka2}) which are summarized in theorems 2.2.2, 2.2.3, 2.2.4 below.\\
One more piece of notation: If $W$ is an $m$-dimensional vector space and $\lambda = (\lambda_1 ,\ldots ,\lambda_m )$ is 
a non-increasing sequence of integers, then $\Sigma^{\lambda} W$ will denote the space of the irreducible representation 
$\varrho_{\lambda} : \mathrm{GL} (W) \to \mathrm{Aut} (\Sigma^{\lambda} W )$ of $\mathrm{GL} (W)\simeq \mathrm{GL}_m
\mathbb{C}$ with highest weight 
$\lambda$. $\Sigma^{\lambda}$ is called the Schur functor associated to $\lambda$; if $\mathcal{E}$ is a rank $m$ vector 
bundle on a variety $Y$, $\Sigma^{\lambda} \mathcal{E}$ will denote the vector bundle $P_{\mathrm{GL}} 
(\mathcal{E}) \times_{\varrho_{\lambda}} \Sigma^{\lambda} (W):= P_{\mathrm{GL}} (\mathcal{E}) \times \Sigma^{\lambda} (W)/
 \{ (f,w) \sim (fg^{-1} , \varrho_{\lambda}(g) w) , \; f\in P_{\mathrm{GL}} (\mathcal{E}),\; w \in \Sigma^{\lambda} W ,\;
 g\in \mathrm{GL}_m \mathbb{C} \}$ where $P_{\mathrm{GL}} (\mathcal{E})$ is the principal $\mathrm{GL}_m \mathbb{C}$-
bundle of local frames in $\mathcal{E}$. 

\begin{theorem}
Let $\mathrm{Grass} (k , V)$ be the Grassmanian of $k$-dimensional subspaces of an $n$-dimensional vector space $V$, 
and let $\mathcal{R}$ be the tautological rank $k$ subbundle on $\mathrm{Grass} (k , V)$. Then the bundles $\Sigma^{
\lambda} \mathcal{R}$ where $\lambda$ runs over $Y(k , n-k)$, the set of Young diagrams with no more than $k$ rows and 
no more than $n-k$ columns, are all exceptional, have no higher extension groups between each other and generate 
$D^b (Coh\, \mathrm{Grass} (k , V))$.\\
Moreover, $\mathrm{Hom} (\Sigma^{\lambda} \mathcal{R} , \Sigma^{\mu} \mathcal{R}) \neq 0$ iff $\lambda_i \ge \mu_i\; 
\forall i=1, \dots ,k$. (Thus these $\Sigma^{\lambda} \mathcal{R}$ form a strong complete exceptional sequence in 
$D^b (Coh\, \mathrm{Grass} (k , V))$ when appropriately ordered).   
\end{theorem}  

\begin{theorem}
If $V$ is an $n$-dimensional vector space, $1\le k_1 < \dots < k_l \le n$ a strictly increasing sequence of integers, 
and $\mathrm{Flag} ( k_1 , \dots ,  k_l ; V)$ the variety of flags of subspaces of type $(k_1 ,\dots , k_l)$ in $V$, 
and if $\mathcal{R}_{k_1} \subset \dots \subset \mathcal{R}_{k_l}$ denotes the tautological flag of subbundles, then 
the bundles 
\[
\Sigma^{\lambda_1}\mathcal{R}_{k_1}\otimes \dots \otimes \Sigma^{\lambda_l} \mathcal{R}_{k_l}
\]
where $\lambda_j$, $j=1,\ldots , l-1$, runs over $Y(k_j , k_{j+1}-k_j)$, the set of Young diagrams with no more than $k_j$
 rows and no more than $k_{j+1}-k_j$ columns, and $\lambda_l$ runs over $Y(k_l , n-k_l)$, form a strong complete 
exceptional sequence in $D^b (Coh \, \mathrm{Flag} ( k_1 , \dots ,  k_l ; V)$ if we order them as follows:\\
Choose a total order $\prec_j$ on each of the sets $Y(k_j , k_{j+1}-k_j)$ and $\prec_l$ on $Y(k_l, n-k_l)$ such that 
if $\lambda \prec_j \mu$ (or $\lambda \prec_t \mu$) then the Young diagram of $\lambda$ is not contained in the 
Young diagram of $\mu$; endow the set $Y= Y(k_l , n-k_l) \times Y(k_{l-1} , k_l- k_{l-1}) \times \dots \times 
Y(k_1 ,k_2-k_1)$ with the resulting lexicographic order $\prec$. Then $\Sigma^{\lambda_1}\mathcal{R}_{k_1}\otimes \dots
 \otimes \Sigma^{\lambda_l} \mathcal{R}_{k_l}$ precedes $\Sigma^{\mu_1}\mathcal{R}_{k_1}\otimes \dots \otimes 
\Sigma^{\mu_l} \mathcal{R}_{k_l}$ iff $(\lambda_l, \dots ,\lambda_1) \prec (\mu_l , \dots , \mu_1)$. 
\end{theorem}

\begin{theorem}
Let $V$ be again an $n$-dimensional vector space and $Q\subset \mathbb{P}(V)$ a nonsingular quadric hypersurface. \\
If $n$ is odd and $\Sigma$ denotes the spinor bundle on $Q$, then the following constitutes a strong complete 
exceptional sequence in $D^b(Coh\, Q)$: 
\[
(\Sigma(-n+2),\: \mathcal{O}_Q(-n+3),\dots , \mathcal{O}_Q(-1), \: \mathcal{O}_Q)
\]
and $\mathrm{Hom}(\mathcal{E}, \mathcal{E}')\neq 0$ for two bundles $\mathcal{E}$, $\mathcal{E}'$ in this sequence 
iff $\mathcal{E}$ precedes $\mathcal{E}'$ in the ordering of the sequence.\\
If $n$ is even and $\Sigma^{+}$, $\Sigma^{-}$ denote the spinor bundles on $Q$, then 
\[
(\Sigma^{+}(-n+2),\: \Sigma^{-}(-n+2),\: \mathcal{O}_Q(-n+3),\dots , \mathcal{O}_Q(-1), \: \mathcal{O}_Q)
\]
is a strong complete exceptional sequence in $D^b(Coh\, Q)$ and $\mathrm{Hom}(\mathcal{E}, \mathcal{E}')\neq 0$ 
for two bundles $\mathcal{E}$, $\mathcal{E}'$ in this sequence 
iff $\mathcal{E}$ precedes $\mathcal{E}'$ in the ordering of the sequence with the one exception that 
$\mathrm{Hom} (\Sigma^{+}(-n+2), \Sigma^{-}(-n+2))=0$. 
\end{theorem}

Here by $\Sigma$ (resp. $\Sigma^{+},\:\Sigma^{-}$), we mean the homogeneous vector bundles on $Q=\mathrm{Spin}_n\mathbb{C}
/ P(\alpha_1)$, $\alpha_1$ the simple root corresponding to the first node in the Dynkin diagram of type $B_m$, $
n=2m+1$, (resp. the Dynkin diagram of type $D_m$, $n=2m$), that are
the duals of the vector bundles associated to the irreducible
representation of $P(\alpha_1)$ with highest weight $\omega_m$
(resp. highest weights $\omega_m$, $\omega_{m-1}$). We will deal more extensively with spinor bundles 
in subsection 3.2 below.\\
First look at theorem 2.2.2. It is well known (cf. \cite{BiLa}, section 3.1) that if one sets
\[
I_{k,n}:= \{ \underline{i}=(i_1,\ldots , i_k)\in \mathbb{N}^k \, |\, 1\le i_1 < \dots < i_k \le n \}
\]
and if $V_i:= \langle v_1, \dots , v_i \rangle$ where $(v_1, \dots , v_n)$ is a basis for $V$, then the Schubert 
varieties in $\mathrm{Grass}(k , V)$ can be identified with 
\[
X_{\underline{i}}:= \left\{ L\in \mathrm{Grass}(k , V)\, |\, \dim (L\cap V_{i_j})\ge j \quad\forall 1\le j\le k\right\} 
, \quad \underline{i}\in I_{k, n} 
\]
and the Bruhat order is reflected by 
\[
X_{\underline{i}} \subseteq X_{\underline{i'}} \iff i_j \le i'_j \quad \forall 1\le j\le k ;
\]
and the $\underline{i}\in I_{k, n}$ bijectively correspond to Young
diagrams in $Y(k , n-k)$ by associating to $\underline
{i}$ the Young diagram $\lambda (\underline{i})$ defined by 
\[
\lambda (\underline{i})_t:= i_{k-t+1}- (k-t+1) \quad \forall 1\le t\le k \, . 
\] 
Then containment of Schubert varieties corresponds to containment of
associated Young diagrams. Thus conjecture  2.2.1 (B) 
is verified by the strong complete exceptional sequence of theorem 2.2.2.\\
In the case of $\mathrm{Flag}(k_1,\ldots ,k_l; V)$ (theorem 2.2.3) one can describe the Schubert subvarieties 
and the Bruhat order as follows (cf. \cite{BiLa}, section 3.2): Define 
\[
I_{k_1, \dots , k_l}:= \left\{ \left(\underline{i}^{(1)}, \dots , \underline{i}^{(l)}\right)\in I_{k_1 , n}\times 
\dots \times I_{k_l, n} \, | \, \underline{i}^{(j)}\subset \underline{i}^{(j+1)}\quad \forall 1\le j\le l-1 \right\} .  
\]
Then the Schubert varieties in $\mathrm{Flag}(k_1,\ldots ,k_l; V)$ can be identified with the 
\begin{gather*}
X_{\left(\underline{i}^{(1)}, \dots , \underline{i}^{(l)}\right)}:= 
\left\{ \left( L_1, \ldots , L_l  \right)\in \mathrm{Flag}(k_1,\ldots ,k_l; V)\subset \mathrm{Grass}(k_1, V)\times 
\ldots \right. \\ \left. \ldots \times \mathrm{Grass}(k_l, V)\, | \, 
L_j \in X_{\underline{i}^{(j)}} \quad \forall 1\le j\le l 
\right\} 
\end{gather*}
for $\left(\underline{i}^{(1)}, \dots , \underline{i}^{(l)}\right)$ running over $I_{k_1, \dots , k_l}$ (keeping the 
preceding notation for the Grassmannian). The Bruhat order on the
Schubert varieties may be identified with the following 
partial order on $I_{k_1, \dots , k_l}$: 
\[
\left(\underline{i}^{(1)}, \dots , \underline{i}^{(l)}\right) \le \left(\underline{j}^{(1)}, \dots ,
 \underline{j}^{(l)}\right) \iff i^{(t)}\le j^{(t)}\quad \forall 1\le t\le l .
\]
To set up a natural bijection between the set $Y$ in theorem 2.2.3 and $I_{k_1, \dots , k_l}$ associate to 
$\mathbf{i}:=\left(\underline{i}^{(1)}, \dots , \underline{i}^{(l)}\right)$ the following Young diagrams: 
$\lambda_l(\mathbf{i}) \in Y(k_l , n-k_l)$ is defined by
\[
(\lambda_l(\mathbf{i}))_t:=i^{(l)}_{k_l-t+1}-(k_l-t+1)\quad \forall 1\le t\le k_l .
\]
Now since $\underline{i}^{(j)}\subset \underline{i}^{(j+1)}\quad \forall 1\le j\le l-1$ one can write
\[
i^{(j)}_s= i^{(j+1)}_{r(s)} ,\; s=1, \ldots , k_j
\]
where $1\le r(1) < \ldots < r(k_j) \le k_{j+1}$. One then defines $\lambda_{j}(\mathbf{i})\in Y(k_j, k_{j+1}-k_j)$ by 
\[
(\lambda_j(\mathbf{i}))_t:=r(k_j-t+1)-(k_j-t+1)\quad \forall 1\le t\le k_j .
\]
However it is not clear to me in this case how to relate the Bruhat
order on $I_{k_1, \dots , k_l}$ with the vanishing or non-vanishing of
$\mathrm{Hom}$-spaces between members of the strong complete
exceptional sequence in theorem 2.2.3 (there is an explicit combinatorial criterion for the non-vanishing of 
\[
\mathrm{Hom}\left( \Sigma^{\lambda_1}\mathcal{R}_{k_1}\otimes \dots \otimes \Sigma^{\lambda_l} \mathcal{R}_{k_l}
, \Sigma^{\mu_1}\mathcal{R}_{k_1}\otimes \dots \otimes \Sigma^{\mu_l} \mathcal{R}_{k_l} \right)
\]
formulated in \cite{Ka3}, 3.12, but if this relates in any perspicuous way to the Bruhat order is not clear). In this 
respect, for the time being, conjecture 2.2.1 (parts (A) and (B)) must remain within the confines of wishful thinking. 
\\
If in the set-up of theorem 2.2.4 $Q\subset \mathbb{P}(V)$, $\dim V=n=2m+1$ odd, is a smooth quadric hypersurface, then 
there are $2m$ Schubert varieties in $Q$ and the Bruhat order on them is linear (cf. \cite{BiLa}, pp. 139/140), so the 
strong complete exceptional sequence of theorem 2.2.4 satisfies
conjecture 2.2.1. (B).\\
The case of a smooth quadric hypersurface $Q\subset \mathbb{P}(V)$ with $\dim V= n= 2m$ even, is more interesting. The 
Bruhat order on the set of Schubert varieties can be depicted in the following way (cf. \cite{BiLa}, p. 142/143): 
\vspace{0.2cm}
\setlength{\unitlength}{1cm}
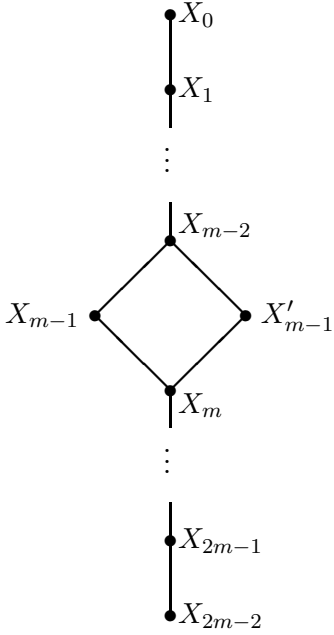
\begin{center}
\begin{picture}(1,8)
\put(0,0){$\bullet$}
\put(0,1){$\bullet$}
\put(0,2){$\vdots$}
\put(0,3){$\bullet$}
\put(-1,4){$\bullet$}
\put(1,4){$\bullet$}
\put(0,5){$\bullet$}
\put(0,6){$\vdots$}
\put(0,7){$\bullet$}
\put(0,8){$\bullet$}

\thicklines
\put(0.1,0.1){\line(0,1){1}}
\put(0.1,1.1){\line(0,1){0.5}}
\put(0.1,3.1){\line(0,-1){0.5}}

\put(0.1,3.1){\line(1,1){1}}
\put(0.1,3.1){\line(-1,1){1}}
\put(-0.9,4.1){\line(1,1){1}}
\put(1.1,4.1){\line(-1,1){1}}

\put(0.1,5.1){\line(0,1){0.5}}
\put(0.1,7.1){\line(0,-1){0.5}}
\put(0.1,7.1){\line(0,1){1}}
\thinlines

\put(0.2,0){$X_{2m-2}$}
\put(0.2,1){$X_{2m-1}$}
\put(0.2,2.8){$X_m$}
\put(-2.1,4){$X_{m-1}$}
\put(1.3,4){$X_{m-1}'$}
\put(0.2,5.2){$X_{m-2}$}
\put(0.2,7){$X_1$}
\put(0.2,8){$X_0$}
\end{picture}
\end{center}
\vspace{0.3cm}
Here $X_0,\dots , X_{m-2},\: X_{m-1},\: X_{m-1}', \: X_m , \ldots , X_{2m-2}$ are labels for the Schubert varieties in 
$Q$ and the subscript denotes the codimension in $Q$. The strong
complete exceptional sequence 
\[
(\Sigma^{+}(-2m+2),\: \Sigma^{-}(-2m+2),\: \mathcal{O}_Q(-2m+3),\dots , \mathcal{O}_Q(-1), \: \mathcal{O}_Q)
\] 
does not verify conjecture 2.2.1 (B), but we claim that there is a strong complete exceptional sequence in the same 
braid group orbit (see thm. 2.1.10) that does. In fact, by \cite{Ott}, theorem 2.8, there are two natural exact 
sequences on $Q$
\begin{gather*}
0 \longrightarrow \Sigma^{+} (-1)\longrightarrow \mathrm{Hom}( \Sigma^{+}(-1),
\mathcal{O}_Q)^{\vee}\otimes \mathcal{O}_Q
 \longrightarrow \Sigma^{-} 
\longrightarrow 0 \\
0 \longrightarrow \Sigma^{-} (-1) \longrightarrow \mathrm{Hom}( \Sigma^{-}(-1),
\mathcal{O}_Q)^{\vee}\otimes \mathcal{O}_Q
 \longrightarrow \Sigma^{+} 
\longrightarrow 0
\end{gather*}  
where the (injective) arrows are the canonical morphisms 
of definition 2.1.9; one also has $\dim \mathrm{Hom}( \Sigma^{+}(-1),
\mathcal{O}_Q)^{\vee}=\dim \mathrm{Hom}( \Sigma^{-}(-1),
\mathcal{O}_Q)^{\vee}=2^{m-1}$. 
(Caution: the spinor bundles in \cite{Ott} are the duals of the bundles that are called spinor 
bundles in this text which is clear from the discussion in \cite{Ott}, p.305!). It follows that if in the above 
strong complete exceptional sequence we mutate $\Sigma^{-}(-2m+2)$ across $\mathcal{O}_Q(-2m+3), \dots , \mathcal{O}_Q
(-m+1)$ to the right and afterwards mutate $\Sigma^{+}(-2m+2)$ across $\mathcal{O}_Q(-2m+3), \dots , \mathcal{O}_Q
(-m+1)$ to the right, we will obtain the following complete exceptional sequences in $D^b(Coh \, Q)$:
\begin{gather*}
\intertext{If $m$ is odd:}
\left( \mathcal{O}_Q(-2m+3), \dots , \mathcal{O}_Q(-m+1), \: \Sigma^{+}(-m+1), \: \Sigma^{-} (-m+1),\right.\\
\left. \: \mathcal{O}_Q(-m+2), 
\ldots , \mathcal{O}_Q(-1), \: \mathcal{O}_Q \right) ,\\
\intertext{if $m$ is even:}
\left( \mathcal{O}_Q(-2m+3), \dots , \mathcal{O}_Q(-m+1), \: \Sigma^{-}(-m+1), \: \Sigma^{+} (-m+1),\right. \\
\left. \: \mathcal{O}_Q(-m+2), 
\ldots , \mathcal{O}_Q(-1), \: \mathcal{O}_Q \right) .
\end{gather*}
One finds (e.g. using theorem 2.2.4 and \cite{Ott}, thm.2.3 and thm. 2.8) that these exceptional sequences are again 
strong and if we let the bundles occurring in them (in the order given by the sequences) correspond to 
$X_0,\dots , X_{m-2},\: X_{m-1},\: X_{m-1}', \: X_m , \ldots , X_{2m-2}$ (in this order), then the above two strong 
complete exceptional sequences verify conjecture 2.2.1. (B).
                                           
\subsection{Information detected on the level of K-theory}

The cellular decomposition of $X$ has the following impact on $D^b(Coh\, X)$.

\begin{proposition}
The structure sheaves $\mathcal{O}_{X_w}$, $w\in W^P$, of Schubert
varieties in $X$ generate $D^b(Coh\, X)$ as a  
triangulated category.
\end{proposition}

Since we have the Bruhat decomposition and each Bruhat cell is isomorphic to an affine space, the proof of the 
proposition will follow from the next lemma.

\begin{lemma}
Let $Y$ be a reduced algebraic scheme, $U\subset Y$ an open subscheme
with $U\simeq \mathbb{A}^d$, for some $d\in \mathbb{N}$, $Z:=Y 
\backslash U$, $i: U\hookrightarrow Y$, $j: Z\hookrightarrow Y$ the natural embeddings. Look at the sequence of 
triangulated categories and functors
\[
\begin{CD}
D^b(Coh\, Z) @>{j_{\ast}}>> D^b(Coh\, Y) @>{i^{\ast }}>> D^b (Coh\, U) 
\end{CD}
\]
(thus $j_{\ast}$ is extension by $0$ outside $Z$ which is exact, and $i^{\ast}$ is the restriction to $U$, 
likewise exact). Suppose $Z_1,\ldots , Z_n\in \mathrm{obj}\, D^b(Coh \, Z)$ generate $D^b(Coh\, Z)$.\\
Then $D^b (Coh\, Y)$ is generated by $j_{\ast} Z_1, \ldots , j_{\ast} Z_n,\, \mathcal{O}_Y$. 
\end{lemma}

\begin{proof}
$D^b(Coh\, Y)$ is generated by $Coh\, Y$ so it suffices to prove that each coherent sheaf $\mathcal{F}$ on $Y$ 
is isomorphic to an object in the triangulated subcategory generated by $j_{\ast} Z_1, \ldots , j_{\ast} Z_n,\, 
 \mathcal{O}_Y$. By the Hilbert syzygy theorem $i^{\ast}\mathcal{F}$ has a resolution
\[
(\ast )\quad 0\to \mathcal{L}_t \to \ldots \to \mathcal{L}_0 \to i^{\ast}\mathcal{F} \to 0
\]
where the $\mathcal{L}_i$ are finite direct sums of $\mathcal{O}_U$. We recall the following facts (cf. \cite{FuLa}
, VI, lemmas 3.5, 3.6, 3.7):
\begin{itemize}
\item[(1)]
For any coherent sheaf $\mathcal{G}$ on $U$ there is a coherent extension $\overline{\mathcal{G}}$ to $Y$.
\item[(2)]
Any short exact sequence of coherent sheaves on $U$ is the restriction of an exact sequence of coherent sheaves 
on $Y$. 
\item[(3)]
If $\mathcal{G}$ is coherent on $U$ and $\overline{\mathcal{G}}_1$, $\overline{\mathcal{G}}_2$ are two coherent 
extensions of $\mathcal{G}$ to $Y$, then there are a coherent sheaf $\overline{\mathcal{G}}$ on $Y$ and 
homomorphisms $\overline{\mathcal{G}}\stackrel{f}{\longrightarrow} \overline{\mathcal{G}}_1$, 
$\overline{\mathcal{G}}\stackrel{g}{\longrightarrow} 
\overline{\mathcal{G}}_2$ which restrict to isomorphisms over $U$.
\end{itemize}
Note that in the set-up of the last item we can write 
\begin{gather*}
0\to \mathrm{ker}( f) \to \overline{\mathcal{G}}\stackrel{f}{\longrightarrow} \overline{\mathcal{G}}_1 \to 
\mathrm{coker}(f) \to 0 \, ,\\ 
0\to \mathrm{ker}( g) \to \overline{\mathcal{G}}\stackrel{g}{\longrightarrow} \overline{\mathcal{G}}_2 \to 
\mathrm{coker}(g) \to 0 
\end{gather*}
and $\mathrm{ker} (f)$, $\mathrm{coker} (f)$, $\mathrm{ker} (g)$, $\mathrm{coker} (g)$ are sheaves with support 
in $Z$, i.e. in the image of $j_{\ast}$. Thus they will be isomorphic to an object in the subcategory generated
 by $j_{\ast} Z_1, \ldots , j_{\ast} Z_n$. In conclusion we see that if one coherent extension $\overline{
\mathcal{G}}_1$ of $\mathcal{G}$ is isomorphic to an object in the subcategory generated by $j_{\ast} Z_1,
 \ldots , j_{\ast} Z_n,\, \mathcal{O}_Y$, the same will be true for any other coherent extension $
\overline{\mathcal{G}}_2$.\\
The rest of the proof is now clear: We split ($\ast$) into short exact sequences and write down extensions of 
these to $Y$ by item (2) above. Since the $\mathcal{L}_{i}$ are finite direct sums of $\mathcal{O}_U$ one 
deduces from the preceding observation that $\mathcal{F}$ is indeed isomorphic to an object in the triangulated 
subcategory generated by $j_{\ast} Z_1, \ldots , j_{\ast} Z_n,\, \mathcal{O}_Y$.
\end{proof}

\begin{remark}
On $\mathbb{P}^n$ it is possible to prove Beilinson's theorem with the help of proposition 2.3.1. Indeed the 
structure sheaves of a flag of linear subspaces $\{ \mathcal{O}_{\mathbb{P}^n}, \: \mathcal{O}_{\mathbb{P}^{n-1}}
, \ldots , \mathcal{O}_{\mathbb{P}^1}, \: \mathcal{O}_{\mathbb{P}^0}\}$ admit the Koszul resolutions
\begin{gather*}
0\to \mathcal{O}(-1)\to \mathcal{O}\to \mathcal{O}_{\mathbb{P}^{n-1}}\to 0 \\
0\to \mathcal{O}(-2) \to \mathcal{O}(-1)^{\oplus 2}\to \mathcal{O}\to \mathcal{O}_{\mathbb{P}^{n-2}}\to 0 \\
\vdots\\
0\to \mathcal{O}(-n) \to \mathcal{O}(-(n-1))^{\oplus n}\to  \ldots \to \mathcal{O}(-1)^{\oplus n}\to 
\mathcal{O}\to \mathcal{O}_{\mathbb{P}^{0}}\to 0 
\end{gather*}
from which one concludes inductively that $(\mathcal{O}(-n),\ldots , \mathcal{O}(-1),\: \mathcal{O})$ generates 
$D^b(Coh\, \mathbb{P}^n)$.
\end{remark}

Next we want to explain a point of view on exceptional sequences that in particular makes obvious the fact that 
the number of terms in any complete exceptional sequence on $X=G/P$ equals the number $| W^P |$ of 
Schubert varieties in $X$. 

\begin{definition}
Let $\mathcal{T}$ be a triangulated category. The \emph{Grothendieck group} $K_{\circ} (\mathcal{T})$ 
\emph{of} $\mathcal{T}$ is the quotient of the free abelian group on the isomorphism classes $[A]$ of 
objects of $\mathcal{T}$ by the subgroup generated by expressions 
\[
[A] - [B] + [C]
\] 
for every distinguished triangle $A\longrightarrow B \longrightarrow C \longrightarrow A[1]$ in $\mathcal{T}$. 
\end{definition}

If $\mathcal{T}=D^b(\mathcal{A})$, $\mathcal{A}$ an Abelian category, then we also have $K_{\circ}(\mathcal{A})$ 
the Grothendieck group of $\mathcal{A}$, i.e. the free abelian group on the isomorphism classes of objects of 
$\mathcal{A}$ modulo relations $[D']-[D]+[D'']$ for every short exact sequence $0\to D'\to D \to D''\to 0$ in 
$\mathcal{A}$, and it is clear that in this case $K_{\circ}(D^b(\mathcal{A}))\simeq K_{\circ }(\mathcal{A})$ (
to a complex $A\in\mathrm{obj}\, D^b(\mathcal{A})$ one associates $\sum_{i\in\mathbb{Z}} (-1)^i [H^i(A)]\in 
K_{\circ}(\mathcal{A})$ which is a map that is additive on distinguished triangles by the long exact cohomology 
sequence and hence descends to a map $K_{\circ}(D^b(\mathcal{A}))\to  K_{\circ }(\mathcal{A})$; the inverse map 
is induced by the embedding $\mathcal{A}\hookrightarrow D^b( \mathcal{A})$).\\
Let now $Y$ be some smooth projective variety. Then to $Z_1 ,\: Z_2\in \mathrm{obj}\, D^b(Coh\, Y)$ one can 
assign the integer $\sum_{i\in\mathbb{Z}} (-1)^i \dim_{\mathbb{C}} \mathrm{Ext}^i (Z_1 , Z_2)$, a map which is 
biadditive on distinguished triangles. Set $K_{\circ}(Y):=K_{\circ}(Coh\, Y)$.
\begin{definition}
The (in general nonsymmetric) bilinear pairing
\begin{gather*}
\chi : K_{\circ}(Y)\times K_{\circ} (Y) \to \mathbb{Z} \\
\hspace{5cm}([Z_1], [Z_2]) \mapsto \sum_{i\in\mathbb{Z}} (-1)^i \dim_{\mathbb{C}} \mathrm{Ext}^i (Z_1 , Z_2)
\end{gather*}
is called the \emph{Euler bilinear form} (cf.\cite{Gor}). 
\end{definition}

\begin{proposition}
Suppose that the derived category $D^b(Coh\, Y)$ of a smooth projective variety $Y$ is generated by an exceptional 
sequence $(E_1, \ldots , E_n)$. Then $K_{\circ }(Y)\simeq
\mathbb{Z}^n$ is a free $\mathbb{Z}$-module of rank $n$ 
with basis given by $([E_1], \ldots , [E_n])$.\\
The Euler bilinear form $\chi$ is unimodular with Gram matrix with respect to the basis $([E_1], \ldots , [E_n])$:
\[
\left( \begin{array}{ccccc}
1 &   &   &  &   \\
0 & 1 &   & \ast  &    \\
0 & 0 & 1 &  &    \\
\vdots & \vdots & \vdots & \ddots &  \\
0 & 0 & 0 & \cdots & 1 
\end{array}
\right)\, ;
\]
in other words, $([E_1], \ldots , [E_n])$ is a semi-orthonormal basis
w.r.t. $\chi$.\\
Moreover, $n=\mathrm{rk}\, K_{\circ}(Y)=\bigoplus_{r=0}^{\dim
  Y}\mathrm{rk}\, A^r(Y)$, where $A^r(Y)$ is the group of codimension
$r$ algebraic cycles on $Y$ modulo rational equivalence (so that
$A(Y)=\bigoplus_r A^r(Y)$ is the Chow ring of $Y$). 
\end{proposition}

\begin{proof}
Since the $E_i$, $i=1,\ldots , n$, generate $D^b(Coh\, Y)$ in the sense of definition 2.1.2 it is clear that the 
$[E_i]$ generate $K_{\circ}(Y)$ (note that for $X$, $X'$, $X''\in\mathrm{obj}\, D^b(Coh\, Y)$ we have 
$[X[n]]=(-1)^n [X]$, $n\in\mathbb{Z}$, $[X'\oplus X'']= [X']+[X'']$ and for every distinguished triangle 
$X'\to X \to X''\to X'[1]$ one has $[X'']= [X]-[X']$).\\
$[E_1]\neq 0$ because $\chi ([E_1], [E_1])=1$ since $E_1$ is exceptional. Assume inductively that $[E_1], \ldots ,
[E_i]$ are linearly independent in $K_{\circ}(Y)\otimes \mathbb{Q}$. We claim $[E_{i+1}]\notin \langle [E_1],
 \ldots
 , [E_i]\rangle_{\mathbb{Q}}$. Indeed otherwise $[E_{i+1}]=\sum_{j=1}^i \lambda_j [E_j]$; since $[E_{i+1}]\neq 0$ 
there is $l:= \mathrm{min}\{ j\, \mid\, \lambda_j\neq 0 \}$. Then \[\chi ([E_{i+1}] , [E_{l}])=\chi (\sum_
{j=l}^i \lambda_j [E_j] , [E_l] )=\lambda_l\neq 0\] (using $\mathrm{Ext}^k(E_j , E_i)=0$ $\forall k\in \mathbb{Z
}\; \forall i< j$) contradicting the fact that\\ $\chi ([E_{i+1}], [E_l])=0$ since $l< i+1$. Thus the $([E_1], 
\ldots , [E_n])$ form a free $\mathbb{Z}$-basis of $K_{\circ}(Y)$. The remaining assertions concerning $\chi$ 
are obvious from the above arguments.\\
The last equality follows from the fact that the Grothendieck Chern
character $\mathrm{ch}$ gives an isomorphism
\[
\mathrm{ch} : K_{\circ}(Y)\otimes \mathbb{Q} \to A(Y)\otimes \mathbb{Q}
\]
(cf. \cite{Ful}, 15.2.16 (b)).
\end{proof}

\begin{corollary}
If $(E_1,\ldots , E_n)$ is an exceptional sequence that generates $D^b(Coh\, X)$, $X$ a rational homogeneous 
variety, then $n=| W^P|$, the number of Schubert varieties $X_w$ in $X$. 
\end{corollary}

\begin{proof}
It suffices to show that the $[\mathcal{O}_{X_w}]$'s likewise form a
free $\mathbb{Z}$-basis of $K_{\circ }(X)$.
 One way to see this is 
as follows: By proposition 2.3.1 it is clear that the $[\mathcal{O}_{X_w}]$ generate $K_{\circ}(X)$. $K_{\circ}
(X)$ is a ring for the product $[\mathcal{F}]\cdot [\mathcal{G}]:=\sum_{i\in\mathbb{Z}} (-1)^i [\mathcal{T}or_i^X
(\mathcal{F}, \mathcal{G}) ]$ and 
\begin{gather*}
\beta : K_{\circ}(X)\times K_{\circ} (X) \to \mathbb{Z} \\
\hspace{5cm}([\mathcal{F}], [\mathcal{G}]) \mapsto \sum_{i\in\mathbb{Z}} (-1)^i h^i (X , [\mathcal{F}]\cdot [
\mathcal{G}])
\end{gather*}
is a symmetric bilinear form. One can compute that $\beta ([\mathcal{O}_{X_x}], [\mathcal{O}_{X^y}(-\partial X^y)])
=\delta_{x}^y$ (Kronecker delta) for $x,\: y\in W^P$, cf. \cite{BL}, proof of lemma 6, for details.  
\end{proof}

It should be noted at this point that the constructions in subsection 2.1 relating to semi-orthogonal 
decompositions, mutations etc. all have their counterparts on the K-theory level and in fact appear more 
natural in that context (cf. \cite{Gor}, \S 1). 

\begin{remark}
Suppose that on $X=G/P$ we have a strong complete exceptional sequence $(E_1,\ldots , E_n)$. Then the Gram 
matrix $G$ of $\chi$ w.r.t. the basis $([E_1], \ldots , [E_n])$ on $K_{\circ}(X)\simeq \mathbb{Z}^n$ is upper 
triangular with ones on the diagonal and $(i,j)$-entry equal to $\dim_{\mathbb{C}} \mathrm{Hom}(E_i,E_j)$. Thus 
with regard to conjecture 2.2.1 it would be interesting to know the Gram matrix $G'$  of $\chi$ in the basis given by 
the $[\mathcal{O}_{X_w}]$'s, $w\in W^P$, since $G$ and $G'$ will be
conjugate.\\
 The following computation was suggested to me by M. Brion. Without
 loss of generality one may reduce to the case $X=G/B$ using the
 fibration $\pi : G/B \to G/P$: Indeed, the pull-back under $\pi$ of
 the Schubert variety $X_{wP}$, $w\in W^P$, is the Schubert variety
 $X_{w\, w_{0, P}}$ in $G/B$ where $w_{0, P}$ is the element of
 maximal length of $W_P$, and $\pi^{\ast}
 \mathcal{O}_{X_{wP}}=\mathcal{O}_{X_{w\, w_{0, P}}}$. Moreover, by
 the projection formula and because
 $R\pi_{\ast}\mathcal{O}_{G/B}=\mathcal{O}_{G/P}$, we have
 $R\pi_{\ast}\circ\pi^{\ast}\simeq \mathrm{id}_{D^b(Coh\, G/P)}$ and 
\begin{gather*}
\chi (\pi^{\ast}\mathcal{E} , \pi^{\ast}\mathcal{F})=\chi (\mathcal{E} ,\mathcal{F})
\end{gather*}
for any $\mathcal{E} ,\:\mathcal{F}\in \mathrm{obj}\, D^b(Coh\, G/P)$.\\
Therefore, let $X=G/B$ and let $x,\: y\in W$.
The first observation is that $X_y=w_0X^{w_0y}$ and $\chi (\mathcal{O}_{X_x}, \mathcal{O}_{X_y}) = \chi (\mathcal{O}
_{X_x}, \mathcal{O}_{X^{w_0y}})$. This follows from the facts that
there is a connected chain of rational curves in $G$ joining $g$ to
$\mathrm{id}_G$ (since $G$ is generated by images of homomorphisms
$\mathbb{C}\to G$ and $\mathbb{C}^{\ast}\to G$) and that flat families
of sheaves indexed by open subsets of $\mathbb{A}^1$ yield the same
class in $K_{\circ}(X)$, thus $[\mathcal{O}_{X^{w_0 y}}]=[\mathcal{O}_{w_0X^{w_0 y}}]$.   
 We have 
\begin{gather*}
R\mathrm{Hom}^{\bullet}(\mathcal{O}_{X_x} , \mathcal{O}_{X^{w_0 y}}) \simeq R\Gamma (X , R\mathcal{H}om^{\bullet} (
\mathcal{O}_{X_x} , \mathcal{O}_{X^{w_0 y}})) \\
\simeq  R\Gamma (X , R\mathcal{H}om^{\bullet} (\mathcal{O}_{X_x} , \mathcal{O}_{X})\otimes^L \mathcal{O}_{X^{w_0y}})
\end{gather*}
(cf. \cite{Ha1}, prop. 5.3/5.14). Now Schubert varieties are Cohen-Macaulay, in fact they have 
rational singularities (cf. \cite{Ra1}), whence 
\begin{gather*}
R\mathcal{H}om^{\bullet} (\mathcal{O}_{X_x},\mathcal{O}_X) \simeq \mathcal{E}xt^{\mathrm{codim} (X_x)} (
\mathcal{O}_{X_x},\mathcal{O}_X) [-\mathrm{codim} (X_x)]\\
\simeq \omega_{X_x}\otimes \omega^{-1}_X [-\mathrm{codim} (X_x)]\, .
\end{gather*}
But $\omega_{X_x}\otimes \omega^{-1}_X\simeq \mathcal{L}(\varrho)|_{X_x}(-\partial X_x)$ ($\mathcal{L}(\varrho )$
 is the line bundle associated to the character $\varrho$), cf. \cite{Ra1}, prop. 2 and thm. 4.
 Now $X_x$ and $X^{w_0y}$ are 
Cohen-Macaulay and their scheme theoretic intersection is proper in $X$ and reduced (\cite{Ra1}, thm. 3) whence
 $\mathcal{T}or_i^X(\mathcal{O}_{X_x}, \mathcal{O}_{X^{w_0y}}) 
=0$ for all $i\ge 1$ (cf. \cite{Bri}, lemma 1). Therefore 
\[
R\mathrm{Hom}^{\bullet}(\mathcal{O}_{X_x} , \mathcal{O}_{X^{w_0y}}) \simeq R\Gamma (X ,\mathcal{L}(\varrho)|_{X_x}
(-\partial X_x)[-\mathrm{codim} (X_x)]\otimes \mathcal{O}_{X^{w_0y}} )
\]
so that setting $X_x^{w_0y}:= X_x\cap X^{w_0y}$ and $(\partial X_x)^{w_0y}:= \partial X_x\cap X^{w_0 y}$
\[
\chi (\mathcal{O}_{X_x}, \mathcal{O}_{X_y})= (-1)^{\mathrm{codim}(X_x)} \chi \left(\mathcal{L}(\varrho)|_{X_x^{w_0y
}} (-(\partial X_x)^{w_0y})\right)\, .  
\] 
This is $0$ unless $w_0y \le x$ (because $X^{w_0y}_x$ is non-empty iff $w_0y\le x$, see \cite{BL}, lemma 1); 
moreover if $w_0 y \le x$ there are no higher $h^i$ in the latter Euler characteristic by \cite{BL}, prop. 2. In 
conclusion
\[
\chi (\mathcal{O}_{X_x}, \mathcal{O}_{X_y}) = \left\{  \begin{array}{l} (-1)^{\mathrm{codim}(X_x)} h^0
 \left(\mathcal{L}(\varrho)|_{X_x^{w_0y}} (-(\partial X_x)^{w_0y}) \right)\quad \mathrm{if}\quad w_0y \le x \\
0 \quad \mathrm{otherwise}
\end{array}                 \right.
\]
though the impact of this on conjecture 2.2.1 ((A) or (B)) is not clear to me.\\
Cf. also \cite{Bri2} for this circle of ideas.
\end{remark}

\section{Fibrational techniques}

The main idea pervading this section is that the theorem of Beilinson on
the structure of the derived category of coherent sheaves on projective
space (\cite{Bei}) and the related results of Kapranov (\cite{Ka1},
\cite{Ka2}, \cite{Ka3}) for Grassmannians, flag varieties and quadrics,
generalize without substantial difficulty from the absolute to the
relative setting, i.e. to projective bundles etc. For projective bundles,
Grassmann and flag bundles this has been done in \cite{Or}. We review
these results in subsection 3.1; the case of quadric bundles is dealt
with in subsection 3.2. Aside from being technically a little more
involved, the result follows rather mechanically combining the techniques
from \cite{Ka3} and \cite{Or}. Thus armed, we deduce information on the
derived category of coherent sheaves on isotropic Grassmannians and flag
varieties in the symplectic and orthogonal cases; we follow an idea first
exploited in \cite{Sa} using successions of projective and quadric
bundles.

\subsection{The theorem of Orlov on projective bundles}

Let $X$ be a smooth projective variety, $\mathcal{E}$ a vector bundle of
rank $r+1$ on $X$. Denote by $\mathbb{P}(\mathcal{E})$ the associated
projective bundle \footnote[2]{Here and in the following
$\mathbb{P}(\mathcal{E})$ denotes
$\mathrm{Proj}(\mathrm{Sym}^{\bullet}(\mathcal{E}^{\vee}))$, i.e. the
bundle of 1-dimensional subspaces in the fibres of $\mathcal{E}$, and
contrary to Grothendieck's notation
\textbf{not} the bundle
$\mathrm{Proj}(\mathrm{Sym}^{\bullet}\,\mathcal{E})$ of hyperplanes in
the fibres of $\mathcal{E}$ which might be less intuitive in the sequel.}
and
$\pi :
\mathbb{P}(\mathcal{E})
\to X$ the projection. Set $D^b(\mathcal{E}):=D^b(Coh
(\mathbb{P}(\mathcal{E})))$, $D^b(X):=D^b(Coh (X))$. There are the functors
$R\pi _{\ast} : D^b(\mathcal{E})\to D^b(X)$ (note that $R\pi_{\ast } :
D^{+}(Coh (\mathbb{P}(\mathcal{E}))) \to D^{+}(Coh (X))$ maps
$D^b(\mathcal{E})$ to
$D^b(X)$ using $R^i\pi_{\ast} (\mathcal{F})=0 \, \forall i> \dim
\mathbb{P}(\mathcal{E})\: \forall \mathcal{F} \in \mathrm{Ob}\, Coh
(\mathbb{P}(\mathcal{E}))$ and the spectral sequence in hypercohomology) and
$\pi^{\ast} : D^b(X)\to D^b(\mathcal{E})$ ($\pi$ is flat, hence $\pi^{\ast}$
is exact and passes to the derived category without taking the left
derived functor).\\
We identify $D^b(X)$ with a full subcategory in
$D^b(\mathcal{E})$ via $\pi^{\ast}$ (cf. \cite{Or}, lemma 2.1). More generally
we denote by $D^b(X)\otimes \mathcal{O}_{\mathcal{E}}(m)$ for
$m\in\mathbb{Z}$ the subcategory of $D^b(\mathcal{E})$ which is the image of
$D^b(X)$ in $D^b(\mathcal{E})$ under the functor $\pi^{\ast}(-)\otimes
\mathcal{O}_{\mathcal{E}}(m)$, where $\mathcal{O}_{\mathcal{E}}(1)$ is the
relative hyperplane bundle on $\mathbb{P}(\mathcal{E})$. Then one has the
following result (cf. \cite{Or}, thm. 2.6):
\begin{theorem}
The categories $D^b(X)\otimes\mathcal{O}_{\mathcal{E}}(m)$ are all admissible
subcategories of $D^b(\mathcal{E})$ and we have a semiorthogonal decomposition
\[
D^b(\mathcal{E})=\left\langle D^b(X)\otimes
\mathcal{O}_{\mathcal{E}}(-r),\ldots ,
D^b(X)\otimes\mathcal{O}_{\mathcal{E}}(-1), D^b(X)
\right\rangle \, .
\]
\end{theorem}
We record the useful
\begin{corollary}
If $D^b(X)$ is generated by a complete exceptional sequence
\[
(E_1,\ldots ,E_n)\, ,
\]
then $D^b(\mathcal{E})$ is generated by the complete exceptional sequence
\[
(\pi^{\ast}E_1\otimes
\mathcal{O}_{\mathcal{E}}(-r),
\ldots ,
\pi^{\ast}E_n\otimes \mathcal{O}_{\mathcal{E}}(-r),\pi^{\ast}E_1\otimes
\mathcal{O}_{\mathcal{E}}(-r+1), \ldots ,\pi^{\ast}E_1,
\ldots ,
\pi^{\ast}E_n) \,.
\]
\end{corollary}
\begin{proof}
This is stated in \cite{Or}, cor. 2.7; for the sake of completeness
and because the method will be used repeatedly in the sequel, we give a
proof. One just checks that
\[
\mathrm{Ext}^k(\pi^{\ast}E_i\otimes \mathcal{O}_{\mathcal{E}}(-r_1),
\pi^{\ast}E_j\otimes \mathcal{O}_{\mathcal{E}}(-r_2))=0 \]
$\forall k,
\forall 1\le i,\, j\le n\, \forall 0\le r_1< r_2\le r
$ and $\forall k,\, \forall 1\le j< i\le n, \, r_1=r_2$. Indeed,
\begin{eqnarray*}
\mathrm{Ext}^k(\pi^{\ast}E_i\otimes \mathcal{O}_{\mathcal{E}}(-r_1),
\pi^{\ast}E_j\otimes \mathcal{O}_{\mathcal{E}}(-r_2))
\simeq \mathrm{Ext}^k(\pi^{\ast}E_i,
\pi^{\ast}E_j\otimes \mathcal{O}_{\mathcal{E}}(r_1-r_2))\\
\simeq \mathrm{Ext}^k(E_i,
E_j\otimes R\pi_{\ast}(\mathcal{O}_{\mathcal{E}}(r_1-r_2)))
\end{eqnarray*}
where for the second isomorphism we use that $R\pi_{\ast}$ is right adjoint to
$\pi^{\ast}$, and the projection formula (cf. \cite{Ha2}, II, prop. 5.6). When
$r_1=r_2$ and $i>j$ then $R\pi_{\ast} \mathcal{O}_{\mathcal{E}} \simeq
\mathcal{O}_X$ and
$\mathrm{Ext}^k(E_i,E_j)=0$ for all $k$ because $(E_1,\ldots ,E_n)$ is
 exceptional. If on the other hand $0\le r_1<r_2\le r$ then $-r\le r_1-r_2<0$
and $R\pi_{\ast} (\mathcal{O}_{\mathcal{E}}(r_1-r_2))=0$.\\
It remains to see that each $\pi^{\ast}E_i\otimes
\mathcal{O}_{\mathcal{E}}(-r_1)$ is exceptional. From the above calculation it
is clear that this follows exactly from the exceptionality of $E_i$.
\end{proof}
\begin{remark}
From the above proof it is clear that even if we start in corollary 3.1.2 with
a strong complete exceptional sequence $(E_1,\ldots ,E_n)$ (i.e.
$\mathrm{Ext}^k(E_i,E_j)=0\: \forall i,\, j\: \forall k\neq 0$), the resulting
exceptional sequence on
$\mathbb{P}(\mathcal{E})$ need not again be strong: For example take
$X=\mathbb{P}^1$ with strong complete exceptional sequence
$(\mathcal{O}(-1), \mathcal{O})$ and $\mathcal{E}=\mathcal{O}\oplus
\mathcal{O}(h)$, $h\ge 2$, so that
$\mathbb{P}(\mathcal{E})=\mathbb{F}_h\stackrel{\pi}{\longrightarrow}
\mathbb{P}^1$ is a Hirzebruch surface. Then
$(\pi^{\ast}\mathcal{O}(-1)\otimes\mathcal{O}_{\mathcal{E}}(-1),
\mathcal{O}_{\mathcal{E}}(-1) ,
\pi^{\ast}\mathcal{O}(-1)\otimes\mathcal{O}_{\mathcal{E}},
\mathcal{O}_{\mathcal{E}})$ is an exceptional sequence on $\mathbb{F}_h$ that
generates $D^b(Coh (\mathbb{F}_h))$, but it is not a strong one since
$\mathrm{Ext}^1(\mathcal{O}_{\mathcal{E}}(-1),
\mathcal{O}_{\mathcal{E}})\simeq$ $H^1(\mathbb{P}^1, \pi_{\ast}\,
\mathcal{O}_{\mathcal{E}}(1))\simeq H^1(\mathbb{P}^1,
\mathcal{O}\oplus\mathcal{O}(-h))\simeq \mathrm{Sym}^{h-2}\,\mathbb{C}^2\neq 0$.
\end{remark}
Analogous results hold for relative Grassmannians and flag varieties.
Specifically, if $\mathcal{E}$ is again a rank $r+1$ vector bundle on a smooth
projective variety $X$, denote by $\mathrm{Grass}_X(k , \mathcal{E})$ the
relative Grassmannian of $k$-planes in the fibres of $\mathcal{E}$ with
projection $\pi : \mathrm{Grass}_X(k , \mathcal{E}) \to X$ and tautological
subbundle $\mathcal{R}$ of rank $k$ in $\pi^{\ast} \mathcal{E}$. Denote by
$Y(k, r+1-k)$ the set of partitions $\lambda =(\lambda_1, \ldots ,\lambda_k)$
with $0\le \lambda_k\le \lambda_{k-1}\le \ldots \le \lambda_1\le r+1-k$ or
equivalently the set of Young diagrams with at most $k$ rows and no more than
$r+1-k$ columns. For $\lambda\in V(k, r+1-k)$ we have the Schur functor
$\Sigma^{\lambda}$ and bundles $\Sigma^{\lambda }\mathcal{R}$ on
$\mathrm{Grass}_X(k , \mathcal{E})$. Moreover, as before we can talk about full
subcategories $D^b(X)\otimes \Sigma^{\lambda }\mathcal{R}$ of $D^b(Coh (\mathrm{Grass}_X(k
,\mathcal{E})))$. Choose a total
order
$\prec$ on
$Y(k ,r+1-k)$ such that if $\lambda \prec \mu$ then the Young diagram of
$\lambda$ is not contained in the Young diagram of $\mu$, i.e. $\exists i :
\mu_i <\lambda_i $. Then one has (cf. \cite{Or}, p. 137):
\begin{theorem}
There is a semiorthogonal decomposition
\[
D^b(Coh (\mathrm{Grass}_X(k
,\mathcal{E})))=\left\langle \ldots , D^b(X)\otimes
\Sigma^{\lambda }\mathcal{R},\ldots ,
D^b(X)\otimes\Sigma^{\mu }\mathcal{R}, \ldots
\right\rangle \,
\]
($\lambda\prec\mu$).\\
If $(E_1,\ldots ,E_n)$ is a complete exceptional sequence in $D^b(X)$, then
\[
\left( \ldots , \pi^{\ast}E_1\otimes \Sigma^{\lambda}\mathcal{R},\ldots ,
\pi^{\ast}E_n\otimes \Sigma^{\lambda}\mathcal{R}, \ldots ,
 \pi^{\ast}E_1\otimes \Sigma^{\mu}\mathcal{R},\ldots ,
\pi^{\ast}E_n\otimes \Sigma^{\mu}\mathcal{R}, \ldots \right)
\]
is a complete exceptional sequence in $D^b(Coh (\mathrm{Grass}_X(k
,\mathcal{E})))$. Here all $\pi^{\ast}E_i\otimes \Sigma^{\lambda}\mathcal{R}$,
$i\in\{1 ,\ldots ,n\}$, $\lambda\in Y(k, r+1-k)$ occur in the list, and
$\pi^{\ast}E_i\otimes \Sigma^{\lambda}\mathcal{R}$ precedes $\pi^{\ast}E_j\otimes
\Sigma^{\mu}\mathcal{R}$ iff $\lambda\prec\mu$ or $\lambda =\mu$ and $i<j$.
\end{theorem}
More generally, we can consider for $1\le k_1<\ldots < k_t\le r+1$ the 
variety $\mathrm{Flag}_X(k_1,\ldots ,k_t; \mathcal{E})$ of relative flags of type $(k_1,
\ldots , k_t)$ in the fibres of $\mathcal{E}$, with projection $\pi$ and
tautological subbundles $\mathcal{R}_{k_1}\subset \ldots \subset
\mathcal{R}_{k_t}\subset \pi^{\ast}\mathcal{E}$. If we denote again by $Y(a,b)$
the set of Young diagrams with at most $a$ rows and $b$ columns, we consider the
sheaves $\Sigma^{\lambda_1} \mathcal{R}_{k_1}\otimes \ldots \otimes
\Sigma^{\lambda_t} \mathcal{R}_{k_t}$ on $\mathrm{Flag}_X(k_1,\ldots ,k_t;
\mathcal{E})$ with $\lambda_k\in Y(k_t, r+1-k_t)$ and $\lambda_j\in Y(k_j ,
k_{j+1}-k_j)$ for $j=1,\ldots ,t-1$ and subcategories
$D^b(X)\otimes\Sigma^{\lambda_1} 
\mathcal{R}_{k_1}\otimes \ldots \otimes
\Sigma^{\lambda_t} \mathcal{R}_{k_t}$ of $D^b(Coh (\mathrm{Flag}_X(k_1,\ldots ,k_t
;\mathcal{E})))$. Choose a
total order
$\prec_j$ on each of the sets $Y(k_j, k_{j+1}-k_j)$ and $ \prec_t$ on $Y(k_t,
r+1-k_t)$ with the same property as above for the relative Grassmannian, and
endow the set
$Y=Y(k_t ,r+1-k_t)\times \ldots \times Y(k_1,k_2-k_1)$ with the resulting
lexicographic order $\prec$. 
\begin{theorem}
There is a semiorthogonal decomposition
\begin{gather*}
D^b(Coh (\mathrm{Flag}_X(k_1,\ldots ,k_t
;\mathcal{E})))=\left\langle \ldots ,D^b(X)\otimes\Sigma^{\lambda_1}
\mathcal{R}_{k_1}\otimes \ldots \otimes
\Sigma^{\lambda_t} \mathcal{R}_{k_t} ,\right. \\
\left. \ldots ,
 D^b(X)\otimes\Sigma^{\mu_1} \mathcal{R}_{k_1}\otimes \ldots \otimes
\Sigma^{\mu_t} \mathcal{R}_{k_t}, \ldots
\right\rangle \, 
\end{gather*}
($(\lambda_t,\ldots ,\lambda_1)\prec (\mu_t,\ldots ,\mu_1)$).\\
If $(E_1,\ldots ,E_n)$ is a complete exceptional sequence in $D^b(X)$, then
\begin{eqnarray*}
\left( \ldots , \pi^{\ast}E_1\otimes
\Sigma^{\lambda_1}\mathcal{R}_{k_1}\otimes \ldots \otimes
\Sigma^{\lambda_t}\mathcal{R}_{k_t},\ldots ,
\pi^{\ast}E_n \otimes \Sigma^{\lambda_1}\mathcal{R}_{k_1}\otimes \ldots \otimes
\Sigma^{\lambda_t}\mathcal{R}_{k_t},\ldots ,\right. \\
\left. \pi^{\ast}E_1\otimes
\Sigma^{\mu_1}\mathcal{R}_{k_1}\otimes \ldots \otimes
\Sigma^{\mu_t}\mathcal{R}_{k_t},\ldots ,
\pi^{\ast}E_n \otimes \Sigma^{\mu_1}\mathcal{R}_{k_1}\otimes \ldots \otimes
\Sigma^{\mu_t}\mathcal{R}_{k_t},\ldots
 \right)
\end{eqnarray*}
is a complete exceptional sequence in $D^b(Coh (\mathrm{Flag}_X(k_1,\ldots ,k_t
;\mathcal{E})))$. Here all $\pi^{\ast}E_i\otimes
\Sigma^{\lambda_1}\mathcal{R}_{k_1}\otimes \ldots \otimes
\Sigma^{\lambda_t}\mathcal{R}_{k_t}$,
$i\in\{1 ,\ldots ,n\}$, $(\lambda_t, \ldots ,\lambda_1) \in Y$ occur in the list,
and
$\pi^{\ast}E_i\otimes
\Sigma^{\lambda_1}\mathcal{R}_{k_1}\otimes \ldots \otimes
\Sigma^{\lambda_t}\mathcal{R}_{k_t}$ precedes $\pi^{\ast}E_j\otimes
\Sigma^{\mu_1}\mathcal{R}_{k_1}\otimes \ldots \otimes
\Sigma^{\mu_t}\mathcal{R}_{k_t}$ iff
$(\lambda_t,\ldots ,\lambda_1)\prec (\mu_t,\ldots ,\mu_1)$ or $(\lambda_t,\ldots
,\lambda_1)=(\mu_t,\ldots ,\mu_1)$ and
$i<j$.
\end{theorem}
\begin{proof}
Apply theorem 3.1.4 iteratively to the succession of Grassmann bundles
\begin{eqnarray*}
\mathrm{Flag}_X(k_1,\ldots ,k_t;
\mathcal{E})=\mathrm{Grass}_{\mathrm{Flag}_X(k_2,\ldots ,k_t; \mathcal{E})}(k_1
, \mathcal{R}_{k_2})\\
\rightarrow \mathrm{Flag}_X(k_2,\ldots ,k_t;
\mathcal{E})=\mathrm{Grass}_{\mathrm{Flag}_X(k_3,\ldots ,k_t; \mathcal{E})}(k_2
, \mathcal{R}_{k_3}) \rightarrow \ldots \rightarrow X
\end{eqnarray*}
\end{proof}  

\subsection{The theorem on quadric bundles}
Let us now work out in detail how the methods of Orlov (\cite{Or}) and Kapranov
(\cite{Ka2}, \cite{Ka3}) yield a result for quadric bundles that is analogous to
theorems 3.1.1, 3.1.4, 3.1.5.\\
As in subsection 3.1, $X$ is a smooth projective variety with a vector bundle
$\mathcal{E}$ of rank $r+1$ endowed with a symmetric quadratic form $q\in\Gamma
(X, \mathrm{Sym}^2\, \mathcal{E}^{\vee})$ which is nondegenerate on each fibre;
$\mathcal{Q}:=\{ q=0\} \subset \mathbb{P}(\mathcal{E})$ is the associated quadric
bundle:
\setlength{\unitlength}{1cm}
\begin{center}
\begin{picture}(4,3)
\put(0.7,2.5){$\mathcal{Q}$}
\put(1.75,2.5){$\mathbb{P}(\mathcal{E})$}
\put(1.75,0.7){$X$}
\put(1.2,2.5){$\hookrightarrow$}
\put(0.9,2.3){\vector(2,-3){0.8}}
\put(1.85,2.3){\vector(0,-1){1.2}}
\put(0.0,1.5){$\pi=\Pi |_{\mathcal{Q}}$}
\put(2.0,1.6){$\Pi$} 
\end{picture}\end{center}
Write $D^b(X):=D^b(Coh\, X)$, $D^b(\mathcal{Q}):=D^b(Coh\, \mathcal{Q})$,
$D^b(\mathcal{E}):=D^b(Coh\, \mathbb{P}(\mathcal{E}))$.
\begin{lemma}
The functor
\[
\pi^{\ast}=L\pi^{\ast} : D^b(X) \to D^b(\mathcal{Q})
\]
is fully faithful.
\end{lemma}
\begin{proof}
Since $\mathcal{Q}$ is a locally trivial fibre bundle over $X$ with
rational homogeneous fibre, we have
$\pi_{\ast}\mathcal{O}_{\mathcal{Q}}=\mathcal{O}_X$ and
$R^i\pi_{\ast}\mathcal{O}_{\mathcal{Q}}=0$ for $i>0$. The right adjoint to
$L\pi^{\ast}$ is $R\pi_{\ast}$, and $R\pi_{\ast}\circ L\pi^{\ast}$ is isomorphic
to the identity on
$D^b(X)$ because of the projection formula and
$R\pi_{\ast}\mathcal{O}_{\mathcal{Q}}=\mathcal{O}_X$. Hence $L\pi^{\ast}$ is
fully faithful (and equal to $\pi^{\ast}$ since $\pi$ is flat). 
\end{proof}
Henceforth $D^b(X)$ is identified with a full subcategory of $D^b(\mathcal{Q})$.\\
We will now define two bundles of graded algebras, $\mathcal{A}=\bigoplus\limits_{n\ge 0}
\mathcal{A}_n$ and $\mathcal{B}=\bigoplus\limits_{n\ge
  0}\mathcal{B}_n$, on $X$. 
Form the tensor algebra $\mathrm{T}^{\bullet}(\mathcal{E}[h])$ where $h$ is an
indeterminate with $\deg h=2$ and germs of sections in $\mathcal{E}$ have degree
$1$ and take the quotient modulo the two-sided ideal $\mathcal{I}$ of relations with
$\mathcal{I}(x):=\langle e\otimes e-q(e)h ,\, e\otimes h-h\otimes e
\rangle_{e\in \mathcal{E}(x)}$, ($x\in X$). This quotient is $\mathcal{A}$, the bundle of
graded Clifford algebras of the orthogonal vector bundle $\mathcal{E}$. On the
other hand, $\mathcal{B}$ is simply defined as $\bigoplus\limits_{n\ge
0}\pi_{\ast}\,\mathcal{O}_{\mathcal{Q}} (n)$, the relative coordinate algebra of
the quadric bundle $\mathcal{Q}$.\\
For each graded left $\mathcal{A}$-module
$\mathcal{M}=\bigoplus_{i\in\mathbb{Z}} \mathcal{M}_i$ with $\mathcal{M}_i$
vector bundles on
$X$ we get a complex
$L^{\bullet}(\mathcal{M})$ of bundles on $\mathcal{Q}$
\[
\begin{CD}
L^{\bullet }(\mathcal{M}):\, \ldots @>>> \pi^{\ast }\mathcal{M}_j
\otimes_{\mathbb{C}}
\mathcal{O}_{\mathcal{Q}}(j) @> d^j >>
\pi^{\ast }\mathcal{M}_{j+1} \otimes_{\mathbb{C}} \mathcal{O}_{\mathcal{Q}}(j+1)
@>>>
\ldots
\end{CD}
\]
with differentials given as follows: For $x\in X$ and $e\in\mathcal{E}(x)$ we get a family of
mappings 
\[
d^j(x,e) : \mathcal{M}_{j}(x) \to \mathcal{M}_{j+1}(x)
\]
given by left multiplication by $e$ on $\mathcal{M}_{j}(x)$ and linear in $e$
which globalize to mappings
$\Pi^{\ast}\mathcal{M}_j\otimes\mathcal{O}_{\mathcal{E}}(j)\to
\Pi^{\ast}\mathcal{M}_{j+1}\otimes\mathcal{O}_{\mathcal{E}}(j+1)$. When
restricted to $\mathcal{Q}$ two successive maps compose to $0$ and we get the
required complex.\\
We recall at this point the relative version of Serre's correspondence (cf. e.g.
\cite{EGA}, II, \S 3):
\begin{theorem}
Let $Mod_{\mathcal{E}}^X$ be the category whose objects are coherent sheaves
over $X$ of graded $\mathrm{Sym}^{\bullet} \mathcal{E}^{\vee}$-modules of finite
type with morphisms
\[
\mathrm{Hom}_{Mod_{\mathcal{E}}^X} (\mathcal{M} ,
\mathcal{N}):=\lim_{\stackrel{\longrightarrow}{n}}
\mathrm{Hom}_{\mathrm{Sym}^{\bullet}\,
\mathcal{E}^{\vee}}(
\bigoplus\limits_{i\ge n} \mathcal{M}_i , \bigoplus\limits_{i\ge n} \mathcal{N}_i
)
\] 
(the direct limit running over the groups of homomorphisms of
sheaves of graded modules over
$\mathrm{Sym}^{\bullet}\,\mathcal{E}^{\vee}$ which are homogeneous of
degree 0). If
$\mathcal{F}\in\mathrm{obj} (Coh (\mathbb{P}(\mathcal{E})))$ set
\[
\alpha (\mathcal{F}) := \bigoplus\limits_{n=0}^{\infty}
\Pi_{\ast}(\mathcal{F}(n))\, .
\]
Then the functor $\alpha : Coh (\mathbb{P}(\mathcal{E}))) \to
Mod_{\mathcal{E}}^X$ is an equivalence of categories with quasi-inverse
$(-)^{\sim}$ which is an additive and exact functor.
\end{theorem}
The key remark is now that $L^{\bullet}(\mathcal{A}^{\vee})$ is exact since it
arises by applying the Serre functor $(-)^{\sim}$ to the complex $P^{\bullet}$
given by
\[
\begin{CD}
\ldots @>d>> \mathcal{A}_2^{\vee}\otimes \mathcal{B}[-2] @>d>> 
\mathcal{A}_1^{\vee}\otimes \mathcal{B}[-1] @>d>> \mathcal{A}_0^{\vee}\otimes
\mathcal{B} @>>> \mathcal{O}_X \rightarrow 0 .
\end{CD}
\]
Here, if $(e_1, \ldots ,e_{r+1})$ is a local frame of $\mathcal{E}=\mathcal{A}_1$
and $(e_1^{\vee}, \ldots , e_{r+1}^{\vee})$ is the corresponding dual frame for
$\mathcal{E}^{\vee}=\mathcal{B}_1$, the differential $d$ is
$\sum_{i=1}^{r+1} l_{e_i}^{\vee}\otimes l_{e_i^{\vee}}$, where $l_{e_i}
:\mathcal{A}[-1]\to\mathcal{A}$ is left multiplication by $e_i$ and analogously
$l_{e_i^{\vee}} : \mathcal{B}[-1] \to \mathcal{B}$. This complex is exact since
it is so fibrewise as a complex of vector bundles; the fibre over a point $x\in
X$ is just Priddy's generalized Koszul complex associated to the dual quadratic
algebras $\mathcal{B}(x)= \oplus_i H^0( \mathcal{Q}(x) ,
\mathcal{O}_{\mathcal{Q}(x)}(i))$ and $\mathcal{A}(x)$, the graded Clifford
algebra of the vector space $\mathcal{E}(x)$. See \cite{Ka3}, 4.1 and \cite{Pri}.
\\
Define bundles $\Psi_i$, $i\ge 0$, on $\mathcal{Q}$ by a twisted
truncation, i.e., by the requirement that
\[
0 \to \Psi_i \to \pi^{\ast}\mathcal{A}_i^{\vee} \to
\pi^{\ast}\mathcal{A}_{i-1}^{\vee} \otimes \mathcal{O}_{\mathcal{Q}}(1) \to
\ldots \to \pi^{\ast}\mathcal{A}_0^{\vee} \otimes \mathcal{O}_{\mathcal{Q}}(i)
\to 0
\]
be exact. Look at the fibre product
\[
\begin{CD}
\Delta\subset \mathcal{Q}\times_X \mathcal{Q} @>p_2>> \mathcal{Q} \\
@Vp_1VV                                  @V\pi VV  \\
\mathcal{Q} @>\pi >>              X
\end{CD}
\]
together with the relative diagonal $\Delta$. The goal is to cook up an
infinite to the left but eventually periodic resolution of the sheaf
$\mathcal{O}_{\Delta}$ on $\mathcal{Q}\times_X\mathcal{Q}$, then truncate it in a
certain degree and identify the remaining kernel explicitly. \\
Write $\Psi_i \boxtimes \mathcal{O}(-i)$ for $p_1^{\ast}\Psi_i \otimes
p_2^{\ast} \mathcal{O}_{\mathcal{Q}}(-i)$ and consider the maps $\Psi_i\boxtimes
\mathcal{O}(-i) \to \Psi_{i-1}\boxtimes \mathcal{O}(-i+1)$ induced by the maps
of complexes
\[
\begin{CD}
(\pi^{\ast}\mathcal{A}_{i}^{\vee}\otimes\mathcal{O}) \boxtimes \mathcal{O}(-i)
@>>>
(\pi^{\ast}\mathcal{A}_{i-1}^{\vee}\otimes\mathcal{O}(1))\boxtimes\mathcal{O}(-i)
@>>> \ldots
\\
 @VVV @VVV 
\\
(\pi^{\ast}\mathcal{A}_{i-1}^{\vee}\otimes\mathcal{O}) \boxtimes \mathcal{O}(-i+1)
@>>>
(\pi^{\ast}\mathcal{A}_{i-2}^{\vee}\otimes\mathcal{O}(1))\boxtimes\mathcal{O}(-i+1)
@>>> \ldots
\end{CD}
\]
where the vertical arrows are given by $\sum_{i=1}^{r+1}
(\pi^{\ast} r_{e_i}^{\vee}\otimes \mathrm{id})\boxtimes
\widetilde{l_{e_i^{\vee}}}$; here again we're using the local frames
$(e_1,\ldots ,e_{r+1})$, resp. $(e_1^{\vee},\ldots , e_{r+1}^{\vee})$, $r_{e_i} :
\mathcal{A}[-1] \to \mathcal{A}$ is right multiplication by $e_i$ and
$\widetilde{l_{e_i^{\vee}}}$ is the map induced by $l_{e_i^{\vee}}:
\mathcal{B}[-1] \to \mathcal{B}$ between the associated sheaves (via the
Serre correspondence).\\
This is truly a map of complexes since right and left Clifford multiplication
commute with each other. Moreover, we obtain a complex, infinite on
the left side
\[
R^{\bullet} : \,\ldots\to \Psi_i\boxtimes
\mathcal{O}(-i) \to \ldots \to\Psi_2\boxtimes \mathcal{O}(-2) \to \Psi_1\boxtimes
\mathcal{O}(-1) \to \mathcal{O}_{\mathcal{Q}\times_X\mathcal{Q}}\, .
\]
\begin{lemma}
The complex $R^{\bullet}$ is a left resolution of $\mathcal{O}_{\Delta}$,
$\Delta\subset\mathcal{Q}\times_X\mathcal{Q}$ being the diagonal.
\end{lemma}
\begin{proof}
Consider $\mathcal{B}^2:=\bigoplus_i
\mathcal{B}_i\otimes_{\mathcal{O}_X}\mathcal{B}_i$, the ``Segre product of
$\mathcal{B}$ with itself'' (i.e. the homogeneous coordinate ring of
$\mathcal{Q}\times_X\mathcal{Q}$ under the (relative) Segre morphism). Look at the
following double complex $D^{\bullet\bullet}$ of $\mathcal{B}^2$-modules:
\setlength{\unitlength}{1cm}
\begin{center}
\begin{picture}(5,4)

\put(-4,0){$\ldots$}
\put(-2.5,0){$\bigoplus_i \mathcal{A}_2^{\vee}\otimes\mathcal{B}_{i}\otimes
\mathcal{B}_{i-2}$}
\put(2,0){$\bigoplus_i \mathcal{A}_1^{\vee}\otimes\mathcal{B}_{i}\otimes
\mathcal{B}_{i-1}$}
\put(6.5,0){$\bigoplus_i \mathcal{B}_i\otimes\mathcal{B}_i$}

\put(-4,1.5){$\ldots$}
\put(-2.7,1.5){$\bigoplus_i \mathcal{A}_1^{\vee}\otimes\mathcal{B}_{i+1}\otimes
\mathcal{B}_{i-2}$}
\put(2,1.5){$\bigoplus_i \mathcal{B}_{i+1}\otimes\mathcal{B}_{i-1}$}

\put(-4,3){$\ldots$}
\put(-2.5,3){$\bigoplus_i \mathcal{B}_{i+2}\otimes\mathcal{B}_{i-2}$}

\put(-3.5,0.1){\vector(1,0){0.7}}
\put(-3.5,1.6){\vector(1,0){0.7}}
\put(-3.5,3.1){\vector(1,0){0.7}}

\put(1,0.1){\vector(1,0){0.7}}
\put(1,1.6){\vector(1,0){0.7}}

\put(5.5,0.1){\vector(1,0){0.7}}

\put(-1.5,0.5){\vector(0,1){0.7}}
\put(3,0.5){\vector(0,1){0.7}}

\put(-1.5,2){\vector(0,1){0.7}}
\end{picture}\end{center}

Here the columns correspond to the right resolutions of
$\Psi_0\boxtimes\mathcal{O}$, $\Psi_1\boxtimes \mathcal{O}(-1)$, $\Psi_2\boxtimes
\mathcal{O}(-2)$ etc. (starting from the right) if we pass from complexes of
coherent sheaves on $\mathcal{Q}\times_X\mathcal{Q}$ to complexes of graded
$\mathcal{B}^2$-modules via Serre's theorem. For example, the left-most column in
the above diagram arises from
\begin{eqnarray*}
(\pi^{\ast}\mathcal{A}_2^{\vee}\otimes \mathcal{O}_{\mathcal{Q}})\boxtimes
\mathcal{O}(-2)
\to
(\pi^{\ast}\mathcal{A}_{1}^{\vee} \otimes \mathcal{O}_{\mathcal{Q}}(1))\boxtimes
\mathcal{O}(-2) \\
\to
(\pi^{\ast}\mathcal{A}_0^{\vee} \otimes \mathcal{O}_{\mathcal{Q}}(2))\boxtimes
\mathcal{O}(-2)
\to 0
\end{eqnarray*}
The horizontal arrows in the above diagram then come from the morphisms of
complexes defining the differentials in $R^{\bullet}$.\\
The associated total complex $Tot^{\bullet}(D^{\bullet\bullet})$ has a natural
augmentation $a :Tot^{\bullet}(D^{\bullet\bullet}) \to \bigoplus_{i}
\mathcal{B}_{2i}$ arising from the multiplication maps
$\mathcal{B}_{i+j}\otimes\mathcal{B}_{i-j}
\to\mathcal{B}_{2i}$ and corresponding to the augmentation
$R^{\bullet}\to\mathcal{O}_{\Delta}$.\\ 
Claim: $a$ is a quasi-isomorphism. 
For this note that $D^{\bullet\bullet}$ is the direct sum over $i$ of double complexes
\setlength{\unitlength}{1cm}
\begin{center}
\begin{picture}(5,4)

\put(-4,0){$\ldots$}
\put(-2.5,0){$\mathcal{A}_2^{\vee}\otimes\mathcal{B}_{i}\otimes
\mathcal{B}_{i-2}$}
\put(2,0){$\mathcal{A}_1^{\vee}\otimes\mathcal{B}_{i}\otimes
\mathcal{B}_{i-1}$}
\put(6.5,0){$\mathcal{B}_i\otimes\mathcal{B}_i$}

\put(-4,1.5){$\ldots$}
\put(-2.7,1.5){$\mathcal{A}_1^{\vee}\otimes\mathcal{B}_{i+1}\otimes
\mathcal{B}_{i-2}$}
\put(2,1.5){$\mathcal{B}_{i+1}\otimes\mathcal{B}_{i-1}$}

\put(-4,3){$\ldots$}
\put(-2.5,3){$\mathcal{B}_{i+2}\otimes\mathcal{B}_{i-2}$}

\put(-3.5,0.1){\vector(1,0){0.7}}
\put(-3.5,1.6){\vector(1,0){0.7}}
\put(-3.5,3.1){\vector(1,0){0.7}}

\put(1,0.1){\vector(1,0){0.7}}
\put(1,1.6){\vector(1,0){0.7}}

\put(5.5,0.1){\vector(1,0){0.7}}

\put(-1.5,0.5){\vector(0,1){0.7}}
\put(3,0.5){\vector(0,1){0.7}}

\put(-1.5,2){\vector(0,1){0.7}}
\end{picture}\end{center}
which are bounded ($\mathcal{B}$ is positively graded) and whose rows are just
Priddy's resolution $P^{\bullet}$ in various degrees and thus the total complex of
the above direct summand of $D^{\bullet\bullet}$ is quasi-isomorphic to
$\mathcal{A}_0^{\vee}\otimes\mathcal{B}_{2i}\otimes
\mathcal{B}_0=\mathcal{B}_{2i}$. Thus $Tot^{\bullet}(D^{\bullet\bullet})$ is
quasi-isomorphic to $\bigoplus_i \mathcal{B}_{2i}$.
\end{proof}
The next step is to identify the kernel of the map $\Psi_{r-2}\boxtimes
\mathcal{O}(-r+2) \to \Psi_{r-3}\boxtimes
\mathcal{O}(-r+3)$. For this we have to talk in more detail about spinor bundles.\\
Let $\mathrm{Cliff}(\mathcal{E})=\mathcal{A}/(h-1)\mathcal{A}$ be the Clifford bundle of the orthogonal vector
bundle $\mathcal{E}$. This is just $\mathrm{Cliff}(\mathcal{E}):=
\mathrm{T}^{\bullet}\,\mathcal{E}/ I(\mathcal{E})$ where $I(\mathcal{E})$ is the
bundle of ideals whose fibre at $x\in X$ is the two-sided ideal
$I(\mathcal{E}(x))$ in $T^{\bullet}(\mathcal{E}(x))$ generated by the elements
$e\otimes e -q(e) 1$ for $e\in \mathcal{E}(x)$. 
$\mathrm{Cliff} (\mathcal{E})$ inherits a $\mathbb{Z}/2$-grading,
$\mathrm{Cliff}(\mathcal{E})=\mathrm{Cliff}^{\mathrm{even}}(\mathcal{E})\oplus
\mathrm{Cliff}^{\mathrm{odd}}(\mathcal{E})$.\\
Let us now make the assumption 
\begin{quote}
(A 1) \quad $H^1(X ; \mathbb{Z}/2\mathbb{Z})=0$.
\end{quote}
E. g., this will hold if $X$ is simply connected. 
Consider the bundle $P_{\mathrm{O}}(\mathcal{E})$ of orthonormal
frames in $\mathcal{E}$. This is the principal
$\mathrm{O}_{r+1}\,\mathbb{C}$-bundle whose fibre at a point $x\in X$
is the set of orthonormal bases of the fibre $\mathcal{E}(x)$. Choose
a principal $\mathrm{SO}_{r+1}\,\mathbb{C}$-subbundle
$P_{\mathrm{SO}}(\mathcal{E})\subset
P_{\mathrm{O}}(\mathcal{E})$. This is possible by (A 1) since then
$H^1(X ; \mathbb{Z}/ 2\mathbb{Z})=0$ and the exact sequence
\[
1\to \mathrm{SO}_{r+1}\, \mathbb{C} \to
\mathrm{O}_{r+1}\,\mathbb{C}\to \mathbb{Z}/2 \mathbb{Z}\to 0
\]
yields the exact sequence of cohomology sets with distinguished
elements
\[
H^0(X; \mathbb{Z}/2\mathbb{Z} )\to H^1 (X ;
\mathrm{SO}_{r+1}\,\mathbb{C}) \to H^1(X ;
\mathrm{O}_{r+1}\,\mathbb{C}) \to H^1(X ; \mathbb{Z}/2\mathbb{Z}) 
\]
where the \v{C}ech cohomology sets $H^1 (X ;
\mathrm{SO}_{r+1}\,\mathbb{C})$ resp. $H^1 (X ;
\mathrm{O}_{r+1}\,\mathbb{C})$ parametrize equivalence classes of
principal $\mathrm{SO}_{r+1}\,\mathbb{C}$-
resp. $\mathrm{O}_{r+1}\,\mathbb{C}$-bundles on $X$.\\
The short exact sequence
\[
0 \to \mathbb{Z}/2 \mathbb{Z} \to \mathrm{Spin}_{r+1}\,\mathbb{C} \to
\mathrm{SO}_{r+1}\, \mathbb{C} \to 1
\]
gives an exact sequence
\[
H^1 (X ; \mathbb{Z}/ 2\mathbb{Z}) \to H^1 (X ;
\mathrm{Spin}_{r+1}\,\mathbb{C}) \to H^1(X ; \mathrm{SO}_{r+1}\,
\mathbb{C}) \stackrel{\delta}{\longrightarrow} H^2 (X ; \mathbb{Z}/
2\mathbb{Z}) .
\]
We make the
additional assumption that
\begin{quote}
(A 2) \quad $\delta \left( [P_{\mathrm{SO}}(\mathcal{E})]  \right)=0$.
\end{quote}
This assumption just means that $\mathcal{E}$ carries a spin
structure, i.e. that $P_{\mathrm{SO}}(\mathcal{E})$ is the
$\mathbb{Z}/2 \mathbb{Z}$-quotient of a principal
$\mathrm{Spin}_{r+1}\,\mathbb{C}$-bundle
$P_{\mathrm{Spin}}(\mathcal{E})$ on $X$. The spin lifting
$P_{\mathrm{Spin}}(\mathcal{E})$ of $P_{\mathrm{SO}}(\mathcal{E})$ is
unique under assumption (A 1) since $H^1 (X ; \mathbb{Z}/
2\mathbb{Z})=0$.\\
In fact assumptions (A 1) and (A 2) will be automatically satisfied in
the applications to rational homogeneous manifolds in subsection 3.4
below, but in the abstract setting one has to make them.\\
If $r+1$ is odd one has the spin representation $S$ of
$\mathrm{Spin}_{r+1}\,\mathbb{C}$
\[
\varrho_S : \mathrm{Spin}_{r+1}\,\mathbb{C} \to \mathrm{Aut}\, S
\]
with $\dim S=2^{r/2}$ (for the description of
$\mathrm{Spin}_{r+1}\,\mathbb{C}$ as a closed subgroup of the group of
units in the even Clifford algebra of an $r+1$-dimensional orthogonal
vector space and the resulting classical construction of $S$ via the
identification of the even Clifford algebra with the algebra of linear
endomorphisms of the exterior algebra of a maximal isotropic subspace cf. \cite{FuHa}, \S 20); we have on $X$ the spinor
bundle $S (\mathcal{E})$ of $\mathcal{E}$:
\begin{gather*}
S(\mathcal{E}):=P_{\mathrm{Spin}}(\mathcal{E}) \times_{\varrho_S} S\\
(= P_{\mathrm{Spin}}(\mathcal{E}) \times S / \{ (p ,s )\sim (pg^{-1} ,
\varrho_S (g) s) , \; p\in P_{\mathrm{Spin}}(\mathcal{E}) ,\; g\in
\mathrm{Spin}_{r+1}\,\mathbb{C} , \; s\in S    \} .) 
\end{gather*}
For $r+1$ even one has the two half spin representations $S^{\pm}$
with $\dim S^{\pm} = 2^{\frac{r+1}{2} -1}$: 
\[
\varrho_{S^{\pm}} : \mathrm{Spin}_{r+1}\,\mathbb{C} \to \mathrm{Aut}\, S^{\pm}
\]
(identifying $\mathrm{Spin}_{r+1}\,\mathbb{C}$ with a closed subgroup of the group of
units in the even Clifford algebra of an $r+1$-dimensional orthogonal
vector space, we have in this case that the even Clifford algebra
splits as a direct sum: One summand is the algebra of linear
endomorphisms of the space consisting of all the even exterior powers
of a maximal isotropic subspace; the other summand is the algebra of linear
endomorphisms of the space consisting of all the odd exterior powers
of this maximal isotropic subspace. We refer again to \cite{FuHa}, \S
20); on $X$ we have the
 associated spinor bundles $S^{\pm}(\mathcal{E})$ of the
orthogonal vector bundle $\mathcal{E}$, i.e. $S^{\pm}(\mathcal{E}):=
P_{\mathrm{Spin}} (\mathcal{E}) \times_{\varrho_{S^{\pm}}} S^{\pm}$.\\
For \fbox{$r+1$ even} $\mathrm{Cliff} (\mathcal{E}) \simeq
\mathcal{E}nd ( S^{+} (\mathcal{E}) \oplus S^{-}(\mathcal{E}) )$
($S^{+} (\mathcal{E}) \oplus S^{-}(\mathcal{E})$ is a bundle of
irreducible $\mathrm{Cliff} (\mathcal{E})$-modules) and 
\[
\mathrm{Cliff}^{\mathrm{even}}(\mathcal{E}) \simeq \mathcal{E}nd
(S^{+}(\mathcal{E})) \oplus
\mathcal{E}nd (S^{-}(\mathcal{E}))\, .
\]
Then
\[
\mathcal{M}^{-}:=S^{-}(\mathcal{E}) \oplus S^{+}(\mathcal{E}) \oplus
S^{-}(\mathcal{E})\oplus \ldots
\]
and
\[
\mathcal{M}^{+}:=S^{+}(\mathcal{E}) \oplus S^{-}(\mathcal{E}) \oplus
S^{+}(\mathcal{E})\oplus \ldots
\]
are graded left $\mathcal{A}$-modules (the grading starting from $0$); one
defines bundles $\Sigma^{+}$, $\Sigma^{-}$ on $\mathcal{Q}$ by the requirement
that
\begin{eqnarray*}
0\to (\Sigma^{\pm})^{\vee}\to L^{\bullet} (\mathcal{M}^{\pm})\;\;
\mathrm{for}\;\; \fbox{$r+1\equiv 0 (\mathrm{mod}\, 4)$} \, ,\\
0\to (\Sigma^{\mp})^{\vee}\to L^{\bullet} (\mathcal{M}^{\pm})\; \;\mathrm{for}\;
\; \fbox{$r+1\equiv 2 (\mathrm{mod}\, 4)$}
\end{eqnarray*}
be exact. \\
For \fbox{$r+1$ odd} we have $\mathrm{Cliff} (\mathcal{E})\simeq
\mathrm{E}nd (S (\mathcal{E})) \oplus \mathrm{E}nd (S (\mathcal{E})')$
with $S(\mathcal{E})'\simeq S(\mathcal{E})$, and
$\mathrm{Cliff}^{\mathrm{even}}(\mathcal{E})\simeq \mathcal{E}nd
(S(\mathcal{E} ))$. Let $\mathcal{M}$ be the graded left
$\mathcal{A}$- module (grading starting from $0$)
\[
\mathcal{M}:= S (\mathcal{E})\oplus S (\mathcal{E})\oplus S
(\mathcal{E})\oplus\ldots
\]
and define the bundle $\Sigma$ on $\mathcal{Q}$ by the requirement that
\[
0\to (\Sigma )^{\vee} \to L^{\bullet} (\mathcal{M})
\]
be exact.\\
From the definition $\Sigma^{\pm}=:\Sigma^{\pm}(\mathcal{O}_{\mathcal{Q}}(-1)^{\perp}/\mathcal{O}_{\mathcal{Q}}(-1))$
resp. $\Sigma =:\Sigma
(\mathcal{O}_{\mathcal{Q}}(-1)^{\perp}/\mathcal{O}_{\mathcal{Q}}(-1))$
are the duals of the spinor bundles associated to the orthogonal
vector bundle $\mathcal{O}_{\mathcal{Q}}(-1)^{\perp}/\mathcal{O}_{\mathcal{Q}}(-1)$ on
$\mathcal{Q}$, but (slightly abusing the language) we will refer to
them just as spinor bundles in the sequel.
\begin{lemma}
\begin{gather*}
\mathrm{ker} (\Psi_{r-2}\boxtimes\mathcal{O}(-r+2) \to
\Psi_{r-3}\boxtimes\mathcal{O}(-r+3) )\\
=\left\{ \begin{array}{r@{, \quad}l}
\Sigma (-1)\boxtimes \Sigma (-r+1) & r+1 \;\; \mathrm{odd} \\
(\Sigma^{+}(-1)\boxtimes\Sigma^{+}(-r+1))\oplus
(\Sigma^{-}(-1)\boxtimes\Sigma^{-}(-r+1)) & r+1\equiv 2 (\mathrm{mod}\, 4) \\
(\Sigma^{+}(-1)\boxtimes\Sigma^{-}(-r+1))\oplus
(\Sigma^{-}(-1)\boxtimes\Sigma^{+}(-r+1)) & r+1\equiv 0 (\mathrm{mod}\, 4) \\
\end{array} \right.
\end{gather*}
\end{lemma}
\begin{proof}
For $i\ge r$
$\mathcal{A}_i\stackrel{\mathrm{mult}(h)}{\longrightarrow}\mathcal{A}_{i+2}$ is
an isomorphism because $(e_{i_1}\cdot\ldots \cdot e_{i_k} h^m)$, $1\le i_1<
\ldots <i_k\le r+1$, $m\in \mathbb{N}$ is a local frame for $\mathcal{A}$ if
$(e_1,\ldots , e_{r+1})$ is one for $\mathcal{E}$, and the map $\mathcal{A}_i\to
\mathrm{Cliff}^{\mathrm{par}(i)}(\mathcal{E})$ induced by
$\mathcal{A}\to \mathcal{A}/(h-1)\mathcal{A}$ is then an isomorphism where
\[
\mathrm{par}(i):= \left\{ \begin{array}{r@{, \quad}l}
\mathrm{even} & i\equiv 0 (\mathrm{mod}\, 2) \\
\mathrm{odd} & i\equiv 1 (\mathrm{mod}\, 2)
\end{array} \right. .
\]
Because $L^{\bullet}(\mathcal{A}^{\vee})$ is exact, $\Psi_i$ is also the cokernel
of
\[ (\ast ) \quad
\ldots \to \pi^{\ast}\mathcal{A}_{i+3}^{\vee}\otimes
\mathcal{O}_{\mathcal{Q}}(-3) \to \pi^{\ast}\mathcal{A}_{i+2}^{\vee}\otimes
\mathcal{O}_{\mathcal{Q}}(-2) \to \pi^{\ast}\mathcal{A}_{i+1}^{\vee}\otimes
\mathcal{O}_{\mathcal{Q}}(-1)
\]
Since $\mathrm{ker} (\Psi_{r-2}\boxtimes\mathcal{O}(-r+2) \to \Psi_{r-3}\boxtimes
\mathcal{O}(-r+3) )= \mathrm{coker} (\Psi_{r}\boxtimes\mathcal{O}(-r) \to
\Psi_{r-1}\boxtimes
\mathcal{O}(-r+1) )$ we conclude that a left resolution of the kernel in lemma
3.2.4 is given by $Tot^{\bullet}(E^{\bullet\bullet})$ where $E^{\bullet\bullet}$
is the following double complex:
\small
\setlength{\unitlength}{1cm}
\begin{center}
\begin{picture}(5,3.5)

\put(-4,0){$\ldots$}
\put(-3,0){$(\pi^{\ast}\mathrm{Cliff}^{\mathrm{par}(r+1)}\mathcal{E}^{\vee}(-1))
\boxtimes\mathcal{O}(-r)$}
\put(3,0){$(\pi^{\ast}\mathrm{Cliff}^{\mathrm{par}(r)}\mathcal{E}^{\vee}
(-1))\boxtimes\mathcal{O}(-r+1)$}

\put(-4,1.5){$\ldots$}

\put(-3,1.5){$(\pi^{\ast}\mathrm{Cliff}^{\mathrm{par}(r+2)}\mathcal{E}^{\vee}
(-2))\boxtimes\mathcal{O}(-r)$}

\put(3,1.5){$(\pi^{\ast}\mathrm{Cliff}^{\mathrm{par}(r+1)}\mathcal{E}^{\vee}(-2))
\boxtimes\mathcal{O}(-r+1)$}

\put(-0.5,3){$\vdots$}
\put(5.5,3){$\vdots$}

\put(-3.5,0.1){\vector(1,0){0.3}}
\put(-3.5,1.6){\vector(1,0){0.3}}

\put(2.2,0.1){\vector(1,0){0.7}}
\put(2.2,1.6){\vector(1,0){0.7}}

\put(-0.5,1.2){\vector(0,-1){0.7}}
\put(5.5,1.2){\vector(0,-1){0.7}}
\put(5.5,2.7){\vector(0,-1){0.7}}

\put(-0.5,2.7){\vector(0,-1){0.7}}
\end{picture}\end{center}
\normalsize
\vspace{0.3cm}
Here the columns (starting from the right) are the left resolutions ($\ast$) of
$\Psi_{r-1}\boxtimes \mathcal{O}(-r+1)$, $\Psi_r\boxtimes\mathcal{O}(-r)$, etc.
and the rows are defined through the morphisms of complexes defining the
differentials $\Psi_{r-1}\boxtimes\mathcal{O}(-r+1) \to\Psi_r\boxtimes
\mathcal{O}(-r)$ etc. in the resolution $R^{\bullet}$. For odd $r+1$ we have
$\mathrm{Cliff}^{odd}(\mathcal{E})\simeq
\mathrm{Cliff}^{even}(\mathcal{E})\simeq \mathcal{E}nd(S
(\mathcal{E}))\simeq S (\mathcal{E})^{\vee}\otimes S (\mathcal{E})$
whence our double complex becomes
\small
\setlength{\unitlength}{1cm}
\begin{center}
\begin{picture}(5,3.5)

\put(-4,0){$\ldots$}

\put(-2.5,0){$\pi^{\ast}S (\mathcal{E})^{\vee}(-1)
\boxtimes\pi^{\ast}S (\mathcal{E})(-r)$}
\put(3,0){$\pi^{\ast}S (\mathcal{E})^{\vee}(-1)
\boxtimes\pi^{\ast}S (\mathcal{E})(-r+1)$}

\put(-4,1.5){$\ldots$}

\put(-2.5,1.5){$\pi^{\ast} S (\mathcal{E})^{\vee}(-2)
\boxtimes\pi^{\ast} S (\mathcal{E})(-r)$}

\put(3,1.5){$\pi^{\ast}   S   (\mathcal{E})^{\vee}(-2)
\boxtimes\pi^{\ast}   S   (\mathcal{E})(-r+1)$}

\put(0,3){$\vdots$}
\put(5.5,3){$\vdots$}

\put(-3.5,0.1){\vector(1,0){0.7}}
\put(-3.5,1.6){\vector(1,0){0.7}}

\put(2,0.1){\vector(1,0){0.7}}
\put(2,1.6){\vector(1,0){0.7}}

\put(0,1.2){\vector(0,-1){0.7}}
\put(5.5,1.2){\vector(0,-1){0.7}}
\put(5.5,2.7){\vector(0,-1){0.7}}

\put(0,2.7){\vector(0,-1){0.7}}
\end{picture}\end{center}
\normalsize
\vspace{0.3cm}
and is thus isomorphic as a double complex to
$L^{\bullet}(\mathcal{M})^{\vee}(-1)\boxtimes
L^{\bullet}(\mathcal{M})^{\vee}(-r+1)$, i.e. quasi-isomorphic to $\Sigma
(-1)\boxtimes \Sigma (-r+1)$. The cases for even $r+1$ are considered similarly.
\end{proof}

\begin{lemma}
Consider the following two ordered sets of sheaves on $\mathcal{Q}$:
\begin{gather*}
\mathfrak{S}=\left\{ \Sigma (-r+1) \prec \mathcal{O}_{\mathcal{Q}}(-r+2)\prec
\ldots \prec
\mathcal{O}_{\mathcal{Q}}(-1)\prec \mathcal{O}_{\mathcal{Q}} \right\} \quad (
r+1\;
\mathrm{odd} )\, , \\
\mathfrak{S}'=\left\{ \Sigma^{+} (-r+1)\prec \Sigma^{-} (-r+1) \prec
\ldots \prec
\mathcal{O}_{\mathcal{Q}}(-1)\prec \mathcal{O}_{\mathcal{Q}} \right\} \quad (
r+1\;
\mathrm{even} ) .
\end{gather*}
If $\mathcal{V}, \,\mathcal{V}_1 ,\, \mathcal{V}_2 \in \mathfrak{S}$ (resp.: $\in
\mathfrak{S}'$) with $\mathcal{V}_1 \prec \mathcal{V}_2$, $\mathcal{V}_1\neq
\mathcal{V}_2$, we have the following identities
\begin{eqnarray*}
R^i\pi_{\ast} (\mathcal{V}\otimes\mathcal{V}^{\vee})=0\, , \quad \forall i\neq 0 ,
\quad  \\
R^i\pi_{\ast} (\mathcal{V}_1\otimes\mathcal{V}_2^{\vee})=0 \; \forall
i\in\mathbb{Z} , \quad R^i\pi_{\ast}
(\mathcal{V}_2\otimes\mathcal{V}_1^{\vee})=0\; \forall i \neq 0 \, .
\end{eqnarray*}
and the canonical morphism $R^0\pi_{\ast}
(\mathcal{V}\otimes\mathcal{V}^{\vee})\to \mathcal{O}_X$ is an isomorphism.
\end{lemma}
\begin{proof}
In the absolute case (where the base $X$ is a point) this is a calculation in
\cite{Ka3} , prop. 4.9., based on Bott's theorem. The general assertion follows
from this because the question is local on $X$ and we can check this on 
affine open sets $U\subset X$ which cover $X$ and over which $\mathcal{Q}$ is trivial using
$H^q(\pi^{-1}(U), \mathcal{F})\simeq \Gamma (U, R^q\pi_{\ast}(\mathcal{F}))$ for
every coherent $\mathcal{F}$ on $\pi^{-1}(U)$ and the K\"unneth formula.
\end{proof}

As in subsection 3.1, for $\mathcal{V}\in\mathfrak{S}$ (resp. 
$\in\mathfrak{S}'$), we can talk about subcategories $D^b(X)\otimes \mathcal{V}$
of $D^b(\mathcal{Q})$ as the images of $D^b(X)$ in $D^b(\mathcal{Q})$ under the
functor $\pi^{\ast}(-)\otimes \mathcal{V}$.

\begin{proposition}
Let $\mathcal{V},\: \mathcal{V}_1 ,\mathcal{V}_2$ be as in lemma 3.2.5. The
subcategories $D^b(X)\otimes \mathcal{V}$ of $D^b(\mathcal{Q})$ are all
admissible subcategories. Moreover, for $A\in
\mathrm{obj} (D^b(X)\otimes\mathcal{V}_2)$, $B\in
\mathrm{obj} (D^b(X)\otimes\mathcal{V}_1)$ we have $R\mathrm{Hom}(A, B)=0$.
\end{proposition}
\begin{proof}
Let $A=\pi^{\ast} A'\otimes \mathcal{V}_2$, $B=\pi^{\ast} B'\otimes
\mathcal{V}_1$. Using lemma 3.2.5 and the projection formula we compute
\begin{eqnarray*}
R^i\mathrm{Hom}(\pi^{\ast}A'\otimes \mathcal{V}_2,
\pi^{\ast}B'\otimes \mathcal{V}_1)
\simeq R^i\mathrm{Hom}(\pi^{\ast} A',
\pi^{\ast} B'\otimes \mathcal{V}_1\otimes\mathcal{V}_2^{\vee})\\
\simeq R^i\mathrm{Hom}(A',
B'\otimes R\pi_{\ast}(\mathcal{V}_1\otimes\mathcal{V}_2^{\vee})) \simeq 0\, .
\end{eqnarray*}
If we repeat the same calculation with $\mathcal{V}$ instead of $\mathcal{V}_1$
and $\mathcal{V}_2$ we find that $R^i\mathrm{Hom}(\pi^{\ast}A'\otimes
\mathcal{V},\pi^{\ast}B'\otimes \mathcal{V})\simeq R^i\mathrm{Hom}(A',B')$. This
shows that the categories $D^b(X)\otimes \mathcal{V}$ are all equivalent to
$D^b(X)$ as triangulated subcategories of $D^b(\mathcal{Q})$. It follows from
\cite{BoKa} , prop. 2.6 and thm. 2.14, together with lemma 3.2.1 that the
$D^b(X)\otimes\mathcal{V}$ are admissible subcategories of $D^b(\mathcal{Q})$.
\end{proof}

\begin{theorem}
Let $X$ be a smooth projective variety, $\mathcal{E}$ an orthogonal
vector bundle on $X$, $\mathcal{Q}\subset \mathbb{P}(\mathcal{E})$ the
associated quadric bundle, and let assumptions (A 1) and (A 2) above
be satisfied.\\
Then there is a semiorthogonal decomposition
\begin{gather*}
D^b(\mathcal{Q})=\left\langle  D^b(X)\otimes \Sigma (-r+1) , D^b(X)\otimes
\mathcal{O}_{\mathcal{Q}}(-r+2) ,\right. \\
\left. \ldots , D^b(X)\otimes \mathcal{O}_{\mathcal{Q}}(-1) , D^b(X) 
\right\rangle
\end{gather*}
for $r+1$ odd and
\begin{gather*}
D^b(\mathcal{Q})=\left\langle  D^b(X)\otimes \Sigma^{+}(-r+1) , D^b(X)\otimes
\Sigma^{-}(-r+1) ,\right. \\ \left. D^b(X)\otimes
\mathcal{O}_{\mathcal{Q}}(-r+2) , \ldots , D^b(X)\otimes
\mathcal{O}_{\mathcal{Q}}(-1) , D^b(X)
\right\rangle 
\end{gather*}
for $r+1$ even.
\end{theorem}
\begin{proof}
By proposition 3.2.6 the categories in question are semiorthogonal and it remains
to see that they generate $D^b(\mathcal{Q})$. For ease of notation we will
consider the case of odd $r+1$, the case of even $r+1$ being entirely similar.\\
From lemmas 3.2.3 and 3.2.4 we know that in the situation of the fibre product
\[
\begin{CD}
\Delta\subset \mathcal{Q}\times_X \mathcal{Q} @>p_2>> \mathcal{Q} \\
@Vp_1VV                                  @V\pi VV  \\
\mathcal{Q} @>\pi >>              X
\end{CD}
\]
we have a resolution
\begin{eqnarray*}
0\to \Sigma (-1) \boxtimes \Sigma (-r+1) \to
\Psi_{r-2}\boxtimes\mathcal{O}_{\mathcal{Q}}(-r+2) \to \ldots \\
\ldots \to
\Psi_1\boxtimes\mathcal{O}(-1) \to \mathcal{O}_{\mathcal{Q}\times_X\mathcal{Q}}
\to \mathcal{O}_{\Delta}\to 0
\end{eqnarray*}
and tensoring this with $p_1^{\ast}\mathcal{F}$ ($\mathcal{F}$ a coherent sheaf
on $\mathcal{Q}$)
\begin{eqnarray*}
0\to (\Sigma (-1)\otimes \mathcal{F}) \boxtimes \Sigma (-r+1) \to
(\Psi_{r-2}\otimes\mathcal{F})\boxtimes\mathcal{O}_{\mathcal{Q}}(-r+2) \to \ldots
\\
\ldots \to
(\Psi_1\otimes\mathcal{F})\boxtimes\mathcal{O}(-1) \to
\mathcal{F}\boxtimes\mathcal{O}_{\mathcal{Q}}
\to p_1^{\ast}\mathcal{F}|_{\Delta}\to 0
\end{eqnarray*}
and applying $R p_{2\,\ast}$ we obtain a spectral sequence
\begin{gather*}
E_1^{ij}=R^ip_{2\, \ast}( (\Psi_{-j}\otimes \mathcal{F})\boxtimes
\mathcal{O}_{\mathcal{Q}}(j))\quad -r+2<j\le 0\\
\quad \quad = R^ip_{2\, \ast}( (\Sigma (-1)\otimes \mathcal{F})\boxtimes
\Sigma (-r+1))\quad j=-r+1\\
=0 \quad \mathrm{otherwise}\hspace{5cm}
\end{gather*}
and $E_1^{ij} \Rightarrow R^{i+j}p_{2\, \ast} (p_1^{\ast}\mathcal{F}|_{\Delta})$
which is $=\mathcal{F}$ for $i+j=0$ and $=0$ otherwise. But since cohomology
commutes with flat base extension (cf. \cite{EGA}, III, \S 1, prop. 1.4.15), we
have
$R^ip_{2\,
\ast} p_1^{\ast}
\mathcal{G}\simeq \pi^{\ast} R^i\pi_{\ast} \mathcal{G}$ for any coherent
$\mathcal{G}$ on $\mathcal{Q}$. This together with the projection formula shows
that all $E^{ij}_1$ belong to one of the admissible subcategories in the
statement of theorem 3.2.7. This finishes the proof because $D^b(\mathcal{Q})$ is
generated by the subcategory $Coh (\mathcal{Q})$.
\end{proof}

\begin{corollary}
If $D^b(X)$ is generated by a complete exceptional sequence
\[
(E_1,\ldots , E_n) \, ,
\]
then $D^b(\mathcal{Q})$ is generated by the complete exceptional sequence
\[
\left( \pi^{\ast}E_1\otimes \Sigma (-r+1), \ldots , \pi^{\ast}E_n\otimes \Sigma
(-r+1),
\pi^{\ast}E_1\otimes \mathcal{O}_{\mathcal{Q}}(-r+2), \ldots , \pi^{\ast}
E_n\right)
\]
for $r+1$ odd and 
\begin{gather*}
\left( \pi^{\ast}E_1\otimes \Sigma^{+} (-r+1), \ldots , \pi^{\ast}E_n\otimes
\Sigma^{+} (-r+1),\ldots , \pi^{\ast}E_1\otimes \Sigma^{-} (-r+1), \right.\\
\left. \ldots ,
\pi^{\ast}E_n\otimes
\Sigma^{-} (-r+1),\pi^{\ast}E_1\otimes \mathcal{O}_{\mathcal{Q}}(-r+2), \ldots ,
\pi^{\ast} E_1 , \ldots , \pi^{\ast} E_n
\right) 
\end{gather*}
for $r+1$ even.
\end{corollary}
\begin{proof}
Using lemma 3.2.5, one proves this analogously to corollary 3.1.2; we omit the
details. 
\end{proof}

\subsection{Application to varieties of isotropic flags in a
  symplectic vector space}

We first fix some notation: Let $V$ be a $\mathbb{C}$-vector space of even dimension $2n$
 with a nondegenerate skew symmetric bilinear form $\langle \cdot , \cdot \rangle $. For 
$1\le k_1 < \ldots < k_t\le n$ we denote $\mathrm{IFlag}(k_1,\ldots ,k_t ; V):= \{ (L_{k_1},
\ldots , L_{k_t})\,\mid \, L_{k_1}\subset \ldots L_{k_t} \subset V$ isotropic subspaces of 
$V$ with $\dim L_{k_j}= k_j , \: 1\le j\le t\}$ the (partial) flag variety of isotropic flags of
 type $(k_1,\ldots ,k_t)$ in $V$; 
moreover, for $1\le k\le n$, put $\mathrm{IGrass}(k , V):=\mathrm{IFlag}(k ; V)$, the  
Grassmann manifold of isotropic $k$-planes in $V$. As usual, we have the tautological flag of
subbundles $\mathcal{R}_{k_1}\subset \ldots \subset \mathcal{R}_{k_t}\subset V\otimes 
\mathcal{O}_{\mathrm{IFlag}(k_1,\ldots ,k_t ; V)}$ on $\mathrm{IFlag}(k_1,\ldots ,k_t ; V)$ and 
the tautological subbundle $\mathcal{R}$ on $\mathrm{IGrass}(k , V)$.
\begin{remark}
Via the projection $\mathrm{IFlag}(k_1,\ldots ,k_t ; V) \to \mathrm{IGrass} (k_t , V)$, the 
variety $\mathrm{IFlag}(k_1,\ldots ,k_t ; V)$ identifies with $\mathrm{Flag}_{\mathrm{IGrass}
 (k_t , V)} (k_1 ,\ldots , k_{t-1} ; \mathcal{R} )$, the relative variety of flags of type 
$(k_1,\ldots ,k_{t-1})$ in the fibres of the tautological subbundle $\mathcal{R}$ on 
$\mathrm{IGrass} (k_t , V)$. Therefore, by theorem 3.1.5, if we want to exhibit complete 
exceptional sequences in the derived categories of coherent sheaves on all possible varieties
of (partial) isotropic flags in $V$, we can reduce to finding them on isotropic Grassmannians. 
Thus we will focus on the latter in the sequel. 
\end{remark}

Now look at the following diagram (the notation will be explained below)
\vspace{0.4cm}
\setlength{\unitlength}{1cm}
\begin{center}
\begin{picture}(2,5.8)
\put(-4,0){$\mathrm{IFlag}(1 ; V)\simeq \mathbb{P}(V)$}
\put(-4,1.5){$\mathrm{IFlag}(1,2 ; V)\simeq \mathbb{P}(\mathcal{E}_1)$}
\put(-4,4.3){$\mathrm{IFlag}(1,2,\ldots , k-1 ; V)\simeq \mathbb{P}(\mathcal{E}_{k-2})$}
\put(-4,5.8){$\mathrm{IFlag}(1,2,\ldots , k-1, k ; V)\simeq \mathbb{P}(\mathcal{E}_{k-1})
\simeq \mathrm{Flag}_{\mathrm{IGrass}(k , V)}(1, 2 ,\ldots , k ; \mathcal{R})$}
\put(-3,2.9){$\vdots$}

\put(-3,1.3){\vector(0, -1){0.8}}
\put(-3,5.6){\vector(0, -1){0.8}}

\put(-3,2.7){\vector(0, -1){0.8}}
\put(-3,4.1){\vector(0, -1){0.8}}

\put(4.5,0){$\mathrm{IGrass}(k , V)$}

\put(4,5.3){\vector(1, -3){1.6}}
\put(4,3){$\pi$}

\put(-2.8 , 0.8){$\pi_1$}
\put(-2.8 , 2.2){$\pi_2$}

\put(-2.8 , 3.6){$\pi_{k-2}$}
\put(-2.8 , 5.1){$\pi_{k-1}$}
\end{picture}\end{center}
\vspace{0.3cm}
Since for $1\le i \le j \le k-1$ the $i$-dimensional tautological subbundle on 
$\mathrm{IFlag}(1 , 2 ,\ldots , j ; V)$ pulls back to the $i$-dimensional tautological 
subbundle on  $\mathrm{IFlag}(1 , 2 ,\ldots , j+1 ; V)$ under $\pi_{j}$, we denote all 
of them by the same symbol $\mathcal{R}_i$ regardless of which space they live on, if no 
confusion can arise. \\
Since any line in $V$ is isotropic, the choice of a $1$-dimensional isotropic $L_1\subset V$ 
comes down to picking a point in $\mathbb{P}(V)$ whence the identification $\mathrm{IFlag}(1 ; V)
\simeq \mathbb{P}(V)$ above; the space $L_1^{\perp}/L_1$ is again a symplectic vector space with 
skew form induced from $\langle \cdot ,\cdot \rangle$ on $V$, and finding an isotropic plane 
containing $L_1$ amounts to choosing a line $L_2/L_1$ in $L_1^{\perp}/L_1$. Thus $\mathrm{IFlag}
(1,2 ; V)$ is a projective bundle $\mathbb{P}(\mathcal{E}_1)$ over $\mathrm{IFlag}(1; V)$ with 
$\mathcal{E}_1=\mathcal{R}_1^{\perp}/\mathcal{R}_1$, and on $\mathbb{P}(\mathcal{E}_1)\simeq 
\mathrm{IFlag}(1 ,2 ; V)$ we have 
$\mathcal{O}_{\mathcal{E}_1}(-1)\simeq \mathcal{R}_2/\mathcal{R}_1$. Of course, $\mathrm{rk}\, 
\mathcal{E}_1= 2n-2$.\\
Continuing this way, we successively build the whole tower of projective bundles over $\mathbb{P}
(V)$ in the above diagram where
\begin{eqnarray*}
\mathcal{E}_j \simeq \mathcal{R}_j^{\perp}/\mathcal{R}_j\, ,\quad j=1 , \ldots , k-1\, , \quad 
\mathrm{rk}\, \mathcal{E}_j=2n-2j \\
\mathrm{and}\quad \mathcal{O}_{\mathcal{E}_j}(-1) \simeq \mathcal{R}_{j+1}/\mathcal{R}_j\, .
\end{eqnarray*} 
Moreover, $\mathrm{IFlag}(1,2,\ldots , k-1, k ; V)$ is isomorphic to $\mathrm{Flag}_{\mathrm{IGrass}
(k , V)}(1, \ldots , k ; \mathcal{R})$, the relative variety of complete flags in the fibres 
of the tautological subbundle $\mathcal{R}$ on $\mathrm{IGrass} (k , V)$; the flag of tautological 
subbundles in $V\otimes\mathcal{O}_{\mathrm{IFlag}(1, \ldots ,k; V)}$ on $\mathrm{IFlag}(1, \ldots ,
k; V)$ and the flag of relative tautological subbundles in $\pi^{\ast}\mathcal{R}$ on $\mathrm{Flag}
_{\mathrm{IGrass}(k , V)}(1, 2 ,\ldots , k ; \mathcal{R})$ correspond to each other under this 
isomorphism, and we do not distinguish them notationally.\\
For $\lambda =(\lambda_1, \ldots , \lambda_k)\in\mathbb{Z}^k$ define the line bundle $\mathcal{L}
(\lambda )$ on the variety $\mathrm{Flag}
_{\mathrm{IGrass}(k , V)}(1,\ldots , k ; \mathcal{R})$ by 
\[
\mathcal{L}(\lambda ):= (\mathcal{R}_1)^{\otimes (-\lambda_1)} \otimes (\mathcal{R}_2/\mathcal{R}_1)
^{\otimes (-\lambda_2) } \otimes \ldots \otimes (\pi^{\ast}\mathcal{R} / \mathcal{R}_{k-1})^{\otimes 
(-\lambda_k)} \, .
\]
Repeatedly applying corollary 3.1.2 to the above tower of projective bundles, we find that the 
following sheaves constitute a complete exceptional sequence in $D^b(Coh \, (\mathrm{Flag}
_{\mathrm{IGrass}(k , V)}(1,\ldots , k ; \mathcal{R})) )$:
\begin{eqnarray*}
(\mathcal{L}(\lambda )) \quad \mathrm{with} \quad -2n+1 \le \lambda_1 \le 0 , \\
 (\sharp )\hspace{5cm}                                              -2n+3 \le \lambda_2 \le 0 ,\\
                                             \vdots \quad \quad \quad\\
                                                -2n +2k -1 \le \lambda_k \le 0 .
\end{eqnarray*}
Here $\mathcal{L}(\lambda )$ precedes $\mathcal{L}(\mu )$ according to the ordering of the 
exceptional sequence iff $(\lambda_k, \lambda_{k-1}, \ldots , \lambda_1)  \prec (\mu_k , 
\mu_{k-1}, \ldots , \mu_1)$ where $\prec$ is the lexicographic order on $\mathbb{Z}^k$.\\
Let us record here the simple
\begin{lemma}
The set of full direct images $R\pi_{\ast } \mathcal{L}(\lambda )$, as
 $\mathcal{L}(\lambda )$ varies  
among the bundles $(\sharp )$, generates the derived category $D^b (Coh (\mathrm{IGrass}(k, V)))$.
\end{lemma}
\begin{proof}
As in lemma 3.2.1, $R\pi_{\ast }\mathcal{O}_{\mathrm{Flag}
_{\mathrm{IGrass}(k , V)}(1,\ldots , k ; \mathcal{R})} \simeq \mathcal{O}_{\mathrm{IGrass}(k, V)}$, 
and $R\pi_{\ast }\circ \pi^{\ast}$ is isomorphic to the identity functor on $D^b (Coh (\mathrm{
IGrass}(k, V)))$ by the projection formula. Thus, since the bundles in $(\sharp )$ generate the 
derived category upstairs, if $E$ is an object in $D^b (Coh (\mathrm{
IGrass}(k, V)))$, $\pi^{\ast} E$ will be isomorphic to an object in the smallest full triangulated 
subcategory containing the objects $(\sharp)$, i.e. starting from the set $(\sharp )$ and 
repeatedly enlarging it by taking finite direct sums, shifting in cohomological degree and 
completing distinguished triangles by taking a mapping cone, we can reach an object isomorphic to 
$\pi^{\ast} E$. Hence it is clear that the objects $Rp_{\ast}\mathcal{L}(\lambda )$ will generate 
the derived category downstairs because $R\pi_{\ast}\pi^{\ast } E\simeq E$. 
\end{proof}
Now the fibre of $\mathrm{Flag}_{\mathrm{IGrass}(k, V)}(1,\ldots ,k; \mathcal{R})$ over a point 
$x\in \mathrm{IGrass}(k , V)$ is just the full flag variety $\mathrm{Flag}(1, \ldots , k ; 
\mathcal{R}(x))$ which is a quotient of $\mathrm{GL}_k \mathbb{C}$ by a Borel subgroup $B$; the 
$\lambda \in\mathbb{Z}^k$ can be identified with weights or characters of a maximal torus $H\subset B$ 
and the restriction of $\mathcal{L}(\lambda )$ to the fibre over $x$ is just the line bundle 
associated to the character $\lambda $, i.e. $\mathrm{GL}_k \mathbb{C} \times_B \mathbb{C}_{-\lambda}$, 
where $\mathbb{C}_{-\lambda }$ is the one-dimensional $B$-module in which the torus $H$ acts via the 
character $-\lambda$ and the unipotent radical $R_u(B)$ of $B$ acts trivially, and $\mathrm{GL}_k \mathbb{C}
 \times_B \mathbb{C}_{-\lambda }:=\mathrm{GL}_k \mathbb{C} \times \mathbb{C}_{-\lambda }/\{ (g,v) \sim 
(gb^{-1} , b v)\, , \: b\in B \}  $. Thus we can calculate the $R\pi_{\ast} \mathcal{L}(\lambda )$ by 
the following (relative) version of Bott's theorem (cf. \cite{Wey}, thm. 4.1.4 or \cite{Akh}, \S 4.3 
for a full statement):\\

Let $\varrho :=(k-1 , k-2 , \ldots , 0)$ (the half sum of the positive roots) and let $W=\mathfrak{S}
_k$, the symmetric group on $k$ letters (the Weyl group), act on $\mathbb{Z}^k$ by permutation of 
entries:
\[
\sigma ((\lambda_1 , \ldots , \lambda_k)) := \left( \lambda_{\sigma (1)}, \ldots , \lambda_{\sigma (k)}
 \right) \, .
\]
The dotted action of $\mathfrak{S}_k$ on $\mathbb{Z}^k$ is defined by 
\[
\sigma^{\bullet} (\lambda ):= \sigma (\lambda +\varrho )- \varrho \, .
\]
Then the theorem of Bott asserts in our case:
\begin{itemize}
\item
\textbf{Either} there exists $\sigma \in \mathfrak{S}_k$, $\sigma \neq \mathrm{id}$, such that 
$\sigma^{\bullet} (\lambda )= \lambda $. Then $R^i\pi_{\ast } \mathcal{L}(\lambda )= 0\; \forall i \in 
\mathbb{Z}$; 

\item
\textbf{or} there exists a unique $\sigma \in \mathfrak{S}_k$ such that $\sigma^{\bullet}( \lambda )= :
 \mu$ 
is non-increasing (i.e., $\mu$ is a dominant weight). Then 
\begin{gather*}
R^i\pi_{\ast} \mathcal{L}(\lambda ) =0 \quad \mathrm{for} \quad i\neq l (\sigma )\, ,\\
R^{l(\sigma )} \pi_{\ast } \mathcal{L}(\lambda )= \Sigma^{\mu } \mathcal{R}^{\vee } \, ,
\end{gather*}
where $l(\sigma )$ is the length of the permutation $\sigma$ (the smallest number of transpositions 
the composition of which gives $\sigma $) and $\Sigma^{\mu }$ is the Schur functor. 
\end{itemize}
As a first consequence, note that the objects $R\pi_{\ast }\mathcal{L}(\lambda )$ all belong -up to 
shift in cohomological degree- to the abelian subcategory of $D^b( Coh (\mathrm{IGrass}(k , V)))$ consisting of coherent 
sheaves. We would like to determine the homogeneous bundles that arise
as direct images of the bundles ($\sharp$) in this way. The following
theorem gives us some information (though it is not optimal).
\begin{theorem}
The derived category $D^b(Coh (\mathrm{IGrass}(k , V)))$ is generated by the bundles $\Sigma^{\nu } 
\mathcal{R}$, where $\nu$ runs over Young diagrams $Y$ which satisfy
\begin{gather*}
\left( \mathrm{number}\;\mathrm{of}\;\mathrm{columns}\;\mathrm{of}\;  Y \right)\le 2n-k\, ,\\
k \ge \left( \mathrm{number}\;\mathrm{of}\;\mathrm{rows}\;\mathrm{of}\;  Y \right) \ge 
\left( \mathrm{number}\;\mathrm{of}\;\mathrm{columns}\;\mathrm{of}\;  Y \right) - 2 (n-k)\, .
\end{gather*}
\end{theorem}
\begin{proof}
Note that if $\lambda $ satisfies the inequalities in ($\sharp$), then for $\delta := \lambda +\varrho$ 
we have 
\begin{eqnarray*}
-(2n-k) \le \delta_1 \le k-1 \, ,\\
(\sharp \sharp ) \hspace{2cm}-(2n-k-1) \le \delta_2 \le k-2 \, ,\\
             \vdots   \quad \quad \quad \quad  \hspace{2mm}            \\
-(2n-2k+1) \le \delta_k \le 0 \, . \quad \hspace{3mm}
\end{eqnarray*}
First of all one remarks that for $\sigma^{\bullet }(\lambda )= \sigma (\delta ) -\varrho$ to be 
non-increasing, it is necessary and sufficient that $\sigma (\delta )$ be strictly decreasing. We 
assume this to be the case in the following. Since the maximum possible value for $\sigma (\delta )_1$ 
is $k-1$, and the minimum possible value for $\sigma (\delta )_k$ is $-(2n-k)$, we find that for 
$\sigma^{\bullet }(\lambda )=: \mu$ 
\[
0 \ge \mu_1 \ge \ldots \ge \mu_k \ge -(2n-k) \, ;
\]
putting $\nu =(\nu_1 ,\ldots , \nu_k ):= (-\mu_k, -\mu_{k-1}, \ldots , -\mu_1)$ and 
noticing that $\Sigma^{\mu}\mathcal{R}^{\vee}\simeq \Sigma^{\nu}\mathcal{R}$, we find that the 
direct images $R^i\pi_{\ast }\mathcal{L}(\lambda )$, $i\in\mathbb{Z}$, $\mathcal{L}(\lambda )$ as in 
($\sharp$), will be a subset of the set of bundles $\Sigma^{\nu }\mathcal{R}$ on $\mathrm{IGrass} 
(k, V)$ where $\nu$ runs over the set of Young diagrams with no more than $2n-k$ columns and no 
more than $k$ rows.\\
But in fact we are only dealing with a proper subset of the latter: Suppose that 
\[
\sigma (\delta )_k = - (2n -k -a +1)\, , \quad 1\le a \le k-1 \, .
\]
Then the maximum possible value for $\sigma (\delta )_a$ is $k-a-1$. For in any case an upper bound 
for $\sigma (\delta )_a$ is $k-a$ because $\sigma (\delta )_1$ can be at most $k-1$ and the sequence 
$\sigma (\delta )$ is strictly decreasing. But in case this upper bound for $\sigma (\delta )_a$ is 
attained, the sequence $\sigma (\delta )$ must start with
\[
\sigma (\delta )_1 = k-1\, , \; \sigma (\delta )_2 =k-2 \, , \ldots , \; \sigma (\delta )_a=k-a \, ,
\]
in other words, we can only have 
\[
\sigma (\delta )_1 =\delta_1 \, ,  \ldots , \; \sigma (\delta )_a=\delta_a  \, .
\]
But this is impossible since $\delta_{a+1}, \ldots , \delta_k$ are all $\ge -(2n-k-a)> -(2n-k-a+1)$ 
and thus we could not have $\sigma (\delta )_k = -(2n-k-a+1)$. Hence $\sigma (\delta )_a$ is at most 
$k-a-1$, that is to say in $\sigma^{\bullet }(\lambda )= \sigma (\delta )-\varrho = \mu$ we have 
$\mu_a =\sigma (\delta )_a - (k-a) < 0$; or in terms of $\nu =( -\mu_k ,\ldots , -\mu_1)$ we 
can say that if the Young diagram $Y(\nu )$ of $\nu$ has $2n-k-a+1$ columns, $1\le a\le k-1$, it
 must have at least $k-a+1$ rows; or that the Young diagram $Y(\nu )$ satisfies 
\[
\left( \mathrm{number}\;\mathrm{of}\;\mathrm{rows}\;\mathrm{of}\;  Y(\nu ) \right) \ge 
\left( \mathrm{number}\;\mathrm{of}\;\mathrm{columns}\;\mathrm{of}\;  Y(\nu ) \right) - 2 (n-k)
\]
where the inequality is meaningless if the number on the right is $\le
0$. Thus by lemma 3.3.2 this concludes the proof of theorem 3.3.3.
\end{proof}

\begin{remark}
By thm. 2.2.2, in $D^b(Coh (\mathrm{Grass}(k, V)))$ there is a complete exceptional 
sequence consisting of the $\Sigma^{\tilde{\nu }} \tilde{\mathcal{R}}$ where 
$\tilde{\mathcal{R}}$ is the tautological subbundle on $\mathrm{Grass}(k , V)$ and $\tilde{\nu }$
 runs over Young 
diagrams with at most $2n-k$ columns and at most $k$ rows. Looking at $\mathrm{IGrass}(k, V)$ as a 
subvariety $\mathrm{IGrass}(k, V)\subset\mathrm{Grass}(k , V)$ we see that the bundles in theorem 
3.3.3 form a proper subset of the restrictions of the $\Sigma^{\tilde{\nu }} \tilde{\mathcal{R}}$
 to $\mathrm{IGrass}(k, V)$.  
\end{remark}

Before making the next remark we have to recall two ingredients in
order to render the following computations transparent:\\
The first is the Littlewood-Richardson rule to decompose $\Sigma^{\lambda }\otimes \Sigma^{\mu}$ into 
irreducible factors where $\lambda$, $\mu$ are Young diagrams (cf. \cite{FuHa}, \S A.1). It says the 
following: Label each box of $\mu$ with the number of the row it belongs to. Then expand the Young 
diagram $\lambda$ by adding the boxes of $\mu$ to the rows of
$\lambda$ subject to the following rules:
\begin{itemize}
\item[(a)]
The boxes with labels $\le i$ of $\mu$ together with the boxes of $\lambda$ form again a Young diagram;

\item[(b)]
No column contains boxes of $\mu$ with equal labels.

\item[(c)]
When the integers in the boxes added are listed from right to left and from top down, then, for any 
$0\le s \le$ (number of boxes of $\mu$), the first $s$ entries of the list satisfy: Each label $l$ 
($1\le l\le$ (number of rows of $\mu$)$-1$ ) occurs at least as many times as the label $l+1$.

\end{itemize}
Then the multiplicity of $\Sigma^{\nu}$ in $\Sigma^{\lambda}\otimes\Sigma^{\mu}$ is the number 
of times the Young diagram $\nu$ can be obtained by expanding $\lambda$ by $\mu$ according to 
the above rules, forgetting the labels.\\ 
The second point is the calculation of the cohomology of the bundles $\Sigma^{\lambda} \mathcal{R}$ on 
the variety $\mathrm{IGrass}(k , V)$, $V$ an $n$-dimensional symplectic vector space (cf. \cite{Wey},
 cor. 4.3.4). Bott's 
theorem gives the following prescription:\\
Look at the sequence 
\[
\mu= (-\lambda_k, -\lambda_{k-1}, \ldots , -\lambda_1, 0 , \ldots , 0) \in \mathbb{Z}^n
\]
considered as a weight of the root system of type $C_n$. Let $W$ be the Weyl group of this root system 
which is a semi-direct product of $(\mathbb{Z}/2\mathbb{Z})^n$ with the symmetric group $\mathfrak{S}_n$ 
and acts on weights by permutation and sign changes of entries. Let $\varrho := (n , n-1, \ldots , 1)$ be
 the half sum of the positive roots for type $C_n$. The dotted action of $W$ on weights is defined as 
above by $\sigma^{\bullet}(\mu ):= \sigma (\mu +\varrho )-\varrho $. Then
\begin{itemize}
\item
either there is $\sigma\in W$, $\sigma \neq \mathrm{id}$, such that $\sigma^{\bullet}(\mu )= \mu$. 
then all cohomology groups
\[
H^{\bullet }( \mathrm{IGrass}(k, V) , \Sigma^{\lambda }\mathcal{R}) = 0 \, .
\]
\item
or there is a unique $\sigma \in W$ such that $\sigma^{\bullet }(\mu )=: \nu$ is dominant (a 
non-increasing sequence of non-negative integers). Then the only non-zero cohomology group is 
\[
H^{l(\sigma ) }( \mathrm{IGrass}(k, V) , \Sigma^{\lambda }\mathcal{R}) = V_{\nu}\, ,
\]
where $l(\sigma )$ is the length of the Weyl group element $\sigma$ and $V_{\nu}$ is the space of the 
irreducible representation of $\mathrm{Sp}_{2n} \mathbb{C}$ with highest weight $\nu$.
\end{itemize} 

\begin{remark}
The $R^i\pi_{\ast}\mathcal{L}(\lambda )$, $i\in \mathbb{Z}$, $\mathcal{L}(\lambda )$ as in ($\sharp$), 
are not in general exceptional: For example, take $k=n=3$, so that we are dealing with $\mathrm{IGrass}
(3, V)$, the Lagrangian Grassmannian of maximal isotropic subspaces in a $6$-dimensional symplectic 
space $V$. Then $\mathcal{L}((0, -3 , 0))$ is in ($\sharp$). Adding $\varrho =(2,1,0)$ to $(0,-3,0)$ 
we get $(2,-2,0)$ and interchanging the last two entries and subtracting $\varrho$ again, we arrive at 
$(0, -1, -2)$ which is non-increasing whence 
\[
R^1\pi_{\ast }\mathcal{L}((0,-3,0) )=\Sigma^{2,1,0}\mathcal{R}\, ,
\]
all other direct images being $0$. To calculate
\begin{eqnarray*}
\mathrm{Ext}^{\bullet} \left( \Sigma^{2,1,0}\mathcal{R}, \Sigma^{2,1,0}\mathcal{R} \right) 
= H^{\bullet}\left( \mathrm{IGrass}(3, V), \Sigma^{2,1,0}\mathcal{R} \otimes \Sigma^{0,-1,-2}
\mathcal{R} \right) \\
= H^{\bullet}\left( \mathrm{IGrass}(3, V), \Sigma^{2,1,0}\mathcal{R} \otimes \Sigma^{2,1,0}
\mathcal{R}\otimes \left( \bigwedge\nolimits^3 \mathcal{R}^{\vee } \right)^{\otimes 2} \right)
\end{eqnarray*}
we use the Littlewood-Richardson rule and Bott's theorem as recalled
above: One gets that 
\begin{gather*}
\Sigma^{2,1,0}\mathcal{R} \otimes \Sigma^{0 , -1, -2}\mathcal{R} = \Sigma^{2,0,-2}\mathcal{R}\oplus 
\Sigma^{2,-1,-1}\mathcal{R}\oplus \Sigma^{1,1,-2}\mathcal{R}\\
\oplus (\Sigma^{1,0,-1}\mathcal{R})^{\oplus 2
}\oplus  \Sigma^{0,0,0}\mathcal{R}
\end{gather*}
in view of the fact that if we expand $\lambda =(2,1,0)$ by $\mu =(2,1,0)$ we get the following Young 
diagrams according to the Littlewood-Richardson rule:
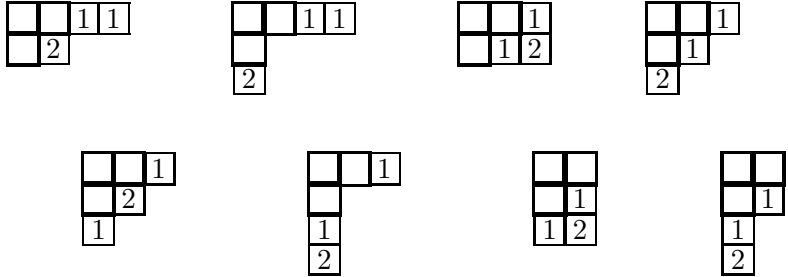
\setlength{\unitlength}{1cm}
\begin{center}
\begin{picture}(4,4.5)
\thicklines
\put(-4,4){\framebox(0.4,0.4){}}
\put(-3.6,4){\framebox(0.4,0.4){}}
\thinlines
\put(-3.2,4){\framebox(0.4,0.4){1}}
\put(-2.8,4){\framebox(0.4,0.4){1}}
\thicklines
\put(-4,3.6){\framebox(0.4,0.4){}}
\thinlines
\put(-3.6,3.6){\framebox(0.4,0.4){2}}

\thicklines
\put(-1,4){\framebox(0.4,0.4){}}
\put(-0.6,4){\framebox(0.4,0.4){}}
\thinlines
\put(-0.2,4){\framebox(0.4,0.4){1}}
\put(0.2,4){\framebox(0.4,0.4){1}}
\thicklines
\put(-1,3.6){\framebox(0.4,0.4){}}
\thinlines
\put(-1,3.2){\framebox(0.4,0.4){2}}

\thicklines
\put(2,4){\framebox(0.4,0.4){}}
\put(2.4,4){\framebox(0.4,0.4){}}
\thinlines
\put(2.8,4){\framebox(0.4,0.4){1}}
\thicklines
\put(2,3.6){\framebox(0.4,0.4){}}
\thinlines
\put(2.4,3.6){\framebox(0.4,0.4){1}}
\put(2.8,3.6){\framebox(0.4,0.4){2}}

\thicklines
\put(4.5,4){\framebox(0.4,0.4){}}
\put(4.9,4){\framebox(0.4,0.4){}}
\thinlines
\put(5.3,4){\framebox(0.4,0.4){1}}
\thicklines
\put(4.5,3.6){\framebox(0.4,0.4){}}
\thinlines
\put(4.9,3.6){\framebox(0.4,0.4){1}}
\put(4.5,3.2){\framebox(0.4,0.4){2}}

\thicklines
\put(-3,2){\framebox(0.4,0.4){}}
\put(-2.6,2){\framebox(0.4,0.4){}}
\thinlines
\put(-2.2,2){\framebox(0.4,0.4){1}}
\put(-3,1.2){\framebox(0.4,0.4){1}}
\thicklines
\put(-3,1.6){\framebox(0.4,0.4){}}
\thinlines
\put(-2.6,1.6){\framebox(0.4,0.4){2}}

\thicklines
\put(0,2){\framebox(0.4,0.4){}}
\put(0.4,2){\framebox(0.4,0.4){}}
\thinlines
\put(0.8,2){\framebox(0.4,0.4){1}}
\put(0,0.8){\framebox(0.4,0.4){2}}
\thicklines
\put(0,1.6){\framebox(0.4,0.4){}}
\thinlines
\put(0,1.2){\framebox(0.4,0.4){1}}

\thicklines
\put(3,2){\framebox(0.4,0.4){}}
\put(3.4,2){\framebox(0.4,0.4){}}
\thinlines
\put(3,1.2){\framebox(0.4,0.4){1}}
\thicklines
\put(3,1.6){\framebox(0.4,0.4){}}
\thinlines
\put(3.4,1.6){\framebox(0.4,0.4){1}}
\put(3.4,1.2){\framebox(0.4,0.4){2}}

\thicklines
\put(5.5,2){\framebox(0.4,0.4){}}
\put(5.9,2){\framebox(0.4,0.4){}}
\thinlines
\put(5.5,0.8){\framebox(0.4,0.4){2}}
\thicklines
\put(5.5,1.6){\framebox(0.4,0.4){}}
\thinlines
\put(5.9,1.6){\framebox(0.4,0.4){1}}
\put(5.5,1.2){\framebox(0.4,0.4){1}}

\end{picture}\end{center}
\vspace{0.3cm}
Thus calculating the cohomology of $\Sigma^{2,1,0}\mathcal{R} \otimes \Sigma^{0,-1,-2}\mathcal{R}$ by 
the version of Bott's theorem recalled above one finds that
\[
\mathrm{Hom}\left( \Sigma^{2,1,0}\mathcal{R}, \Sigma^{2,1,0}\mathcal{R}\right) = \mathbb{C}\, \quad 
\mathrm{Ext}^1\left( \Sigma^{2,1,0}\mathcal{R}, \Sigma^{2,1,0}\mathcal{R}\right) =V_{1,1,0} \oplus 
V_{2,0,0} \neq 0
\]
the other Ext groups being $0$.
\end{remark}
Next we want to show by some examples that, despite the fact that theorem 3.3.3 does not give a complete  
exceptional sequence on $\mathrm{IGrass}(k , V)$, it is sometimes -for small values of $k$ and $n$- 
not so hard to find one with its help. 
\begin{example}
Choose $k=n=2$, i.e. look at $\mathrm{IGrass}(2 , V)$, $\dim V=4$. Remarking that $\mathcal{O}(1)$ on 
$\mathrm{IGrass}(2 , V)$ in the Pl\"ucker embedding equals $\bigwedge^{\mathrm{top}} \mathcal{R}^{\vee}$ 
and applying theorem 3.3.3 one finds that the following five sheaves generate $D^b(Coh (\mathrm{IGrass}
(2 , V)))$:
\[
\mathcal{O}\, , \; \mathcal{R}\, ,\; \bigwedge\nolimits^2 \mathcal{R}=\mathcal{O}(-1)\, , \; 
\Sigma^{2,1}\mathcal{R}= \mathcal{R}(-1)\, ,
\; \mathcal{O}(-2)\, ;
\] 
The real extra credit that one receives from working on the Lagrangian Grassmannian $\mathrm{IGrass}
(2 , V)$ is that $\mathcal{R}=\mathcal{R}^{\perp}$ and the tautological factor bundle can be identified 
with $\mathcal{R}^{\vee }\simeq \mathcal{R}(1)$, i.e. one has an exact sequence
\[
\begin{CD}
0 @>>> \mathcal{R} @>>> V\otimes \mathcal{O} @>>> \mathcal{R}(1) @>>> 0 \, .
\end{CD}
\] 
Twisting by $\mathcal{O}(-1)$ in this sequence shows that of the above five sheaves, $\mathcal{R}(-1)$ 
is in the full triangulated subcategory generated by the remaining
four; moreover, it is a straightforward computation with Bott's theorem that 
\[
\left(  \mathcal{O}(-2) , \: \mathcal{O}(-1) , \:\mathcal{R} , \:  \mathcal{O} \right)
\]
is a strong exceptional sequence in $D^b(Coh (\mathrm{IGrass}(2 , V)))$; but this is also complete, 
i.e., it generates this derived category by the preceding considerations. In fact, this does not come as a 
surprise. $\mathrm{IGrass}(2 , V)$ is isomorphic to a quadric hypersurface in $\mathbb{P}^4$, more 
precisely it is a hyperplane section of the Pl\"ucker quadric $\mathrm{Grass}(2 , V)\subset \mathbb{P}^5$. 
By \cite{Ott}, thm. 1.4 and ex. 1.5, the spinor bundles on the Pl\"ucker quadric are the 
dual of the tautological 
subbundle and the tautological factor bundle on $\mathrm{Grass}(2 , V)$ and these both restrict to 
the spinor bundle $\mathcal{R}^{\vee}$ on $\mathrm{IGrass}(2 , V)\subset \mathbb{P}^4$ (let us renew here 
the warning from subsection 2.2 that the spinor bundles in \cite{Ott} are the duals of the bundles that we 
choose to call spinor bundles in this work). 
 We thus recover the result 
of \cite{Ka3}, \S 4, in a special case. Note that the identification of $\mathrm{IGrass}(2 , V)$ with a 
quadric hypersurface in $\mathbb{P}^4$ also follows more conceptually from the isomorphism of marked 
Dynkin diagrams
\setlength{\unitlength}{1cm}
\begin{center}
\begin{picture}(4,1.5)

\put(-3,0.5){\circle{0.4}}
\put(-1,0.5){\circle*{0.4}}
\put(-2.8,0.57){\line(1,0){1.6}}
\put(-2.8,0.43){\line(1,0){1.6}}
\put(-2.2, 0.5){\line(2,1){0.36}}
\put(-2.2, 0.5){\line(2,-1){0.36}}

\put(-3.2,1){$\alpha_1$}
\put(-1.2,1){$\alpha_2$}

\put(-2.2,-0.5){$C_2$}

\put(1.5,0.5){$\simeq$}

\put(4,0.5){\circle*{0.4}}
\put(6,0.5){\circle{0.4}}
\put(4.2,0.57){\line(1,0){1.6}}
\put(4.2,0.43){\line(1,0){1.6}}
\put(5.2, 0.5){\line(-2,1){0.36}}
\put(5.2, 0.5){\line(-2,-1){0.36}}

\put(3.8,1){$\alpha_1'$}
\put(5.8,1){$\alpha_2'$}

\put(4.8,-0.5){$B_2$}

\end{picture}\end{center}
\vspace{0.5cm}
corresponding to the isomorphism $\mathrm{Sp}_4\mathbb{C}/ P(\alpha_2)\simeq \mathrm{Spin}_5 \mathbb{C}/ 
P(\alpha_1')$ (cf. \cite{Stei}, prop. p. 16 and \cite{FuHa}, \S 23.3). Recalling the one-to-one 
correspondence between conjugacy classes of parabolic subgroups of a simple complex Lie group $G$ and subsets of the set of 
simple roots, the notations $P(\alpha_2)$ resp. $P(\alpha_1')$ are self-explanatory.
\end{example}

\begin{example}
Along the same lines which are here exposed in general, A. V. Samokhin treated in \cite{Sa} the particular 
case of $\mathrm{IGrass}(3, V)$, $\dim V=6$, and using the identification of the tautological factor 
bundle with $\mathcal{R}^{\vee}$ on this Lagrangian Grassmannian and the exact sequence
\[
\begin{CD}
0 @>>> \mathcal{R} @>>> V\otimes \mathcal{O} @>>> \mathcal{R}^{\vee} @>>> 0 \, .
\end{CD}
\] 
together with its symmetric and exterior powers found the following strong complete exceptional 
sequence for $D^b(Coh (\mathrm{IGrass}(3 ,V)))$: 
\[
\left( \mathcal{R}(-3),\: \mathcal{O}(-3),\: \mathcal{R}(-2),\: \mathcal{O}(-2),\: 
\mathcal{R}(-1),\: \mathcal{O}(-1),\: \mathcal{R},\: \mathcal{O}      \right)
\]
and we refer to \cite{Sa} for details of the computation.
\end{example}

In general I conjecture that on any Lagrangian Grassmannian
$\mathrm{IGrass}(n , V)$, $\dim V=2n$, every ``relation'' between the
bundles in theorem 3.3.3 in the derived category $D^b(Coh \,
\mathrm{IGrass}(n , V))$ (that is to say that one of these bundles is
in the full triangulated subcategory generated by the remaining ones)
should follow using the Schur complexes (cf. \cite{Wey}, section 2.4)
derived from the exact sequence $0\to \mathcal{R}\to V\otimes
\mathcal{O}\to \mathcal{R}^{\vee}\to 0$ (and the Littlewood-Richardson
rule). \\
Let us conclude this subsection by giving an example which, though we
do not find a complete exceptional sequence in the end, may help to
convey the sort of combinatorial difficulties that one encounters in
general.

\begin{example}
For a case of a non-Lagrangian isotropic Grassmannian, look at
$\mathrm{IGrass}(2 , V)$, $\dim V=6$. Theorem 3.3.3 says that $D^b(
Coh\, \mathrm{IGrass}( 2 , V))$ is generated by the following $14$
bundles:
\[
(\ast )\quad\quad \mathrm{Sym}^a \, \mathcal{R} (-b) \, ,\quad 0\le a
\le 3\, , \quad 0\le
b\le 4-a \, .
\]
By corollary 2.3.7, the number of terms in a complete exceptional
sequence must be $12$ in this case (in general for
$\mathrm{IGrass}(k , V)=\mathrm{Sp}_{2n}\,\mathbb{C}/P(\alpha_k)$,
$\dim V=2n$, one has that $W^{P(\alpha_k)}$, the set of minimal
representatives of the quotient $W/W_{P(\alpha_k)}$ can be identified
with $k$-tuples of integers $(a_1,\ldots , a_k)$ such that 
\begin{gather*}
1\le a_1 < a_2 < \ldots < a_k \le 2n\; \mathrm{and} \\
\mathrm{for}\; 1\le i\le 2n\, , \; \mathrm{if}\; i\in \{ a_1,
\ldots , a_k\} \quad \mathrm{then}\; 2n+1-i \notin \{ a_1 ,\ldots , a_k\} 
\end{gather*}
(see \cite{BiLa}, \S 3.3) and these are 
\[
\frac{2n (2n-2)\dots (2n-2(k-1))}{1\cdot 2\cdot \dots \cdot k}=
2^k{n\choose k}
\]
in number). Without computation, we know by a theorem of Ramanan
(cf. \cite{Ot2}, thm 12.3) that the bundles in ($\ast$) are all simple
since they are associated to irreducible representations of the
subgroup $P(\alpha_2)\subset \mathrm{Sp}_{6}\,\mathbb{C}$.\\
Moreover the bundles 
\[
\Sigma^{c_1 , c_2}\mathcal{R} \;\; \mathrm{and}\;\; \Sigma^{d_1 ,
  d_2}\mathcal{R}\;\; \mathrm{with}\; 0\le c_2 \le c_1 \le 3 \, ,\;
0\le d_2 \le d_1\le 3 
\]
have no higher extension groups between each other: By the
Littlewood-Richardson rule, every irreducible summand $\Sigma^{e_1
  ,e_2 }\mathcal{R}$ occurring in the decomposition of $\Sigma^{d_1 ,
  d_2}\mathcal{R}\otimes \Sigma^{c_1 , c_2}\mathcal{R}^{\vee}$
satisfies $-3\le e_2 \le e_1 \le 3$ and hence for $\mu := (-e_2 ,
-e_1, 0)\in\mathbb{Z}^3$ and $\varrho =(3,2,1)$ we find that $\mu
+\varrho$ is either a strictly decreasing sequence of positive
integers or two entries in $\mu+\varrho$ are equal up to sign or one
entry in $\mu+\varrho$ is $0$. In each of these cases, Bott's
theorem as recalled before remark 3.3.5 tells us that
$H^i(\mathrm{IGrass}(2, V), \Sigma^{e_1 , e_2} \mathcal{R})=0\;
\forall i> 0$. Combining this remark with the trivial observation that
for $\mathrm{Sym}^a\mathcal{R}(-b)$, $\mathrm{Sym}^c\mathcal{R}(-d)$
in the set ($\ast$) with $b, d\ge 1$ we have 
\[
\mathrm{Ext}^i(\mathrm{Sym}^a\mathcal{R}(-b),
\mathrm{Sym}^c\mathcal{R}(-d))=\mathrm{Ext}^i(\mathrm{Sym}^a\mathcal{R}(-b+1),
\mathrm{Sym}^c\mathcal{R}(-d+1)) \;\;\forall i  
\]
we infer that for $\mathcal{A}$, $\mathcal{B}$ bundles in the set
($\ast$) we can only have 
\[
\mathrm{Ext}^j (\mathcal{A} , \mathcal{B})\neq 0 , \;\;
\mathrm{some}\;\; j> 0
\]
if $\mathcal{A}$ occurs in the set 
\[
S_1:=\{\mathcal{O}(-4) ,\; \mathcal{R}(-3) ,\;
\mathrm{Sym}^2\mathcal{R}(-2), \; \mathrm{Sym}^3 \mathcal{R}(-1) \}
\]
and $\mathcal{B}$ is in the set 
\[
S_2:=\{\mathcal{O}, \; \mathcal{R} ,\; \mathrm{Sym}^2\mathcal{R} ,\;
\mathrm{Sym}^3 \mathcal{R}  \}
\]
or vice versa. By explicit calculation (which amounts to applying
Bott's theorem another 32 more times) we find that the only non-vanishing
higher extension groups between two bundles in ($\ast$) are the following: 
\begin{gather*}
\mathrm{Ext}^1( \mathrm{Sym}^3 \mathcal{R} , \mathcal{R}(-3)
)=\mathbb{C} , \;\; \mathrm{Ext}^1( \mathrm{Sym}^2 \mathcal{R} , \mathrm{Sym}^2\mathcal{R}(-2)
)=\mathbb{C} \\
\mathrm{Ext}^1( \mathrm{Sym}^3 \mathcal{R} , \mathrm{Sym}^2\mathcal{R}(-2)
)=V , \;\; \mathrm{Ext}^1( \mathcal{R} , \mathrm{Sym}^3\mathcal{R}(-1)
)=\mathbb{C}\\
\mathrm{Ext}^1( \mathrm{Sym}^2 \mathcal{R} , \mathrm{Sym}^3\mathcal{R}(-1)
)=V , \;\; \mathrm{Ext}^1( \mathrm{Sym}^3 \mathcal{R} , \mathrm{Sym}^3\mathcal{R}(-1)
)=V_{2,0,0} .
\end{gather*}
Thus in this case the set of bundles ($\ast$) does not contain a
strong complete exceptional sequence. It does not contain a complete
exceptional sequence, either, since 
\begin{gather*}
\mathrm{Hom}(\mathcal{R}(-3), \mathrm{Sym}^3 \mathcal{R}
)\neq 0 , \;\; \mathrm{Hom}(\mathrm{Sym}^2\mathcal{R}(-2), \mathrm{Sym}^2 \mathcal{R}
)\neq 0 \\
\mathrm{Hom}( \mathrm{Sym}^2\mathcal{R}(-2), \mathrm{Sym}^3 \mathcal{R}
)\neq 0 , \;\; \mathrm{Hom}(\mathrm{Sym}^3\mathcal{R}(-1), \mathcal{R}
)\neq 0\\
\mathrm{Hom}(\mathrm{Sym}^3\mathcal{R}(-1), \mathrm{Sym}^2 \mathcal{R}
)\neq 0 , \;\; \mathrm{Hom}( \mathrm{Sym}^3\mathcal{R}(-1), \mathrm{Sym}^3 \mathcal{R}
)\neq 0 .
\end{gather*}
On the other hand one has on $\mathrm{IGrass}(2 , V)$ the following
exact sequences of vector bundles:
\begin{eqnarray}
0\to \mathcal{R}^{\perp} \to V\otimes \mathcal{O} \to
\mathcal{R}^{\vee}\to 0 \\
0\to \mathcal{R} \to \mathcal{R}^{\perp} \to
\mathcal{R}^{\perp}/\mathcal{R}\to 0 .
\end{eqnarray}
The second exterior power of the two term complex $0\to \mathcal{R}\to
\mathcal{R}^{\perp}$ gives an acyclic complex resolving $\bigwedge^2
(\mathcal{R}^{\perp}/\mathcal{R})$ which is isomorphic to
$\mathcal{O}_{\mathrm{IGrass}( 2, V)}$ via the mapping induced by the
symplectic form $\langle \cdot , \cdot \rangle $. Thus we get the
exact sequence 
\begin{equation}
0\to \mathrm{Sym}^2 \mathcal{R} \to \mathcal{R}\otimes
\mathcal{R}^{\perp} \to \bigwedge\nolimits^2 \mathcal{R}^{\perp} \to
\mathcal{O}\to 0 .
\end{equation}
The second symmetric power of the two term complex $0\to
\mathcal{R}^{\perp}\to V\otimes \mathcal{O}$ yields the exact sequence
\begin{equation}
0\to \bigwedge\nolimits^2 \mathcal{R}^{\perp} \to
\mathcal{R}^{\perp}\otimes V \to \mathrm{Sym}^2\, V\otimes
\mathcal{O} \to \mathrm{Sym}^2\, \mathcal{R}^{\vee} \to 0 .  
\end{equation}
Note also that $\mathcal{R}^{\vee} \simeq \mathcal{R}(1)$ and
$\mathrm{Sym}^2 \, \mathcal{R}^{\vee}\simeq \mathrm{Sym}^2 \,
\mathcal{R}(2)$. Since $\mathcal{R}\otimes \mathcal{R}(-1) \simeq
\mathrm{Sym}^2 \, \mathcal{R} (-1)\oplus \mathcal{O}(-2)$ sequence (1)
gives 
\begin{equation}
 0\to \mathcal{R}^{\perp}\otimes \mathcal{R}(-2)  \to V\otimes \mathcal{R}(-2) \to
\mathrm{Sym}^2 \, \mathcal{R} (-1)\oplus \mathcal{O}(-2)\to 0
\end{equation}
and
\begin{equation}
0\to \mathcal{R}^{\perp}(-2) \to V\otimes \mathcal{O}(-2) \to
\mathcal{R}(-1) \to 0 .
\end{equation}
Moreover twisting by $\mathcal{O}(-2)$ in (3) and (4) yields 
\begin{eqnarray}
0\to \mathrm{Sym}^2 \mathcal{R}(-2) \to \mathcal{R}\otimes
\mathcal{R}^{\perp}(-2) \to \bigwedge\nolimits^2 \mathcal{R}^{\perp}(-2) \to
\mathcal{O}(-2)\to 0 \\
0\to \bigwedge\nolimits^2 \mathcal{R}^{\perp}(-2) \to
\mathcal{R}^{\perp}\otimes V(-2) \to \mathrm{Sym}^2\, V\otimes
\mathcal{O}(-2) \to \mathrm{Sym}^2\, \mathcal{R} \to 0 .
\end{eqnarray}
What (5), (6), (7), (8) tell us is that $\mathrm{Sym}^2
\,\mathcal{R}(-2)$ is in the full triangulated subcategory generated
by $\mathcal{O}(-2) ,\; \mathrm{Sym}^2\,\mathcal{R}(-1) ,\;
\mathcal{R}(-2) ,\; \mathrm{Sym}^2\,\mathcal{R} ,\;
\mathcal{R}(-1)$. Thus the derived category $D^b (Coh\,
\mathrm{IGrass} (2 , V))$ is generated by the bundles in ($\ast$)
without $\mathrm{Sym}^2\,\mathcal{R}(-2)$, which makes a total of 13
bundles.\\
But even in this simple case I do not know how to pass on to a
complete exceptional sequence because there is no method at this
point to decide which bundles in ($\ast$) should be thrown away and
what extra bundles should be let in to obtain a complete exceptional
sequence. 
\end{example}

\subsection{Calculation for the Grassmannian of isotropic 3-planes in
  a 7-dimensional orthogonal vector space }

In this section we want to show how the method of subsection 3.3 can be
adapted -using theorem 3.2.7 on quadric bundles- to produce sets of
vector bundles that generate the derived categories of coherent
sheaves on orthogonal Grassmannians (with the ultimate goal to obtain
(strong) complete exceptional sequences on them by appropriately
modifying these sets of bundles). Since the computations are more
involved than in the symplectic case, we will restrict ourselves to
illustrating the method by means of a specific example:\\
Let $V$ be a 7-dimensional complex vector space equipped with a
non-degenerate symmetric bilinear form $\langle \cdot , \cdot \rangle
$. $\mathrm{IFlag}(k_1, \dots , k_t ; V)$ denotes the flag variety of
isotropic flags of type $(k_1 ,\dots , k_t)$, $1\le k_1 < \dots < k_t
\le 3$, in $V$ and $\mathrm{IGrass} (k , V)$, $1\le k \le 3$, the
Grassmannian of isotropic $k$-planes in $V$; again in this setting we
have the tautological flag of subbundles  
\[
\mathcal{R}_{k_1} \subset \dots \subset \mathcal{R}_{k_t} \subset
V\otimes \mathcal{O}_{\mathrm{IFlag}(k_1, \dots , k_t ; V)}
\]
on $\mathrm{IFlag}(k_1, \dots , k_t ; V)$ and the tautological
subbundle $\mathcal{R}$ on $\mathrm{IGrass}(k , V)$.\\
Now consider $\mathrm{IGrass}(3 , V)$ which sits in the diagram (D)
\begin{center}
\setlength{\unitlength}{1cm}
\begin{picture}(2 , 4)
\put(-4,0){$\mathbb{P}^6\supset Q \simeq \mathrm{IFlag}(1 ; V)$}
\put(-4, 1.5){$\mathbb{P}(\mathcal{E}_1)\supset \mathcal{Q}_1\simeq
  \mathrm{IFlag}(1, 2 ; V)$}
\put(-4 ,3){$\mathbb{P}(\mathcal{E}_2)\supset \mathcal{Q}_2\simeq
  \mathrm{IFlag}(1, 2 , 3 ; V)\simeq \mathrm{Flag}_{\mathrm{IGrass}(3
    , V)}(1, 2, 3 ; \mathcal{R})$}
\put(-1 , 1.3){\vector(-1,-3){0.3}} \put(-0.8 , 0.8){$\pi_1$}
\put(-0.7 , 2.7){\vector(-1,-3){0.3}}  \put(-0.6 , 2.2){$\pi_2$}
\put(3 , 0){$\mathrm{IGrass}(3 , V)$}
\put(3 , 2.7){\vector(1 , -3){0.8}}
\put(3.8 , 1.5){$\pi$}
\put(6 , 1.5){(D)}
\end{picture}
\end{center}
\vspace{0.5cm}
The rank $i$ tautological subbundle on $\mathrm{IFlag}(1, \dots , j ;
V)$ pulls back to the rank $i$ tautological subbundle on $\mathrm{IFlag}(1, \dots , j+1 ;
V)$ under $\pi_j$, $1\le i\le j\le 2$, and for ease of notation it will
be denoted by $\mathcal{R}_i$ with the respective base spaces being tacitly
understood in each case.\\
The choice of an isotropic line $L_1$ in $V$ amounts to picking a
point in the quadric hypersurface $Q=\{ [v]\in \mathbb{P}(V)\, | \,
\langle v , v\rangle =0 \} \subset \mathbb{P}^6$. An isotropic plane
$L_2$ containing $L_1$ is nothing but an isotropic line $L_2/L_1$ in
the orthogonal vector space $L_1^{\perp}/L_1$. Thus $\mathrm{IFlag}(1
, 2 ; V)$ is a quadric bundle $\mathcal{Q}_1$ over $\mathrm{IFlag}(1 ;
V)$ inside the projective bundle $\mathbb{P}(\mathcal{E}_1)$ of the
orthogonal vector bundle
$\mathcal{E}_1:=\mathcal{R}_1^{\perp}/\mathcal{R}_1$ on
$\mathrm{IFlag}(1 ; V)$. Similarly, $\mathrm{IFlag}(1 , 2 ,3 ; V)$ is
a quadric bundle $\mathcal{Q}_2\subset \mathbb{P}(\mathcal{E}_2)$ over
$\mathrm{IFlag}(1, 2 ; V)$ where $\mathcal{E}_2:=
\mathcal{R}_2^{\perp}/\mathcal{R}_2$, and at the same time
$\mathrm{IFlag}(1 , 2, 3; V)$ is isomorphic to the relative variety of
complete flags $\mathrm{Flag}_{\mathrm{IGrass}(k , V)}(1, 2, 3 ;
\mathcal{R})$ in the fibres of the tautological subbundle
$\mathcal{R}$ on $\mathrm{IGrass}(3 , V)$.\\
Moreover, $\mathcal{O}_{Q}(-1)\simeq \mathcal{R}_1$,
$\mathcal{O}_{\mathcal{Q}_1}(-1)\simeq \mathcal{R}_2/\mathcal{R}_1$,
$\mathcal{O}_{\mathcal{Q}_2}(-1) \simeq
\mathcal{R}_3/\mathcal{R}_2$. By means of the constructions of section
3.2 we have on $Q$ the spinor bundle $\Sigma
(\mathcal{O}_Q(-1)^{\perp}/ \mathcal{O}_Q(-1))$ for the orthogonal
vector bundle $\mathcal{O}_Q(-1)^{\perp}/ \mathcal{O}_Q(-1)$, and on
the quadric bundles $\mathcal{Q}_1$ resp. $\mathcal{Q}_2$ the spinor
bundles\\ $\Sigma
(\mathcal{O}_{\mathcal{Q}_1}(-1)^{\perp}/
\mathcal{O}_{\mathcal{Q}_1}(-1))$ resp. $\Sigma
(\mathcal{O}_{\mathcal{Q}_2}(-1)^{\perp}/
\mathcal{O}_{\mathcal{Q}_2}(-1))$ for the orthogonal vector bundles $\mathcal{O}_{\mathcal{Q}_1}(-1)^{\perp}/
\mathcal{O}_{\mathcal{Q}_1}(-1)$ resp. $\mathcal{O}_{\mathcal{Q}_2}(-1)^{\perp}/
\mathcal{O}_{\mathcal{Q}_2}(-1)$.\\
Note that under the identifications $Q\simeq \mathrm{IFlag}(1 ; V)$,
$\mathcal{Q}_1\simeq \mathrm{IFlag}(1 , 2 ; V)$
resp. $\mathcal{Q}_2\simeq \mathrm{IFlag}(1, 2, 3 ;V)$ we get the
isomorphisms of orthogonal vector bundles
\begin{gather*}
\mathcal{R}_1^{\perp}/\mathcal{R}_1 \simeq \mathcal{O}_Q(-1)^{\perp}/
\mathcal{O}_Q(-1)\, , \ \mathcal{R}_2^{\perp}/\mathcal{R}_2 \simeq \mathcal{O}_{\mathcal{Q}_1}(-1)^{\perp}/
\mathcal{O}_{\mathcal{Q}_1}(-1)\\
\mathrm{resp.}\ \mathcal{R}_3^{\perp}/\mathcal{R}_3 \simeq \mathcal{O}_{\mathcal{Q}_2}(-1)^{\perp}/
\mathcal{O}_{\mathcal{Q}_2}(-1) . 
\end{gather*}
Therefore, by theorem 3.2.7, we get that the following set of bundles
on $\mathrm{IFlag}(1, 2, 3; V)$ generates $D^b (Coh \,
\mathrm{IFlag}(1, 2, 3 ; V))$, and in fact forms a complete
exceptional sequence when appropriately ordered: The bundles
\[
(\heartsuit )\ \quad \mathcal{A}\otimes \mathcal{B}\otimes \mathcal{C}
\]
where $\mathcal{A}$ runs through the set
\[
A:=\{ \Sigma (\mathcal{R}_1^{\perp}/ \mathcal{R}_1)\otimes
\mathcal{R}_1^{\otimes 5} ,\; \mathcal{R}_1^{\otimes 4} , \;
\mathcal{R}_1^{\otimes 3} , \; \mathcal{R}_1^{\otimes 2} , \;
\mathcal{R}_1 , \; \mathcal{O}    \}
\]
and $\mathcal{B}$ runs through
\[
B:= \{ \Sigma (\mathcal{R}_2^{\perp}/ \mathcal{R}_2)\otimes
(\mathcal{R}_2/\mathcal{R}_1)^{\otimes 3} , \; 
(\mathcal{R}_2/\mathcal{R}_1)^{\otimes 2} , \; 
\mathcal{R}_2/\mathcal{R}_1 , \; \mathcal{O}  \}
\]
and $\mathcal{C}$ runs through
\[
C:= \{  \Sigma (\mathcal{R}_3^{\perp} / \mathcal{R}_3)\otimes
(\mathcal{R}_3/\mathcal{R}_2) , \; \mathcal{O}  \}  .
\]
We know that the full direct images under $\pi$ of the bundles in ($\heartsuit$)
will generate $D^b(Coh\, \mathrm{IGrass}(3 , V))$ downstairs; moreover
$\Sigma (\mathcal{R}_3^{\perp}/\mathcal{R}_3)$ is the pull back
$\pi^{\ast}\Sigma (\mathcal{R}^{\perp}/\mathcal{R})$ of the spinor
bundle $\Sigma (\mathcal{R}^{\perp}/\mathcal{R})$ on the base
$\mathrm{IGrass}(3 , V)$. When one wants to apply Bott's theorem to
calculate direct images the trouble is that $\Sigma
(\mathcal{R}_1^{\perp}/\mathcal{R}_1)$ and $\Sigma
(\mathcal{R}_2^{\perp}/\mathcal{R}_2)$, though homogeneous vector
bundles on $\mathrm{IFlag}(1 , 2 , 3; V)= \mathrm{Spin}_7 \mathbb{C} /
B$, are not defined by irreducible representations, i.e. characters of,
the Borel subgroup $B$. Therefore, one has to find Jordan-H\"older
series for these, i.e. filtrations
\[
0=\mathcal{V}_0 \subset \mathcal{V}_1 \subset \mathcal{V}_2 \subset
\dots \subset \mathcal{V}_M= \Sigma (\mathcal{R}_1^{\perp}/ \mathcal{R}_1)
\]
and 
\[
0=\mathcal{W}_0 \subset \mathcal{W}_1 \subset
\dots \subset \mathcal{W}_N= \Sigma (\mathcal{R}_2^{\perp}/ \mathcal{R}_2)
\]
by homogeneous vector subbundles $\mathcal{V}_i$ resp. $\mathcal{W}_j$
such that the quotients $\mathcal{V}_{i+1}/\mathcal{V}_i$, $i=0, \dots
, M-1$, resp. $\mathcal{W}_{j+1}/\mathcal{W}_j$, $j=0, \dots , N-1$,
are line bundles defined by characters of $B$.\\
For this, put $G:=\mathrm{Spin}_7\,\mathbb{C}$ and turning to the
notation and set-up introduced at the beginning of subsection 2.2,
rewrite diagram (D) as 
\begin{center}
\setlength{\unitlength}{1cm}
\begin{picture}(2, 4)
\put(-4 ,0){$G/P(\alpha_1)$}
\put(-4 , 1.5){$G/P(\alpha_1 , \alpha_2)$}
\put(-4 ,3){$G/P(\alpha_1, \alpha_2 ,\alpha_3)=G/B$}
\put(1, 0){$G/P(\alpha_3)$}
\put(-3, 1.3){\vector(0 , -1){0.8}} \put(-2.8, 0.7){$\pi_1$}
\put(-3, 2.8){\vector(0 ,-1){0.8}} \put(-2.8, 2.2){$\pi_2$}
\put(0, 2.8){\vector(2,-3){1.6}} \put(0.5, 1.5){$\pi$}
\put(3 , 1.5){(D')}
\end{picture}
\end{center}
\vspace{0.5cm}
In this picture, the spinor bundles $\Sigma (\mathcal{R}_i^{\perp}/ \mathcal{R}_i)$,
$i=1 ,\: 2, \: 3$, on $\mathrm{IFlag}(1, 2, 3; V)$ are the pull-backs
under the projections $G/B\to G/P(\alpha_i)$ of the vector bundles on
$G/P(\alpha_i)$ which are the duals of the homogeneous vector bundles
associated to the irreducible representations $r_i$ of $P(\alpha_i)$
with highest weight the fundamental weight $\omega_3$.\\
Recall that in terms of an orthonormal basis $\epsilon_1, \dots ,
\epsilon_r$ of $\mathfrak{h}^{\ast}$ we can write the fundamental
weights for $\mathfrak{so}_{2r+1}\,\mathbb{C}$ as
$\omega_i=\epsilon_1+\dots \epsilon_i$, $1\le i< r$,
$\omega_r=(1/2)(\epsilon_1 +\dots +\epsilon_r)$, and simple roots as
$\alpha_i=\epsilon_i-\epsilon_{i+1}$, $1\le i < r$,
$\alpha_r=\epsilon_r$, and that (cf. \cite{FuHa}, \S 20.1) the weights
of the spin representation of $\mathfrak{so}_{2r+1}\,\mathbb{C}$ are
just given by
\[
\frac{1}{2}(\pm \epsilon_1 \pm \dots \pm \epsilon_r)
\]
(all possible $2^r$ sign combinations).\\
Therefore, on the level of Lie algebras, the weights of $dr_1$,
$dr_2$, and $dr_3$ are given by:
\begin{gather*}
dr_1 :\quad \frac{1}{2}(\epsilon_1 \pm \epsilon_2 \pm \epsilon_3) ,
\quad dr_2 : \quad  \frac{1}{2}(\epsilon_1 + \epsilon_2 \pm
\epsilon_3) , \\
dr_3 : \quad \frac{1}{2}(\epsilon_1 + \epsilon_2 + \epsilon_3) .
\end{gather*}
(Indeed, if $v_{\omega_3}$ is a highest weight vector in the irreducible $G$-module of highest weight 
 $\omega_3$, then the span of $P(\alpha_i)\cdot v_{\omega_3}$, $i=1, \dots , 3$, is the irreducible
$P(\alpha_i)$-module of highest weights $\omega_3$, and its weights are therefore those weights of the ambient
irreducible $G$-module wich can be written as $\omega_3-\sum_{j\neq i} c_j \alpha_j$, $c_j\in\mathbb{Z}^+$).
\\  
Therefore, the spinor bundle $\Sigma (\mathcal{R}_3^{\perp}
/\mathcal{R}_3)$ on $G/B$ is just the line bundle
$\mathcal{L}(\omega_3)=\mathcal{L}(1/2 ,1/2 , 1/2 )$ associated to
$\omega_3$ (viewed as a character of $B$),  $\Sigma (\mathcal{R}_2^{\perp}
/\mathcal{R}_2)$ has a Jordan-H\"older filtration of length 2 with
quotients\\ $\mathcal{L}(1/2 , 1/2 , \pm  1/2)$, and  $\Sigma (\mathcal{R}_1^{\perp}
/\mathcal{R}_1)$ has a Jordan-H\"older filtration of length 4 with
quotients the line bundles $\mathcal{L}(1/2 , \pm 1/2 , \pm 1/2)$.\\
In conclusion we get that $D^b(Coh\, G/B)$ is generated by the line
bundles 
\[
(\heartsuit' )\ \quad \mathcal{A}'\otimes \mathcal{B}'\otimes \mathcal{C}'
\]
where $\mathcal{A}'$ runs through the set
\begin{gather*}
A':=\left\{ \mathcal{L}( \frac{1}{2} ,\frac{1}{2} ,
\frac{1}{2}) \otimes \mathcal{L}(-5, 0, 0) , \: \mathcal{L}( \frac{1}{2} ,-\frac{1}{2} ,
-\frac{1}{2}) \otimes \mathcal{L}(-5, 0, 0) ,\right. \\ \mathcal{L}( \frac{1}{2} ,-\frac{1}{2} ,
\frac{1}{2}) \otimes \mathcal{L}(-5, 0, 0),\: 
\mathcal{L}( \frac{1}{2} ,\frac{1}{2} ,
-\frac{1}{2}) \otimes \mathcal{L}(-5, 0, 0)     , \: \mathcal{L}(-4,
0, 0) ,\\
\left. \mathcal{L}(-3, 0, 0),\: \mathcal{L}(-2, 0, 0) , \: \mathcal{L}(-1, 0,
0), \: \mathcal{L}(0, 0, 0) \right\}   
\end{gather*}
and $\mathcal{B}'$ runs through
\begin{gather*}
B':= \left\{ \mathcal{L}( \frac{1}{2} ,\frac{1}{2} ,
\frac{1}{2}) \otimes \mathcal{L}(0, -3, 0) , \: \mathcal{L}( \frac{1}{2} ,\frac{1}{2} ,
-\frac{1}{2}) \otimes \mathcal{L}(0, -3, 0) ,  \right. \\
\left. \mathcal{L}(0, -2, 0) , \: \mathcal{L}(0, -1,
0), \: \mathcal{L}(0, 0, 0)
\right\}
\end{gather*}
and $\mathcal{C}'$ runs through
\[
C':= \left\{ \mathcal{L}( \frac{1}{2} ,\frac{1}{2} ,
\frac{1}{2}) \otimes \mathcal{L}(0, 0, -1) , \: \mathcal{L}(0, 0, 0)   \right\}  .
\]

Then we can calculate
$R\pi_{\ast}(\mathcal{A}'\otimes\mathcal{B}'\otimes\mathcal{C}')$ by
applying the relative version of Bott's theorem as explained in
subsection 3.3 to each of the 90 bundles
$\mathcal{A}'\otimes\mathcal{B}'\otimes\mathcal{C}'$; here of course
one takes into account that $\mathcal{L}(1/2, 1/2, 1/2)=\pi^{\ast} L$,
where for simplicity we denote by $L$ the line bundle on
$G/P(\alpha_3)$ defined by the one-dimensional representation of
$P(\alpha_3)$ with weight $-\omega_3$, i.e. $L=\Sigma
(\mathcal{R}^{\perp}/\mathcal{R})$, and one uses the projection
formula. After a lengthy calculation one thus arrives at the following
\begin{theorem}
The derived category $D^b(Coh\, \mathrm{IGrass}(3 , V))$ is generated
as triangulated category by the following 22 vector bundles:
\begin{gather*}
\bigwedge\nolimits^2 \mathcal{R} (-1), \:\mathcal{O}(-2),\:
\mathcal{R}(-2)\otimes L, \: \mathrm{Sym}^2 \mathcal{R}(-1)\otimes L ,
\:
\mathcal{O}(-3)\otimes L ,\\ \bigwedge\nolimits^2
\mathcal{R}(-2)\otimes L ,\: \Sigma^{2 ,1}\mathcal{R}(-1)\otimes L ,
\: \mathcal{R}(-1) ,\: \mathcal{O}(-2)\otimes L ,\: \mathcal{O}(-1)
,\\
\mathcal{R}(-1)\otimes L ,\; \bigwedge\nolimits^2
\mathcal{R}(-1)\otimes L , \: \Sigma^{2 , 1}\mathcal{R}\otimes L , \:
\mathrm{Sym}^2 \mathcal{R}^{\vee} (-2)\otimes L ,\:
\bigwedge\nolimits^2 \mathcal{R} , \: \mathcal{O} , \\
\Sigma^{2 , 1}\mathcal{R} , \: \mathrm{Sym}^2 \mathcal{R}^{\vee}(-2) ,
\: \mathcal{O}(-1)\otimes L , \: \mathrm{Sym}^2 \mathcal{R}^{\vee}
(-1) , \: \bigwedge\nolimits^2 \mathcal{R}\otimes L , \:
\mathcal{R}\otimes L .
\end{gather*}
\end{theorem}  
One should remark that the expected number of vector bundles in a
complete exceptional sequence is 8 in this case since there are 8
Schubert varieties in $\mathrm{IGrass}(3, V)$ (cf. \cite{BiLa}, \S
3). 

\section{Degeneration techniques}

Whereas in the preceding section a strategy for proving existence of complete exceptional sequences on 
rational homogeneous varieties was exposed which was based on the
method of fibering them into simpler ones of the same type, 
here we propose to explain an idea for a possibly alternative approach to tackle this problem. It relies on 
a theorem due to M. Brion that provides a degeneration of the diagonal $\Delta\subset X\times X$, $X$ 
rational homogeneous, into a union (over the Schubert varieties in $X$) of the products of a Schubert variety 
with its opposite Schubert variety.\\
We will exclusively consider the example of $\mathbb{P}^n$ and the main goal will be to compare resolutions 
of the structure sheaves of the diagonal and its degeneration product in this case. This gives a way of 
proving Beilinson's theorem on $\mathbb{P}^n$ without using a resolution of $\mathcal{O}_{\Delta}$ but only 
of the structure sheaf of the degeneration.\\

\subsection{A theorem of Brion}

The notation concerning rational homogeneous varieties introduced at the beginning of subsection 2.2 is 
retained.\\
The following theorem was proven by M. Brion (cf. \cite{Bri}, thm. 2).

\begin{theorem}
Regard the simple roots $\alpha_1,\dots ,\alpha_r$ as characters of the maximal torus $H$ and put 
\begin{gather*}
\mathfrak{X}:= \mathrm{closure}\: \mathrm{of}\: \{ (hx,\, x,\, \alpha_1(h), \dots , \alpha_r(h))\, |\, x\in 
X=G/P , \: h\in H\}\\
 \:\mathrm{in}\: X\times X\times \mathbb{A}^r 
\end{gather*}
with its projection $\mathfrak{X}\stackrel{\pi}{\longrightarrow} \mathbb{A}^r$. If $H$ acts on $\mathfrak{X}$ 
via its action on the ambient $X\times X\times \mathbb{A}^r$ given by
\[
h\cdot (x_1 ,\, x_2,\, t_1, \dots , t_r):= (hx_1 , \, x_2,\, \alpha_1(h) t_1 ,\dots , \alpha_r(h) t_r ) 
\]
and acts in $\mathbb{A}^r$ with weights $\alpha_1, \dots ,\alpha_r$, then $\pi$ is equivariant, surjective, flat
 with reduced fibres such that
\[
\mathfrak{X}_0:= \pi^{-1} ( (0,\dots , 0))\simeq \bigcup\limits_{w\in W^P} X_w\times X^w \, ,
\]
and is a trivial fibration over $H\cdot (1, \dots , 1)$, the complement of the union of all coordinate 
hyperplanes, with fibre the diagonal $\Delta=\Delta_X \subset X\times X$. 
\end{theorem}

Now the idea to use this result for our purpose is as follows: In \cite{Bei}, Beilinson proved his theorem 
using an explicit resolution of $\mathcal{O}_{\Delta_{\mathbb{P}^n}}$. However, on a general rational 
homogeneous variety $X$ a resolution of the structure sheaf of the diagonal is hard to come up with. The hope 
may  be therefore that a resolution of $\mathfrak{X}_0$ is easier to manufacture (by combinatorial 
methods) than one for $\mathcal{O}_{\Delta }$, and that one could afterwards lift the resolution of $\mathcal{
O}_{\mathfrak{X}_0}$ to one of $\mathcal{O}_{\Delta}$ by flatness.\\
If we denote by $p_1$ resp. $p_2$ the projections of $X\times X$ to the first resp. second factor, the 
preceding hope is closely connected to the problem of comparing the functors $Rp_{2 \ast}(p_1^{\ast} ( -) 
\otimes^L \mathcal{O}_{\mathfrak{X}_0} )$ and $Rp_{2 \ast}(p_1^{\ast} ( -) \otimes^L \mathcal{O}_{\Delta} )
\simeq \mathrm{id}_{D^b(Coh\, X)}$. In the next subsection we will present the computations to clarify these 
issues for projective space. 

\subsection{Analysis of the degeneration of the Beilinson functor on
$\mathbb{P}^n$}

Look at two copies of $\mathbb{P}^n$, one with homogeneous coordinates $x_0, \ldots , x_n$, the other with 
homogeneous coordinates $y_0,\ldots , y_n$. In this case $\mathfrak{X}_0=\bigcup_{i=0}^n \mathbb{P}^i \times 
\mathbb{P}^{n-i}$, and $\mathfrak{X}_0$ is defined by the ideal $J=(x_i y_j)_{0\le i< j\le n}$ and the diagonal 
by the ideal $I=(x_iy_j -x_j y_i)_{0\le i< j\le n}$. \\
Consider the case of $\mathbb{P}^1$. The first point that should be noticed is that $Rp_{2 \ast}(p_1^{\ast} ( -) 
\otimes^L \mathcal{O}_{\mathfrak{X}_0} )$ is no longer isomorphic to the identity: By Orlov's representability 
theorem (cf. \cite{Or2}, thm. 3.2.1) the identity functor is represented uniquely by the structure sheaf of the 
diagonal on the product. Here one can also see this in an easier way as follows. For $d>> 0$ the sheaf 
$p_1^{\ast} \mathcal{O}(d) \otimes \mathcal{O}_{\mathfrak{X}_0}$ is $p_{2\ast}$-acyclic and $p_{2\ast}$ 
commutes with base extension whence $\dim_{\mathbb{C}}\left(  p_{2\ast}( p_1^{\ast}\mathcal{O}(d) \otimes 
\mathcal{O}_{\mathfrak{X}_0}) \otimes \mathbb{C}_P\right) =d+1$ if $P$ is the point $\{ y_1=0 \}$ and $=1$ 
otherwise:
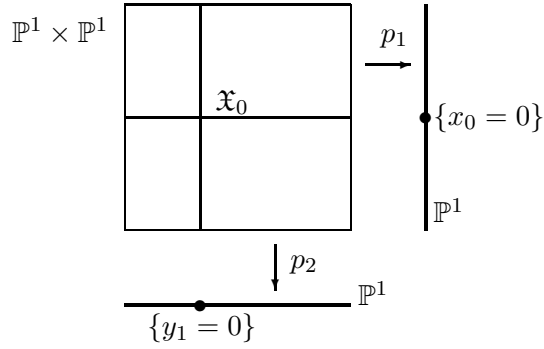
\begin{center} 
\setlength{\unitlength}{1cm}
\begin{picture}(2,4.5)
\put(-2,1){\line(1,0){3}}
\put(-2,1){\line(0,1){3}}
\put(-2,4){\line(1,0){3}}
\put(1,4){\line(0,-1){3}}

\thicklines
\put(-2,0){\line(1,0){3}}
\put(2,4){\line(0,-1){3}}

\put(-1,1){\line(0,1){3}}
\put(-2,2.5){\line(1,0){3}}

\thinlines
\put(-1.1,-0.1){$\bullet$}
\put(1.9, 2.4){$\bullet$}
\put(-3.5,3.5){$\mathbb{P}^1\times\mathbb{P}^1$}
\put(0,0.8){\vector(0,-1){0.6}}
\put(0.2, 0.5){$p_2$}
\put(1.2,3.2){\vector(1,0){0.6}}
\put(1.4, 3.5){$p_1$}
\put(-0.8,2.6){$\mathfrak{X}_0$}
\put(2.1,1.1){$\mathbb{P}^1$}
\put(1.1,0){$\mathbb{P}^1$}
\put(2.1, 2.4){$\{ x_0=0\}$}
\put(-1.7,-0.4){$\{ y_1=0\}$}
\end{picture}
\end{center}
\vspace{0.5cm}
Thus $p_{2\ast}( p_1^{\ast}\mathcal{O}(d) \otimes 
\mathcal{O}_{\mathfrak{X}_0})=Rp_{2\ast}( p_1^{\ast}\mathcal{O}(d) \otimes^L 
\mathcal{O}_{\mathfrak{X}_0})$ is not locally free in this case. We will give a complete description 
of the functor $Rp_{2\ast}( p_1^{\ast}(-) \otimes^L 
\mathcal{O}_{\mathfrak{X}_0})$ below for $\mathbb{P}^n$. 
If one compares the resolutions of $\mathcal{O}_{\mathfrak{X}_0}$ and $\mathcal{O}_{\Delta}$
on $\mathbb{P}^1$:
\begin{gather*}
\begin{CD}
0 @>>> \mathcal{O}(-1,-1) @> (x_0y_1)>> \mathcal{O}_{\mathbb{P}^1\times\mathbb{P}^1} @>>>
\mathcal{O}_{\mathfrak{X}_0} @>>> 0
\end{CD}\\
\begin{CD}
0 @>>> \mathcal{O}(-1,-1) @> (x_0y_1-x_1y_0)>> \mathcal{O}_{\mathbb{P}^1\times\mathbb{P}^1} @>>>
\mathcal{O}_{\Delta} @>>> 0
\end{CD}
\end{gather*}
and on $\mathbb{P}^2$:
\begin{gather*}
0 \to \mathcal{O}(-2, -1)\oplus \mathcal{O}(-1,-2)\stackrel{A'}{\longrightarrow}
\mathcal{O}(-1,-1)^{\oplus 3}  \stackrel{B'}{\longrightarrow}
\mathcal{O}_{\mathbb{P}^2\times\mathbb{P}^2}
\to
\mathcal{O}_{\mathfrak{X}_0} \to 0\\
\begin{CD}
0 \to \mathcal{O}(-2, -1)\oplus \mathcal{O}(-1,-2) \stackrel{A}{\longrightarrow}
\mathcal{O}(-1,-1)^{\oplus 3} \stackrel{B}{\longrightarrow}
\mathcal{O}_{\mathbb{P}^2\times\mathbb{P}^2} \to
\mathcal{O}_{\Delta} \to 0
\end{CD}
\end{gather*}
where
\begin{gather*}
A=\left(\begin{array}{cc} x_0 & y_0 \\
x_1 & y_1 \\ x_2 & y_2
\end{array}\right) \quad
A'=\left( \begin{array}{cc} x_0 & 0 \\
x_1 & y_1 \\ 0 & y_2
\end{array}\right) \\
B=(x_2y_1-x_1y_2 , \: x_0y_2 - x_2 y_0 ,\: x_1 y_0 -x_0 y_1 ) \quad
B=(-x_1y_2 , \: x_0y_2  ,\:  -x_0 y_1 )
\end{gather*}
(these being Hilbert-Burch type resolutions; here $\mathfrak{X}_0$ is no longer
a local complete intersection!) one may wonder if on $\mathbb{P}^n$ there exist resolutions
of $\mathcal{O}_{\mathfrak{X}_0}$ and $\mathcal{O}_{\Delta}$ displaying an analogous
similarity. This is indeed the case, but will require some work.\\
Consider the matrix
\[
\left(
\begin{array}{ccc}
x_0 & \dots & x_n \\
y_0 & \dots & y_n
\end{array}
\right)
\]
as giving rise to a map between free bigraded modules $F$ and $G$ over \\$\mathbb{C}[x_0, \dots ,
x_n;\: y_0 , \dots , y_n]$ of rank $n+1$ and $2$ respectively. Put $K_h:= \bigwedge^{h+2}F
\otimes \mathrm{Sym}^h\, G^{\vee}$ for $h=0, \dots , n-1$. Choose bases $f_0, \dots , f_n$
resp. $\xi ,\: \eta$ for $F$ resp. $G^{\vee}$. Define maps $d_h : K_h\to K_{h-1}$,
$h=1, \ldots, n-1$ by 
\begin{gather*}
d_h \left( f_{j_1}\wedge \dots \wedge f_{j_{h+2}} \otimes
\xi^{\mu_1}\eta^{\mu_2}\right):=\sum\limits_{l=1}^{h+2} (-1)^{l+1} x_{j_l} f_{j_1}\wedge\dots
\wedge \hat{f}_{j_l}\wedge \dots \wedge f_{j_{h+2}} \\
\otimes \xi^{-1}(\xi^{\mu_1}\eta^{\mu_2}) + \sum\limits_{l=1}^{h+2} (-1)^{l+1} y_{j_l} f_{j_1}
\wedge\dots\wedge \hat{f}_{j_l}\wedge \dots \wedge f_{j_{h+2}}\otimes
\eta^{-1}(\xi^{\mu_1}\eta^{\mu_2})
\end{gather*}
where $0\le j_1 < \dots < j_{h+2} \le n$, $\mu_1+\mu_2=h$ and the homomorphism $\xi^{-1}$ (resp.
$\eta^{-1}$) is defined by
\[
\xi^{-1} (\xi^{\mu_1} \eta^{\mu_2}):= \left\{ \begin{array}{l}
\xi^{\mu_1-1} \eta^{\mu_2} \quad \mathrm{if}\quad \mu_1\ge 1 \\
\quad 0 \quad\quad\;\quad \mathrm{if}\quad \mu_1=0
\end{array}\right.
\]
(resp. analogously for $\eta^{-1}$). Then 
\[
\begin{CD}
0 @>>> K_{n-1} @>d_{n-1}>> \dots @>d_1>> K_0 @>>> I @>>> 0
\end{CD}
\]
is a resolution of $I$ which is the Eagon-Northcott complex in our special case (cf.
\cite{Nor}, appendix C).
 
\begin{proposition}
The ideal $J$ has a resolution
\[
\begin{CD}
0 @>>> K_{n-1} @>d'_{n-1}>> \dots @>d'_1>> K_0 @>>> J @>>> 0
\end{CD}
\]
where the differential $d'_h : K_h\to K_{h-1}$ is defined by
\begin{gather*}
d'_h \left( f_{j_1}\wedge \dots \wedge f_{j_{h+2}} \otimes
\xi^{\mu_1}\eta^{\mu_2}\right):=\sum\limits_{l=1}^{h-\mu_2+1} (-1)^{l+1} x_{j_l}
f_{j_1}\wedge\dots
\wedge \hat{f}_{j_l}\wedge \dots \wedge f_{j_{h+2}} \\
\otimes \xi^{-1}(\xi^{\mu_1}\eta^{\mu_2}) + \sum\limits_{l=\mu_1+2}^{h+2} (-1)^{l+1} y_{j_l}
f_{j_1}
\wedge\dots\wedge \hat{f}_{j_l}\wedge \dots \wedge f_{j_{h+2}}\otimes
\eta^{-1}(\xi^{\mu_1}\eta^{\mu_2})\, .
\end{gather*}
\end{proposition}

Intuitively the differentials $d'_h$ are gotten by degenerating the differentials $d_h$. To
prove proposition 4.2.1 we will use the fact that $J$ is a monomial ideal. There is a
combinatorial method  for sometimes writing down resolutions for these by looking at simplicial
or more general cell complexes from topology. The method can be found in \cite{B-S}. We
will recall the results we need in the following. Unfortunately the resolution of proposition
4.2.1 is not supported on a simplicial complex, one needs a more general cell complex.\\
Let $X$ be a finite regular cell complex. This is a non-empty topological space $X$ with a
finite set $\Gamma$ of subsets of $X$ (the \emph{cells} of $X$) such that
\begin{itemize}
\item[(a)]
$X=\bigcup\limits_{e\in\Gamma} e$,
\item[(b)]
the $e\in\Gamma$ are pairwise disjoint,
\item[(c)]
$\emptyset\in\Gamma$,
\item[(d)]
for each non-empty $e\in\Gamma$ there is a homeomorphism between a
closed $i$-dimensional ball and the closure $\bar{e}$ which maps the interior of the ball onto
$e$ (i.e. $e$ is an \emph{open} $i$-\emph{cell}).
\end{itemize}
We will also call the $e\in \Gamma$ \emph{faces}. We will say that $e'\in \Gamma$ is a \emph{face
of}
$e\in
\Gamma$, $e\neq e'$, or that $e$ contains $e'$ if $e'\subset \bar{e}$. The maximal faces of $e$
under containment are called its
\emph{facets}. $0$- and $1$-dimensional faces will be called \emph{vertices} and \emph{edges}
respectively. The set of vertices is denoted $\mathfrak{V}$. A subset $\Gamma'\subset\Gamma$
such that for each $e\in\Gamma'$ all the faces of $e$ are in $\Gamma'$ determines a
\emph{subcomplex} 
$X_{\Gamma'}=\bigcup_{e\in\Gamma'} e$ of $X$. Moreover we assume in addition
\begin{itemize}
\item[(e)]
If $e'$ is a codimension $2$ face of $e$ there are exactly two facets $e_1$, $e_2$ of $e$
containing $e'$.
\end{itemize}
The prototypical example of a finite regular cell complex is the set
of faces of a convex polytope for which property (e) is fulfilled. In
general (e) is added as a kind of regularity assumption.\\
Choose an \emph{incidence function} $\epsilon(e , e')$ on pairs of faces of $e$, $e'$. This
means that $\epsilon$ takes values in $\{ 0 , +1 ,-1\}$, $\epsilon (e ,e')=0$ unless $e'$ is a
facet of $e$, $\epsilon (v,\emptyset )=1$ for all vertices $v\in\mathfrak{V}$ and moreover
\[
\epsilon (e ,e_1)\epsilon (e_1, e') +\epsilon (e, e_2) \epsilon (e_2 , e') =0
\] 
for $e$, $e_1$, $e_2$, $e'$ as in (e).\\
Let now $M=(m_v)_{v\in\mathfrak{V}}$ be a monomial ideal ($m_v$ monomials) in the polynomial
ring $k[T_1,\dots , T_N]$, $k$ some field. For multi-indices $\underline{a}$,
$\underline{b}\in\mathbb{Z}^N$ we write $\underline{a}\le \underline{b}$ to denote $a_i\le b_i$
for all $i=1,\dots ,N$. $T^{\underline{a}}$ denotes $T_1^{a_1}\cdot\dots \cdot T_N^{a_N}$.\\
The oriented chain complex $\tilde{C}(X, k)=\bigoplus_{e\in \Gamma} k e$ (the
homological grading is given by dimension of faces) with differential 
\[
\partial e:= \sum\limits_{e'\in\Gamma} \epsilon (e, e') e'
\]
computes the \emph{reduced cellular homology groups} $\tilde{H}^i(X , k)$ of $X$.\\
Think of the vertices $v\in\mathfrak{V}$ as labelled by the corresponding monomials $m_v$. Each
non-empty face $e\in\Gamma$ will be identified with its set of vertices and will be labelled
by the least common multiple $m_e$ of its vertex labels. The \emph{cellular complex} $F_{X, M}$
associated to $(X, M)$ is the $\mathbb{Z}^N$-graded $k[T_1,\ldots ,T_N]$-module
$\bigoplus_{e\in\Gamma ,e\neq \emptyset} k[T_1, \ldots , T_N] e$ with differential
\[
\partial e:= \sum\limits_{e'\in\Gamma , e'\neq \emptyset} \epsilon (e, e') \frac{m_e}{m_{e'}}e'
\] 
(where again the homological grading is given by the dimension of the faces).\\
For each multi-index $\underline{b}\in\mathbb{Z}^N$ let $X_{\le \underline{b}}$ be the subcomplex
of $X$ consisting of all the faces $e$ whose labels $m_e$ divide $T^{\underline{b}}$. We have

\begin{proposition}
$F_{X, M}$ is a free resolution of $M$ if and only if $X_{\le \underline{b}}$ is acyclic over
$k$ for all $\underline{b}\in\mathbb{Z}^N$ (i.e.$\tilde{H}_i(X_{\le
\underline{b}}, k)=0$ for all $i$).
\end{proposition}

We refer to \cite{B-S}, prop. 1.2, for a proof.\\
\\
Next we will construct appropriate cell complexes $Y^n$, $n=1, \; 2, \dots $, that via the
procedure described above give resolutions of $J=(x_iy_j)_{0\le i<j\le n}$. We will apply
proposition 4.2.2 by showing that for all $\underline{b}\in\mathbb{Z}^{2n+2}$ the subcomplexes
$Y^n_{\le \underline{b}}$ are contractible.\\
It is instructive to look at the pictures of $Y^1$, $Y^2$, $Y^3$, $Y^4$ with their labellings
first: 
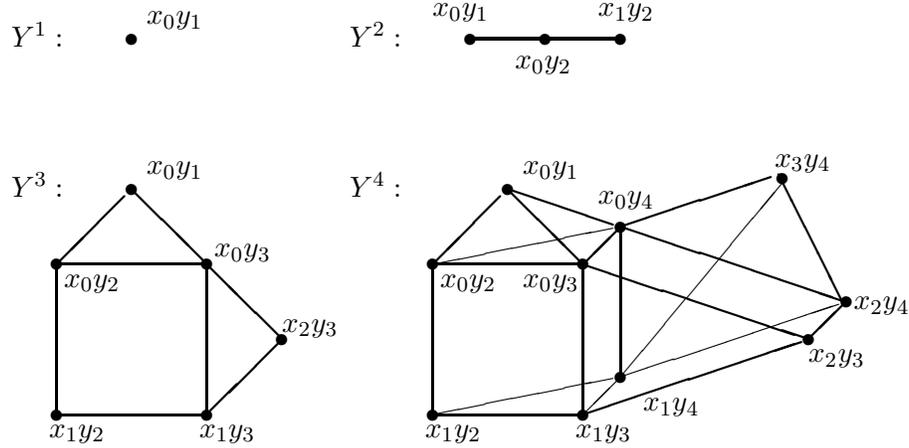
\setlength{\unitlength}{1cm}
\begin{center}
\begin{picture}(5,8)
\thicklines
\put(-3.5,7){$Y^1:$} \put(-2, 7){$\bullet$} \put(-1.7,7.3){$x_0y_1$}

\put(1 ,7){$Y^2:$} \put(2.5,7){$\bullet$} \put(2.6,7.1){\line(1,0){1}}
\put(3.5,7){$\bullet$}  \put(3.6,7.1){\line(1,0){1}} \put(4.5,7){$\bullet$}
\put(2.1, 7.4){$x_0y_1$}
\put(3.2, 6.7){$x_0y_2$} \put(4.3,7.4){$x_1y_2$}

\put(-3.5,5){$Y^3:$}
\put(-2,5){$\bullet$}\put(-1.7,5.3){$x_0y_1$}
\put(-3,4){$\bullet$}  \put(-1, 4){$\bullet$}
\put(-1.85,5.05){\line(1,-1){0.9}} \put(-2,5){\line(-1,-1){0.9}}
\put(-2.9,4.1){\line(1,0){2}} \put(-2.9, 4.1){\line(0,-1){2}}\put(-0.9, 4.1){\line(0,-1){2}}
\put(-3,2){$\bullet$}  \put(-1, 2){$\bullet$}
\put(-2.9,2.1){\line(1,0){2}}
\put(0, 3){$\bullet$}
\put(0.1,3.1){\line(-1,-1){0.95}} \put(0.1,3.1){\line(-1,1){0.95}} 
\put(-2.8,3.8){$x_0y_2$} \put(-0.8,4.2){$x_0y_3$}
\put(-3,1.8){$x_1y_2$}  \put(-1, 1.8){$x_1y_3$}
\put(0.1, 3.2){$x_2y_3$}

\put(1,5){$Y^4:$}
\put(3,5){$\bullet$}\put(3.3,5.3){$x_0y_1$}
\put(2,4){$\bullet$}  \put(4, 4){$\bullet$}
\put(3.15,5.05){\line(1,-1){0.9}} \put(3,5){\line(-1,-1){0.9}}
\put(2.1,4.1){\line(1,0){2}} \put(2.1, 4.1){\line(0,-1){2}}\put(4.1, 4.1){\line(0,-1){2}}
\put(2,2){$\bullet$}  \put(4, 2){$\bullet$}
\put(2.1,2.1){\line(1,0){2}}
\put(2.2,3.8){$x_0y_2$} \put(3.3,3.8){$x_0y_3$}
\put(2,1.8){$x_1y_2$}  \put(4, 1.8){$x_1y_3$}
\put(4.1,4.1){\line(1,1){0.5}}
\put(4.5,4.5){$\bullet$}
\put(4.5,2.5){$\bullet$}
\put(7.1,3.1){\line(1,1){0.5}}
\put(7.5,3.5){$\bullet$}
\put(7,3){$\bullet$}
\put(4.1,4.1){\line(3,-1){3}}
\put(4.6,4.6){\line(3,-1){3}}
\put(4.1,2.1){\line(3,1){2.9}}
\put(3.1,5.1){\line(3,-1){1.4}}
\put(6.65,5.15){$\bullet$}
\put(4.6 , 4.6){\line(3,1){2}}
\put(6.7,5.3){\line(1,-2){0.85}}
\thinlines
\put(4.6,2.6){\line(5,6){2.16}}
\put(4.6,2.6){\line(0,1){2}}
\put(4.6,2.6){\line(3,1){2.9}}
\put(2.1,4.1){\line(5,1){2.4}}
\put(2.1,2.1){\line(5,1){2.4}}
\put(4.1,2.1){\line(1,1){0.5}}
\put(4.3,4.9){$x_0y_4$}
\put(4.9,2.15){$x_1y_4$}
\put(7.7,3.5){$x_2y_4$}
\put(7.1,2.8){$x_2y_3$}
\put(6.65,5.4){$x_3y_4$}
\end{picture}
\end{center}
\vspace{0.3cm}

The general procedure for constructing $Y^n$ geometrically is as follows: In $\mathbb{R}^{n-1}$
take the standard $(n-1)$-simplex $P^1$ on the vertex set $\{x_0y_1,\: x_0y_2, \ldots ,$ $\:
x_0y_n\}$. Then take an $(n-2)$-simplex on the vertex set $\{x_1y_2, \ldots , x_1y_n\}$, viewed
as embedded in the same $\mathbb{R}^{n-1}$, and join the vertices $x_1y_2, \ldots , x_1y_n$,
respectively, to the vertices $x_0y_2,\ldots , x_0y_n$,, respectively, of $P^1$ by drawing an
edge between $x_0y_i$ and $x_1y_i$ for $i=2, \ldots , n$. This describes the process of
attaching a new $(n-1)$-dimensional polytope $P^2$ to the facet of $P^1$ on the vertex set
$\{x_0y_2 , \ldots , x_0y_n \}$.\\
Assume that we have constructed inductively the $(n-1)$-dimensional polytope $P^i$, $2\le i\le
n-1$, with one facet on the vertex set $\{ x_{\mu} y_{\nu} \}_{\stackrel{0\le \mu\le
i-1}{i+1\le \nu\le n}}$. Then take an $(n-i-1)$-simplex on the vertex set $\{ x_iy_{i+1},
\ldots , x_iy_n\}$, viewed as embedded in the same $\mathbb{R}^{n-1}$, and for every $\alpha$
with $1\le \alpha \le n-i$ join the vertex $x_iy_{i+\alpha}$ of this simplex to the vertices
$x_{\mu}y_{i+\alpha}$, $0\le\mu\le i-1$, of $P^i$ by an edge. This corresponds to attaching a
new $(n-1)$-dimensional polytope $P^{i+1}$ to the facet of $P^i$ on the vertex set $\{ x_{\mu}
 y_{\nu} \}_{\stackrel{0\le \mu\le i-1}{i+1\le \nu\le n}}$.\\
In the end we get $(n-1)$-dimensional polytopes $P_1,\dots , P_n$ in $\mathbb{R}^{n-1}$
 where $P_j$ and $P_{j+1}$, $j=1, \dots ,n-1$, are glued along a common facet. These will make
up our labelled cell complex $Y^n$.\\
The $h$-dimensional faces of $Y^n$ will be called $h$\emph{-faces} for short. We need a more
convenient description for them:

\begin{lemma}
There are natural bijections between the following sets:
\begin{itemize}
\item[(i)]
$\{ h$-faces of $Y^n\}$
\item[(ii)]
matrices
\[
\left( \begin{array}{ccc}
x_{i_1}y_{i_{\mu_1+2}} & \cdots  & x_{i_1} y_{i_{h+2}} \\
\vdots            &   \ddots    &   \vdots             \\
x_{i_{\mu_1+1}} y_{i_{\mu_1+2}} & \cdots  & x_{i_{\mu_1+1}} y_{i_{h+2}}
\end{array}
\right)
\]
with vertex labels of $Y^n$ as entries, where $0 \le i_1 < i_2 < \ldots < i_{\mu_1+1} <
i_{\mu_1+2} < \ldots < i_{h+2} \le n$ and $0\le \mu_1 \le h$ are integers. The $(\kappa ,
\lambda )$-entry of the above matrix is thus $x_{i_{\kappa}}y_{i_{\mu_1+\lambda+1}}$.
\item[(iii)]
standard basis vectors
\begin{gather*}
f_{i_1}\wedge \ldots \wedge f_{i_{\mu_1+1}}\wedge f_{i_{\mu_1+2}}\wedge \ldots \wedge
f_{i_{h+2}}\otimes \xi^{\mu_1}\eta^{\mu_2} ,  \quad 0\le \mu_1 \le h , \\
0 \le i_1 < i_2 < \ldots < i_{\mu_1+1} <
i_{\mu_1+2} < \ldots < i_{h+2} \le n
\end{gather*}
of $\bigwedge^{h+2}F\otimes \mathrm{Sym}^h \, G^{\vee}$.
\end{itemize} 
\end{lemma}

\begin{proof}
The bijection between the sets in $(ii)$ and $(iii)$ is obvious: To $f_{i_1}\wedge \ldots \wedge
f_{i_{h+2}}\otimes \xi^{\mu_1}\eta^{\mu_2}$ in $\bigwedge^{h+2}F\otimes \mathrm{Sym}^h \,
G^{\vee}$ one associates the matrix
\[
\left( \begin{array}{ccc}
x_{i_1}y_{i_{\mu_1+2}} & \cdots  & x_{i_1} y_{i_{h+2}} \\
\vdots            &   \ddots    &   \vdots             \\
x_{i_{\mu_1+1}} y_{i_{\mu_1+2}} & \cdots  & x_{i_{\mu_1+1}} y_{i_{h+2}}
\end{array}
\right) .
\]
To set up a bijection between the sets under (i) and (ii) the idea is to identify an $h$-face
$e$ of $Y^n$ with its vertex labels and collect the vertex labels in a matrix of the form given
in (ii). We will prove by induction on $j$ that the $h$-faces $e$ contained in the polytopes
$P_1 ,\ldots , P_j$ are exactly those whose vertex labels may be collected in a matrix of the
form written in (ii) satisfying the additional property $i_{\mu_1+1}\le j-1$. This will prove
the lemma.\\
$P_1$ is an $(n-1)$-simplex on the vertex set $\{ x_0y_1 ,\ldots , x_0y_n\}$ and its $h$-faces
$e$ can be identified with the subsets of cardinality $h+1$ of $\{ x_0y_1 ,\ldots , x_0y_n\}$.
We can write such a subset in matrix form
\[
\left( x_0y_{i_2} \; x_0y_{i_3} \; \cdots \; x_0y_{i_{h+2}} \right)
\]
with $0\le i_2 < i_3 < \ldots < i_{h+2} \le n$. This shows that the preceding claim is true for
$j=1$.\\
For the induction step assume that the $h$-faces of $Y^n$ contained in $P_1, \ldots , P_j$
are exactly those whose vertex labels may be collected in a matrix as in (ii) with
$i_{\mu_1+1}\le j-1$. Look at the $h$-faces $e$ contained in $P_1, \ldots , P_{j+1}$. If $e$ is
contained in $P_1,\ldots , P_j$ (which is equivalent to saying that none of the vertex labels
of $e$ involves the indeterminate $x_j$) then there is nothing to show. Now there are two types
of $h$-faces contained in $P_1, \ldots , P_{j+1}$ but not in $P_1,\ldots , P_j$: The first type
corresponds to $h$-faces $e$ entirely contained in the simplex on the vertex set $\{x_jy_{j+1},
\ldots , x_j y_n\}$. These correspond to matrices 
\[
\left( x_jy_{i_2} \; \cdots \; x_jy_{i_{h+2}} \right) ,
\]
$0\le i_2 < i_3 < \ldots < i_{h+2} \le n$, of the form given in (ii) which involve the
indeterminate $x_j$ and have only one row.\\
The second type of $h$-faces $e$ is obtained as follows: We take an $(h-1)$-face $e'$ contained
in the facet on the vertex set $\{ x_ay_b\}_{\stackrel{0\le a\le j-1}{j+1\le b\le n}}$ which
$P_j$ and $P_{j+1}$ have in common; by induction $e'$ corresponds to a matrix
\[
\left( \begin{array}{ccc}
x_{i_1}y_{i_{\mu_1+2}} & \cdots  & x_{i_1} y_{i_{h+1}} \\
\vdots            &   \ddots    &   \vdots             \\
x_{i_{\mu_1+1}} y_{i_{\mu_1+2}} & \cdots  & x_{i_{\mu_1+1}} y_{i_{h+1}}
\end{array}
\right)
\]
with $0\le i_1 < i_2 < \ldots < i_{\mu_1+1} \le j-1$ and $j+1\le i_{\mu_1+2} < \ldots < i_{h+1}
\le n$, $0\le \mu_1 \le h-1$. Then by construction of $P_{j+1}$ there is a unique $h$-face $e$
in $P_1,\ldots , P_{j+1}$, but not in $P_1, \ldots , P_j$, which contains the $(h-1)$-face $e'$: 
It is the $h$-face whose vertex labels are the entries of the preceding matrix together with
$\{ x_jy_{i_{\mu_1+2}} ,\: x_j y_{i_{\mu_1+3}} ,\ldots , \: x_j y_{i_{h+1}} \}$. Thus $e$
corresponds to the matrix
\[
\left( \begin{array}{ccc}
x_{i_1}y_{i_{\mu_1+2}} & \cdots  & x_{i_1} y_{i_{h+1}} \\
\vdots            &   \ddots    &   \vdots             \\
x_{i_{\mu_1+1}} y_{i_{\mu_1+2}} & \cdots  & x_{i_{\mu_1+1}} y_{i_{h+1}}\\
x_j y_{i_{\mu_1+2}} & \cdots   &   x_j y_{i_{h+1}} 
\end{array}
\right) .
\]
This proves the lemma.
\end{proof}

Now we want to define an incidence function $\epsilon (e, e')$ on pairs of faces $e,\: e'$ of
$Y^n$. Of course if $e'$ is not a facet of $e$ , we put $\epsilon (e , e')=0$ and likewise put
$\epsilon (v ,\emptyset ):= 1$ for all vertices $v$ of $Y^n$. Let now $e$ be an $h$-face.
Using lemma 4.2.3 it corresponds to a matrix
\[
M(e)= \left( \begin{array}{ccc}
x_{i_1}y_{i_{\mu_1+2}} & \cdots  & x_{i_1} y_{i_{h+2}} \\
\vdots            &   \ddots    &   \vdots             \\
x_{i_{\mu_1+1}} y_{i_{\mu_1+2}} & \cdots  & x_{i_{\mu_1+1}} y_{i_{h+2}}
\end{array}
\right) .
\]
A facet $e'$ of $e$ corresponds to a submatrix $M(e')$ of $M(e)$ obtained from $M(e)$ by erasing
either a row or a column. We define $\epsilon (e, e'):= (-1)^l$ if $M(e')$ is obtained from
$M(e)$ by erasing the $l$th row; we define $\epsilon (e ,e'):= (-1)^{\mu_1+j}$ if $M(e')$ is
obtained from $M(e)$ by erasing the $j$th column.\\
One must check that then $\epsilon (e, e_1)\epsilon (e_1 , e'') +\epsilon (e, e_2)\epsilon (e_2
, e'') =0$ for a codimension $2$ face $e''$ of $e$ and $e_1$, $e_2$ the two facets of $e$
containing $e''$. This is now a straightforward computation. There are 3 cases: The matrix
$M(e'')$ is obtained from $M(e)$ by (i) deleting two rows, (ii) deleting two columns, (iii)
erasing one row and one column: 
\begin{itemize}
\item[(i)]
Let $l_1< l_2$ and assume that $M(e_1)$ is $M(e)$ with $l_1$th row erased and $M(e_2)$ is $M(e)$
with $l_2$th row erased. Then
\begin{gather*}
\epsilon (e ,e_1)= (-1)^{l_1} , \quad \epsilon (e ,e_2)= (-1)^{l_2} , \quad \epsilon (e_1 ,e'')=
(-1)^{l_2-1}\\ \epsilon (e_2 ,e'')= (-1)^{l_1} .
\end{gather*}

\item[(ii)]
This is the same computation as for (i) with the roles of rows and columns interchanged.

\item[(iii)]
Assume that $M(e_1)$ is $M(e)$ with $l$th row erased and $M(e_2)$ is $M(e)$
with $j$th column erased. Then
\begin{gather*}
\epsilon (e ,e_1)= (-1)^{l} , \quad \epsilon (e ,e_2)= (-1)^{\mu_1+j} , \quad \epsilon (e_1
,e'')= (-1)^{\mu_1-1+j}\\ \epsilon (e_2 ,e'')= (-1)^{l} .
\end{gather*}

\end{itemize}

Thus $\epsilon$ is an incidence function on $Y^n$. Now one has to compute the cellular complex
$F_{Y^n , J}$: Indeed by lemma 4.2.3 we know that its term in homological degree $h$ identifies
with $\bigwedge^{h+2}F\otimes \mathrm{Sym}^h \, G^{\vee}$. If $e$ is an $h$-face recall that the
differential $\partial$ of $F_{Y^n , J}$ is given by
\[
\partial e = \sum\limits_{e'\; \mathrm{a}\; \mathrm{facet}\; \mathrm{of}\; e ,\; e'\neq
\emptyset} \epsilon(e , e') \frac{m_e}{m_{e'}} e'
\]
and if $e$ corresponds to $f_{i_1}\wedge\ldots\wedge f_{i_{h+2}}\otimes
\xi^{\mu_1}\eta^{\mu_2}\in \bigwedge^{h+2}F\otimes \mathrm{Sym}^h \, G^{\vee}$ we find 

\begin{gather*}
\partial \left( f_{j_1}\wedge \dots \wedge f_{j_{h+2}} \otimes
\xi^{\mu_1}\eta^{\mu_2}\right)=\sum\limits_{l=1}^{h-\mu_2+1} (-1)^{l+1} x_{j_l}
f_{j_1}\wedge\dots
\wedge \hat{f}_{j_l}\wedge \dots \wedge f_{j_{h+2}} \\
\otimes \xi^{-1}(\xi^{\mu_1}\eta^{\mu_2}) + \sum\limits_{l=\mu_1+2}^{h+2} (-1)^{l+1} y_{j_l}
f_{j_1}
\wedge\dots\wedge \hat{f}_{j_l}\wedge \dots \wedge f_{j_{h+2}}\otimes
\eta^{-1}(\xi^{\mu_1}\eta^{\mu_2})\, .
\end{gather*}

Thus the complex $F_{Y^n , J}$ is nothing but the complex in proposition 4.2.1. Thus to prove
proposition 4.2.1 it is sufficient in view of proposition 4.2.2 to prove the following

\begin{lemma}
For all $\underline{b}\in\mathbb{Z}^{2n+2}$ the subcomplexes $Y^n_{\le\underline{b}}$ of $Y^n$
are contractible. 
\end{lemma}

\begin{proof}
Notice that it suffices to prove the following: If 
\[
x_{i_1}\dots x_{i_k}y_{j_1}\dots y_{j_l}\quad 0\le i_1 < \dots i_k \le n , \quad 0\le j_1 < \dots < j_l \le n
\]
is a monomial that is the least common multiple of some subset of the vertex labels of $Y^n$ then the 
subcomplex $\tilde{Y}^n$ of $Y^n$ that consists of all the faces $e$ whose label divides $x_{i_1}\dots 
x_{i_k}y_{j_1}\dots y_{j_l}$ is contractible. This can be done as follows:\\ Put $\kappa (i_d):= \mathrm{min}\, 
\{ t\, : \, j_t > i_d \}$ for $d=1 , \dots ,k$. Note that we have $\kappa (i_1)=1$ and $\kappa (i_1) \le
 \kappa (i_2) \le \dots \le \kappa (i_k)$. Choose a retraction of the face $e^0$ of $\tilde{Y}^n$ corresponding
 to the matrix
\[
\left( \begin{array}{ccc}
x_{i_1}y_{j_{\kappa(i_k)}} & \dots &  x_{i_1}y_{j_l} \\
\vdots &  \ddots  & \vdots \\
x_{i_k}y_{j_{\kappa(i_k)}} & \dots &  x_{i_k}y_{j_l}
\end{array}\right)
\]
onto its facet $e^{0 '}$ corresponding to 
\[
\left( \begin{array}{ccc}
x_{i_1}y_{j_{\kappa(i_k)}} & \dots &  x_{i_1}y_{j_l} \\
\vdots &  \ddots  & \vdots \\
x_{i_{k-1}}y_{j_{\kappa(i_k)}} & \dots &  x_{i_{k-1}}y_{j_l}
\end{array}\right) .
\]
Then choose a retraction of the face $e^1$ corresponding to 
\[
\left( \begin{array}{ccc}
x_{i_1}y_{j_{\kappa(i_{k-1})}} & \dots &  x_{i_1}y_{j_l} \\
\vdots &  \ddots  & \vdots \\
x_{i_{k-1}}y_{j_{\kappa(i_{k-1})}} & \dots &  x_{i_{k-1}}y_{j_l}
\end{array}\right) 
\]
onto its facet $e^{1 '}$ corresponding to 
\[
\left( \begin{array}{ccc}
x_{i_1}y_{j_{\kappa(i_{k-1})}} & \dots &  x_{i_1}y_{j_l} \\
\vdots &  \ddots  & \vdots \\
x_{i_{k-2}}y_{j_{\kappa(i_{k-1})}} & \dots &  x_{i_{k-2}}y_{j_l}
\end{array}\right) .
\]
Notice that $e^{0 '}$ is contained in $e^1$. Continuing this pattern, one can finally retract the face 
 corresponding to 
\[
\left( \begin{array}{cccc}
x_{i_1}y_{j_{\kappa(i_{1})}} & x_{i_1}y_{j_{\kappa(i_{1})+1}}  & \dots &  x_{i_1}y_{j_l} \\
\end{array}\right) ,
\]
i.e. a simplex, onto one of its vertices. Composing these retractions, one gets a retraction of $\tilde{Y}^n$ 
onto a point. 
\end{proof}

In conclusion what we get from proposition 4.2.1 is that on $\mathbb{P}^n\times \mathbb{P}^n$ the sheaf 
$\mathcal{O}_{\mathfrak{X}_0}$ has a resolution
\begin{gather*}
\begin{CD}
(\ast) \quad  0 @>>> \bigoplus_{\stackrel{i+j=n-1}{i,j\ge 0}}\mathcal{O}(-i-1, -j-1) @>d_{n-1}'>> \dots
\end{CD}\quad\quad\quad\quad \\\begin{CD}\dots  @>d_{h+1}'>> \left(\bigoplus_{\stackrel{i+j=h}{i, j\ge 0}}
\mathcal{O}(-i-1, -j-1)\right)^{\oplus {n+1 \choose h+2} } @> d_h'>> \dots 
\end{CD}\\
\begin{CD}
\dots @>d_1' >> \mathcal{O}(-1 ,-1)^{\oplus {n+1 \choose 2} } @>>> \mathcal{O}_{\mathbb{P}^n\times\mathbb{P}^n}
@>>> \mathcal{O}_{\mathfrak{X}_0} @>>> 0
\end{CD}
\end{gather*}
where the differentials can be identified with the differentials in the complex of proposition 4.2.1, and 
$\mathcal{O}_{\Delta}$ has a resolution
\begin{gather*}
\begin{CD}
(\ast\ast) \quad  0 @>>> \bigoplus_{\stackrel{i+j=n-1}{i,j\ge 0}}\mathcal{O}(-i-1, -j-1) @>d_{n-1}>> \dots
\end{CD}\quad\quad\quad\quad \\\begin{CD}\dots  @>d_{h+1}>> \left(\bigoplus_{\stackrel{i+j=h}{i, j\ge 0}}
\mathcal{O}(-i-1, -j-1)\right)^{\oplus {n+1 \choose h+2} } @> d_h>> \dots 
\end{CD}\\
\begin{CD}
\dots @>d_1 >> \mathcal{O}(-1 ,-1)^{\oplus {n+1 \choose 2} } @>>> \mathcal{O}_{\mathbb{P}^n\times\mathbb{P}^n}
@>>> \mathcal{O}_{\Delta} @>>> 0
\end{CD}
\end{gather*}
which is an Eagon-Northcott complex.\\
The next theorem gives a complete description of the functor $Rp_{2\ast}( p_1^{\ast}(-) \otimes^L 
\mathcal{O}_{\mathfrak{X}_0}) : D^b(Coh \, \mathbb{P}^n) \to D^b(Coh\, \mathbb{P}^n)$ (recall that in $D^b(Coh\, 
\mathbb{P}^n)$ one has the strong complete exceptional sequence $(\mathcal{O}, \: \mathcal{O}(1),\dots , \: 
\mathcal{O}(n))$ ).

\begin{theorem}
Let $\{ pt\}=L_0\subset L_1\subset \dots \subset
L_n=\mathbb{P}^n$ be a full flag of projective linear subspaces of
$\mathbb{P}^n$ (the Schubert varieties in $\mathbb{P}^n$) and denote
by $L^j$ the Schubert variety opposite to $L_j$.\\
For $d\ge 0$ one has in $D^b(Coh\, \mathbb{P}^n)$
\begin{gather*}
Rp_{2\ast}( p_1^{\ast}(\mathcal{O}(d)) \otimes^L 
\mathcal{O}_{\mathfrak{X}_0})\simeq \bigoplus_{j=0}^n
\mathcal{O}_{L_j}\otimes H^0(L^j, \mathcal{O}(d))^{\vee}/ H^0(L^{j+1},
\mathcal{O}(d))^{\vee} \, .
\end{gather*}
In terms of the coordinates $x_0, \dots ,x_n$, $y_0,\dots , y_n$
introduced above: 
\begin{gather*}
Rp_{2\ast}( p_1^{\ast}(\mathcal{O}(d)) \otimes^L 
\mathcal{O}_{\mathfrak{X}_0})\simeq \mathcal{O}\oplus (\mathcal{O}/(y_n))^{\oplus d} \oplus (\mathcal{O}/(y_n ,
y_{n-1}))^{\oplus \frac{d(d+1)}{2}}\oplus \dots\\
\dots \oplus (\mathcal{O}/(y_n, \dots , y_{n-i}) )^{\oplus{d+i \choose d-1}} \oplus \dots \oplus 
(\mathcal{O}/(y_n, \dots , y_{1}) )^{\oplus{d+n-1 \choose d-1}} .
\end{gather*}
Moreover for the map $\mathcal{O}(e) \stackrel{\cdot x_k}{\longrightarrow} \mathcal{O}(e+1)$ ($e\ge 0$, $0\le k 
\le  n$) one can describe the induced map $Rp_{2\ast}( p_1^{\ast}( \cdot x_k) \otimes^L 
\mathcal{O}_{\mathfrak{X}_0})$ as follows: \\
For each $d\ge 0$ and each $i=-1 ,\ldots , n-1$ choose a bijection between the set of monomials $M_i^d$ 
in the variables $x_{n-1-i}, \: x_{n-i} ,\dots , x_n$ of the form $x_{n-i-1}^{\alpha_1}x_{n-i}^{\alpha_2}
 \dots x_n^{\alpha_{i+2}}$ with $\alpha_1 > 0$ and $\sum \alpha_j=d$, and the set of copies of $\mathcal{O}
/(y_n,\dots ,y_{n-i})$ occuring in the above expression for $Rp_{2\ast}( p_1^{\ast}(\mathcal{O}(d))\otimes^L 
\mathcal{O}_{\mathfrak{X}_0})$. Then the copy of $\mathcal{O}
/(y_n,\dots ,y_{n-i})$ corresponding to a monomial $m\in M_i^e$ is mapped under 
 $Rp_{2\ast}( p_1^{\ast}( \cdot x_k) 
\otimes^L \mathcal{O}_{\mathfrak{X}_0})$ identically  to the copy of $\mathcal{O}/(y_n, \dots , y_{n-i})$
 corresponding to the monomial $x_km$ iff $x_k$ occurs in $m$. If $x_k$ does not occur in $m$ then the copy
 of $\mathcal{O}
/(y_n,\dots ,y_{n-i})$ corresponding to the monomial $m\in M_i^e$ is mapped to the copy of $\mathcal{O}
/(y_n,\dots ,y_{k+1})$ corresponding to $x_km$ via the natural surjection
\[
\mathcal{O}/(y_n,\dots ,y_{n-i}) \to \mathcal{O}/(y_n,\dots ,y_{k+1})  .
\]
\end{theorem}

\begin{proof}
For $d\ge 0$ one tensors the resolution ($\ast$) of $\mathcal{O}_{\mathfrak{X}_0}$ by $p_1^{\ast} \mathcal{O}(d)$ 
and notes that then all the bundles occuring in the terms of ($\ast$)$\otimes p_1^{\ast} \mathcal{O}(d)$ are 
$p_{2 \ast}$-acyclic whence $Rp_{2\ast}( p_1^{\ast}(\mathcal{O}(d)) \otimes^L 
\mathcal{O}_{\mathfrak{X}_0})$ is (as a complex concentrated in degree $0$) the cokernel of the map 
\[
\Phi\,  :\, \left( H^0(\mathbb{P}^n , \mathcal{O}(d-1))\otimes \mathcal{O}(-1) \right)^{\oplus{n+1\choose 2}}
\to H^0(\mathbb{P}^n , \mathcal{O}(d))\otimes \mathcal{O} 
\]
which on the various summands $H^0(\mathbb{P}^n , \mathcal{O}(d-1))\otimes\mathcal{O}(-1)$ of the domain is given by the maps
\begin{gather*}
H^0(\mathbb{P}^n , \mathcal{O}(d-1))\otimes\mathcal{O}(-1) \to H^0(\mathbb{P}^n , \mathcal{O}(d))\otimes \mathcal{O} \\
\quad\quad  m\otimes \sigma \mapsto x_i m\otimes y_j \sigma 
\end{gather*}
for $0\le i< j\le n$. For $i=-1, \ldots , n-1$ let  $M_i^d$ be as above the set of monomials in 
$x_{n-1-i}, \: x_{n-i},
 \dots , x_n$ of the form $x_{n-i-1}^{\alpha_1}x_{n-i}^{\alpha_2} \dots x_n^{\alpha_{i+2}}$ with $\alpha_1 > 0$
 and $\sum \alpha_j=d$. Then we have the identification
\[
H^0(\mathbb{P}^n , \mathcal{O}(d))\otimes \mathcal{O} \simeq \bigoplus_{i=-1}^{n-1} \left( \bigoplus_{m\in M_i^d}
\mathcal{O} \right) ;
\]
For given $m\in M_i^d$ write cont$(m)$ for the subset of the variables $x_0,\dots , x_n$ that occur in $m$. Then 
the map $\Phi$ above is the direct sum of maps 
\[
\bigoplus_{x_i\in\mathrm{cont}(m)} \bigoplus_{n\ge j > i}\mathcal{O}(-1) \to \mathcal{O}  \quad \forall i=-1, 
\dots , n-1 \:\forall m\in M_i^d
\]
which on the summand $\mathcal{O}(-1)$ on the left side of the arrow corresponding to $x_{i_0}\in\mathrm{cont}(m)$ 
and $j_0 > i_0$ are multiplication by $y_{j_0}$. Since $M_i^d$ has cardinality 
\[
{d+i \choose d-1} = {d+i+1\choose d} - {d+i \choose d}
\]
one finds that the cokernel of $\Phi$ is indeed
\[
\mathcal{O}\oplus (\mathcal{O}/(y_n))^{\oplus d} \oplus 
\dots \oplus (\mathcal{O}/(y_n, \dots , y_{n-i}) )^{\oplus{d+i \choose d-1}} \oplus \dots \oplus 
(\mathcal{O}/(y_n, \dots , y_{1}) )^{\oplus{d+n-1 \choose d-1}}
\]
as claimed.\\
The second statement of the theorem is now clear because $Rp_{2\ast}( p_1^{\ast}( \cdot x_k) \otimes^L 
\mathcal{O}_{\mathfrak{X}_0})$ is induced by the map $H^0(\mathbb{P}^n , \mathcal{O}(e))\otimes \mathcal{O}
\to  H^0(\mathbb{P}^n , \mathcal{O}(e+1))\otimes \mathcal{O}$ which is multiplication by $x_k$. 
\end{proof}

\begin{remark}
It is possible to prove Beilinson's theorem on $\mathbb{P}^n$ using only knowledge of the resolution ($\ast$) 
of $\mathcal{O}_{\mathfrak{X}_0}$: Indeed by theorem 4.1.1 one knows a priori that one can lift the resolution 
($\ast$) of $\mathcal{O}_{\mathfrak{X}_0}$ to a resolution of $\mathcal{O}_{\Delta}$ of the form $(\ast\ast)$ by 
flatness (cf. e.g. \cite{Ar}, part I, rem. 3.1). Since the terms in the resolution ($\ast$) are direct sums 
of bundles $\mathcal{O}(-k, -l)$, $0\le k, l\le n$, we find by the
standard argument from \cite{Bei}(i.e., the decomposition
$\mathrm{id}\simeq R p_{2\ast}(p_1^{\ast} (-)\otimes^L
\mathcal{O}_{\Delta })$) that 
$D^b(Coh\, \mathbb{P}^n)$ is generated by $(\mathcal{O}(-n), \dots , \mathcal{O})$. 
\end{remark}

Finally it would be interesting to know if one could find a resolution of $\mathcal{O}_{\mathfrak{X}_0}$ on $X
\times X$ for any rational homogeneous $X=G/P$ along the same lines as in this subsection, i.e. by first 
finding a ``monomial description'' of $\mathfrak{X}_0$ inside $X\times X$ (e.g. using standard monomial theory, 
cf. \cite{BiLa}) and then using the method of cellular resolutions from \cite{B-S}. Thereafter it would be 
even more important to see if one could obtain valuable information about $D^b(Coh \, X)$ by lifting the 
resolution of $\mathcal{O}_{\mathfrak{X}_0}$ to one of $\mathcal{O}_{\Delta}$.

\end{document}